\documentclass[10pt]{article}
\usepackage{latexsym,amsfonts,amssymb,amsthm,upref,amscd}
\usepackage{times}
\newtheorem{Theorem}{Theorem}[section]
\newtheorem{definition}[Theorem]{Definition}
\newtheorem{lemma}[Theorem]{Lemma}

\newtheorem{remark}[Theorem]{Remark}
\newtheorem{Open Problem}[Theorem]{Open Problem}

\makeatletter \@addtoreset{equation}{section} \makeatother

\def\be{{}\begin{equation}}
\def\ee{{}\end{equation}}

\begin{document}

\title{\bf  Multi-bump Solutions for  a Strongly Indefinite Semilinear
Schr\"odinger Equation Without Symmetry or convexity Assumptions
 }

\author{Shaowei Chen \thanks{The author acknowledges the support of
NNSF of China (No. 10526041). }\\ \\
\small LMAM, School of Mathematical Sciences, Peking
University, \\
\small Beijing 100871, P. R. China\\
\small and\\
\small  Department of Mathematics, Fujian Normal
University, \\
\small Fuzhou 350007, P. R.  China\\
 \small E-mail: swchen6@yahoo.com.cn\\
}

\date{}
\maketitle

\begin{minipage}{13cm}
{\small {\bf Abstract:} In this paper, we study the following
semilinear Schr\"odinger equation with periodic coefficient:
$$-\triangle u +V(x)u=f(x,u),\ u\in H^{1}(\mathbb{R}^{N}).$$
The functional corresponding to this equation possesses strongly
indefinite structure. The nonlinear term $f(x,t)$ satisfies some
superlinear growth conditions and  need not be  odd or increasing
strictly in $t$. Using a new variational reduction method  and a
generalized Morse theory, we proved that this equation has
infinitely many geometrically different solutions. Furthermore, if
the solutions of this equation under some energy level are isolated,
then we can show that this equation has infinitely many $m-$bump
solutions for any positive integer $m\geq 2.$

\medskip {\bf Key words:}\ semilinear
Schr\"odinger equation, multi-bump solutions, critical group, reduction methods. \\
\medskip 2000 Mathematics Subject Classification:  35J20, 35J70}
\end{minipage}

\section{Introduction and main results}
\hspace*{\parindent}In this paper, we consider the following problem
\begin{equation}\label{1.1}
-\triangle u+V(x)u=f(x,u),\ u\in H^{1}(\mathbb{R}^{N}),
\end{equation}
where $x=(x_{1},x_{2},\cdots,x_{N})\in\mathbb{R}^{N} \ (N\geq 1)$
and $V$, $f$ satisfy the following conditions:
\begin{description}
\item{$(\bf{V_{1}}).$} $V\in L^{\infty}(\mathbb{R}^{N})$ is
$1-$periodic in each $x_{i},\ i=1,\cdots, N.$
\item{$(\bf{V_{2}}).$} The linear operator
$L:H^{2}(\mathbb{R}^{N})\rightarrow L^{2}(\mathbb{R}^{N}),\
u\mapsto -\triangle u+Vu$ is invertible and $0$ lies in a gap of
the spectrum of $L$. \item{$(\bf{f_{1}}).$} $f(x,t)$ is a
Caratheodory function and is $1-$periodic in each $x_{i},\
i=1,\cdots, N.$ $f'_{t}(x,t)$ exists for every $t\in\mathbb{R}$
and for almost all $x\in\mathbb{R}^{N}.$ And $f'_{t}(x,t)$ is a
Caratheodory function.

 \item{$(\bf{f_{2}}).$}
For some $2<q<p<2^{*}:=\left\{
\begin{array}{l}
\frac{2N}{N-2}, \ N\geq 3\\
\infty, \quad N=1,2\\
\end{array} \right.$ and $C>0,$
$$|f'_{t}(x,t)|\leq C(|t|^{q-2}+|t|^{p-2}), \ \mbox{for any}\ (x,t)\in\mathbb{R}^{N}\times \mathbb{R}.$$
\item{$(\bf{f_{3}}).$} There exists $\gamma>2$ such that for every
$t\neq 0$ and $x\in\mathbb{R}^{N}$,
$$0<\gamma F(x,t)\leq tf(x,t),$$
where $F(x,t)=\int^{t}_{0}f(x,s)ds.$
\end{description}

\medskip

Note that  the power nonlinearity  $f(x,t)=h(x)|u|^{p-2}u,$ with
positive $1-$periodic $h\in L^{\infty}(\mathbb{R}^{N}),$ $
h(x)\geq 0,\ h\not\equiv 0$ and $2<p<2^{*}$, satisfies all the
assumptions.

 Under the assumptions $\bf (V_{1})-(V_{2})$ and $\bf (f_{1})-(f_{3})$ the functional
\begin{equation}\label{fanon7htbvzp0o}
J(u)=\frac{1}{2}\int_{\mathbb{R}^{N}}(|\nabla
u|^{2}+V(x)u^{2})dx-\int_{\mathbb{R}^{N}}F(x,u)dx
\end{equation}
is of class $C^{2}$ on the Sobolev space $H^{1}(\mathbb{R}^{N})$
and critical points of (\ref{fanon7htbvzp0o}) correspond to weak
solutions of Equation (\ref{1.1}).

The operator $L=-\triangle +V$ (on $L^{2}(\mathbb{R}^{N})$) has
purely continuous spectrum which is bounded below and consists of
closed disjoint intervals (\cite[Theorem XIII.100]{RS}). We denote
by $|L|^{1/2}$ the square root of the absolute value of $L.$ The
domain of $|L|^{1/2} $ is the space $X:=H^{1}(\mathbb{R}^{N})$. On
$X,$ we choose the inner product
$(u,v)_{X}=\int_{\mathbb{R}^{N}}|L|^{1/2}u\cdot |L|^{1/2}v dx$ and
the corresponding norm $||u||=\sqrt{ (u,u)_{X}}.$ There exists an
orthogonal decomposition $X=Y\oplus Z$ such that $Z$ and $Y$ are the
positive and negative spaces corresponding to the spectral
decomposing of $L$.  They are invariant under the action of
$\mathbb{Z}^{N},$ i.e., for any $u\in Y$ or $u\in Z$ and for any
$\mathbf{k}=(n_{1},\cdots,n_{N})\in \mathbb{Z}^{N},$
$u(\cdot-\mathbf{k})$ is also in $Y$ or $Z.$ Furthermore,
\begin{equation}\label{hnbxxc5r8871z} \forall u\in Y,\
\int_{\mathbb{R}^{N}}(|\nabla
u|^{2}+Vu^{2})dx=-(u,u)_{X}=-||u||^{2},\end{equation}
\begin{equation}\label{un77715490kjh76}
\forall u\in Z,\ \int_{\mathbb{R}^{N}}(|\nabla
u|^{2}+Vu^{2})dx=(u,u)_{X}=||u||^{2}.
\end{equation}
Since  $0$ lies in a gap of spectrum of $L$,  the dimensions of $Y$
and $Z$ are both infinity. In this case, Equation (\ref{1.1}) is
called strongly indefinite.
 By (\ref{hnbxxc5r8871z}) and (\ref{un77715490kjh76}), we have
\begin{equation}\label{1.2}
J(u)={1\over{2}}||Qu||^{2}-{1\over{2}}||Pu||^{2}-\int_{\mathbb{R}^{N}}F(x,u)dx,\
u\in X.
\end{equation}

It is easy to verify that  if $v$ is a solution of Equation
(\ref{1.1}), then $v(\cdot-\mathbf{k})$ is also a solution of
Equation (\ref{1.1}) for any $\mathbf{k}\in\mathbb{Z}^{N}.$ Let $u$
and $v$ be two solutions of Equation (\ref{1.1}). They are called
geometrically different if $u(\cdot-\mathbf{k})\neq v$ for any
$\mathbf{k}\in\mathbb{Z}^{N}.$ Let $v_{1},\cdots,v_{n}$ be solutions
of Equation (\ref{1.1}) such that their barycenter are sufficiently
separated. Solutions of Equation (\ref{1.1}) that are close to
$\sum^{m}_{i=1}v_{i}$ are called $m-$bump solutions. The main result
of this paper is the following theorem:

\begin{Theorem}\label{theorem1}
Assume $\bf {(V_{1})-(V_{2})}$ and $\bf{(f_{1})-(f_{3})}.$ Equation
(\ref{1.1}) has infinitely many, geometrically different solutions.
Furthermore, if the condition {\bf{(*)}} (see Section 2 for its
definition) holds, then for any positive integer $m\geq 2,$ Equation
$(\ref{1.1})$ has infinitely many, geometrically different, $m-$bump
solutions.
\end{Theorem}

Equation (\ref{1.1}) arises from studying of steady state and
standing wave solutions of time-independent nonlinear
Schr\"odinger equations. Readers can consult \cite{Pan} for more
physical background and applications of Equation (\ref{1.1}).
Semilinear Schr\"odinger equation with periodic potential has been
studied by many authors in the past decade. In the celebrate
papers \cite{Rabi} and \cite{Rabi2}, Coti Zelati and Rabinowitz
used a variational gluing methods to obtain multi-bump type
solutions for Hamilton ODE and elliptic PDEs with periodic
potential. The linear parts of the Hamilton ODE and elliptic PDEs
they studied are positive definite and the functionals
corresponding to these Hamilton ODE and elliptic PDEs have
Mountain Pass structures. Coti Zelati and Rabinowitz used the
solutions obtained by Mountain Pass theorem as basic building
blocks to construct multi-bump solutions. Readers can consult
\cite{A2}, \cite{LW}, \cite{LW1} and references therein for more
recent development in this direction. In \cite{S}, S\'er\'e
considered some Hamiltonian systems whose linear parts are
strongly indefinite, i.e., the dimensions of the positive and
negative spaces corresponding to the spectral decomposing are both
infinity. He constructed multi-bump solutions for these
Hamiltonian systems.  But he imposed some convexity conditions on
the nonlinear terms of these Hamiltonian systems and then
transform them by dual variational methods into some equivalent
systems whose variational functionals are bounded below. For the
strong indefinite semilinear Schr\"odinger equation (\ref{1.1}),
Alama and Li constructed multi-bump solutions in \cite{AL} by dual
variational methods under the assumption that $f(x,t)$ increases
strictly in $t,$ i.e., the function $F(x,t)$ is  convex. The first
work of directly dealing with Equation (\ref{1.1}) without the
convexity assumptions on nonlinear term $f$ was done by C.
Troestler and M. Willem in \cite{T}. They obtained the results on
the existence of nontrivial solutions to Equation (\ref{1.1}).
 In \cite{KS},  Kryszewski and  Szulkin
 obtained the result that there exist infinitely many geometrically different  solutions to Equation
 (\ref{1.1}) whenever $f$ is odd in $u.$ In a very recent paper
 \cite{A}, Ackermann provided an interesting abstract framework in which
 multi-bump solutions can be obtained in many situations. It
 reduces the problem of constructing multi-bump solutions to the
 problem of finding an isolated solution with nontrivial topology
 in a specific sense. Using the abstract results in \cite{A}, a
 very general result on the existence of multi-bump solutions to
 strong indefinite periodic semilinear Schr\"odinger equations (even with nonlocal nonlinearities)
 is obtained. However in \cite{A}, the result on the existence of multi-bump solutions for   Equation
 (\ref{1.1}) was obtained under the assumptions that the nonlinear term $f(x,t)$
 is
 $C^{2}$ and convex in $t$. Therefore, the question on
 the existence of infinitely many geometrically different solutions for Equation
 (\ref{1.1}) without the assumption that $f$ is odd  in $t$ or
 $F(x,t)$ is  convex in $t$
 was still left open.

 Theorem \ref{theorem1} of the present paper gives an affirmative answer to this open problem.
  In the present paper we shall show that Equation (\ref{1.1}) has a
 solution $u_{0}\neq 0$ which has the properties that after
 reducing the corresponding functional of Equation (\ref{1.1}) in a
 neighborhood of $u_{0},$ the critical group of the reduction
 function in the critical point $u_{0}$ is nontrivial. Then using
 $u_{0}$ as a basic building block, we constructed multi-bump
 solutions for Equation (\ref{1.1}) by a perturbation technique
stemed from
 Chang and Ghoussoub (see \cite{CG}). To obtain such $u_{0}$, we
 consider  the approximation problem firstly:
 \begin{equation}\label{huannuzhixueyuan6t}
 -\triangle u+V(x)u=f(x,u),\ u\in H^{1}_{per}(Q_{k}),
 \end{equation}
where $Q_{k}$ is cube  of $\mathbb{R}^{N}$ with edge length
$k\in\mathbb{N}$ and $H^{1}_{per}(Q_{k})$ denotes the space of
$H^{1}(Q_{k})-$functions which are $k-$periodic in $x_{i}$,
$i=1,2,\cdots,N.$ The variational functional corresponding to
(\ref{huannuzhixueyuan6t}) satisfies Palais-Smale condition and has
linking structure (see \cite{P}). Secondly, using the linking
theorem (one can see \cite{Will} for reference), we can get a
solution $u_{k}$ of (\ref{huannuzhixueyuan6t}) which satisfies that
there exist finite many nontrivial solutions $u^{i},\ i=1,\cdots,n$
of Equation (\ref{1.1}) and  sequences $\{b^{i}_{k}\},\
i=1,2,\cdots,n$ such that $|b^{i}_{k}-b^{j}_{k}|\rightarrow\infty,\
i\neq j$ as $k\rightarrow\infty$ and
$$||u_{k}-\sum^{n}_{i=1}u^{i}(\cdot-b^{i}_{k})||_{H^{1}(Q_{k})}\rightarrow 0.$$
Finally, we show that at least one of $u^{i},\ i=1,\cdots,n$ has the
properties that its critical group of the reduction
 function  is nontrivial.

 This paper is organized as follows:  From Section 2  to Section 4, we use
 an
 approximation method, reduction methods and critical point theory to
 obtain the existence of a special nontrivial solution of Equation of
 (\ref{1.1}) which has the properties we mentioned above. In
 Section 5, we give the proof of Theorem \ref{theorem1}. In Section
 6, we provide the detail proofs of some Lemmas stated in Section
 3.

\medskip

\noindent{\bf Notation.} $\mathbb{R},\ \mathbb{Z}$ and $\mathbb{N}$
denote the sets of real number, integer and positive integer
respectively. $B_{E}(a,\rho)$ denotes the open ball in $E$ centered
at $a$ and having radius $\rho$. The closure of a set $A$ is denoted
by $\overline{A}$ or $cl(A).$ By $\rightarrow$ we denote the strong
and by $\rightharpoonup$ the weak convergence. $\mbox{dist}(a, A)$
denotes the distance from  the point $a$ to the set $A$.
$\mbox{diam}(A)$ denotes the diameter of the set $A.$  By $\ker A$
 denotes the null space of the operator $A.$ If $f$ is a $C^{2}$ functional defined on a Hilbert
space $H$, $\nabla f$ (or $df$) and $\nabla^{2}f$ denote the
gradient of $f$ and the second differential of $f$ respectively. And
for $a, b\in\mathbb{R},$ we denote $f^{a}:=\{u\in H\ :\ f(u)\leq
a\}$ and $f_{b}:=\{u\in H\ :\ f(u)\geq b\}$ the sub- and superlevel
sets of the functional $f, $ moreover, $f^{a}_{b}:=\{u\in H\ :\
b\leq f(u)\leq a\}.$ $\delta_{i,j}$ denotes the Kronecker notation:
$ \delta_{i,j}=\left\{
\begin{array}{l}
1,  i=j\\
0,  i\neq j.\\
\end{array} \right.$ If $H$ is a Hilbert space and $W$ is a closed
subspace of $H,$ we denote the orthogonal complement space of $W$ in
$H$ by $W^{\bot}.$ For a subset $A\subset H,$ $\mbox{span}\{A\}$
denotes the subspace of $H$  generated by $A.$

\section{A periodic  approximation problem}

\hspace*{\parindent}Associated with Equation (\ref{1.1}), we study
the approximation problem in cubes $Q_{k}$ of $\mathbb{R}^{N}$
with edge length $k\in\mathbb{N}$
\begin{equation}\label{jjb}
-\triangle u+V(x)u=f(x,u)\ \mbox{in}\ Q_{k}, \ u\in
E_{k}:=H^{1}_{per}(Q_{k}),
\end{equation}
where $H^{1}_{per}(Q_{k})$ denotes the space of
$H^{1}(Q_{k})-$functions which are $k-$periodic in $x_{i}$,
$i=1,2,\cdots,N.$ The operator $-\triangle+V$ on
$L^{2}_{per}(Q_{k})$ has discrete spectrum with eigenvalues
$\lambda_{k,1}\leq\lambda_{k,2}\leq\cdots\rightarrow +\infty$ and
there exists a finite $\min\{i :  \lambda_{k,i}>0\}.$  Moreover,
every eigenvalue $\lambda_{k,i}$ is contained in the spectrum of
$-\triangle+V$ on the whole space. This follows from the spectral
gap around $0$ assumed in Reed and Simon (see \cite{RS}). Therefore,
if $(-\alpha,\beta),\ \alpha,
 \beta>0$ denotes the
spectral gap around $0$ assumed in $\bf (V_{2})$. We claim that
$\lambda_{k,i}\not\in (-\alpha,\beta)$ for every $k,i\in\mathbb{N}.$
We denote by $\phi_{k,i}$ the corresponding eigenfunctions.

Let $j(k)=\min\{i :  \lambda_{k,i}>0\}-1.$ Now we define an
orthogonal decomposition of $E_{k}$ by $E_{k}=Y_{k}\oplus Z_{k},$
where $$Y_{k}=\mbox{span} \{\phi_{k,1},\cdots, \phi_{k,j(k)}\},\
Z_{k}=Y^{\bot}_{k}.$$ The associated energy functional to
(\ref{jjb}) is
\begin{equation}\label{ncbcvddfdf}
J_{k}(u)=\frac{1}{2}\int_{Q_{k}}(|\nabla
u|^{2}+V(x)u^{2})dx-\int_{Q_{k}}F(x,u)dx.
\end{equation}
We may define a new inner product $(\cdot,\cdot)_{k}$ on $E_{k}$
with corresponding norm $||\cdot||_{k}$ such that
\begin{equation}\label{chuenju}
\forall u\in Y_{k},\ \int_{Q_{k}}(|\nabla
u|^{2}+Vu^{2})dx=-(u,u)_{k}=-||u||_{k}^{2},
\end{equation}
\begin{equation}\label{jnadsewe}
\forall u\in Z_{k},\ \int_{Q_{k}}(|\nabla
u|^{2}+Vu^{2})dx=(u,u)_{k}=||u||_{k}^{2}.
\end{equation}
If we denote by $P_{k}:E_{k}\rightarrow Y_{k}$ and
$T_{k}:E_{k}\rightarrow Z_{k}$ the orthogonal projections, our
functional becomes
\begin{equation}\label{vvvvvl}
J_{k}(u)=\frac{1}{2}(||T_{k}u||_{k}^{2}-||P_{k}u||_{k}^{2})-\int_{Q_{k}}F(x,u)dx.
\end{equation}

For convenience, we assume that
$Q_{k}=(-\frac{k}{2},\frac{k}{2})^{N},\ k\in\mathbb{N},$ then
\begin{eqnarray}\label{mabsdefrt}
&&H^{1}_{per}(Q_{k}) \nonumber\\
&=&\{u\in H^{1}(Q_{k}) : u(x_{1},\cdots,x_{i-1},
-\frac{k}{2},x_{i+1},\cdots,x_{N})=u(x_{1},\cdots,
x_{i-1},\frac{k}{2},x_{i+1},\cdots,x_{N}),\nonumber\\
&& (x_{1},\cdots,x_{i-1},x_{i+1},\cdots,x_{N})\in
(-\frac{k}{2},\frac{k}{2})^{N-1},\ i=1,\cdots,N\}.\nonumber
\end{eqnarray}
Let $\mathbb{Z}_{k}=\mathbb{Z}/k\mathbb{Z}$. For
$\mathbf{b}\in\mathbb{Z}^{N}_{k}$ and $u\in E_{k}$, the action of
$\mathbf{b}$ on $u$, we still denote it by $u(\cdot+\mathbf{b})$, is
defined by the following way: For $\mathbf{b}=(0,\cdots,0,
\stackrel{i}{1}, 0\cdots,0)\in \mathbb{ Z}_{k}^{N}$ and $u\in
E_{k},$
$$ u(x+\mathbf{b})=\left\{
\begin{array}{l}
u(x_{1},\cdots,x_{i-1},x_{i}+1,x_{i+1},\cdots,x_{N}),  \quad\quad-\frac{k}{2}\leq x_{i}\leq\frac{k}{2}-1\\
u(x_{1},\cdots,x_{i-1},x_{i}+1-k,x_{i+1},\cdots,x_{N}),  \quad\frac{k}{2}-1\leq x_{i}\leq\frac{k}{2},\\
\end{array} \right.$$
$$ u(x-\mathbf{b})=\left\{
\begin{array}{l}
u(x_{1},\cdots,x_{i-1},x_{i}-1,x_{i+1},\cdots,x_{N}),  \quad\quad\ \ -\frac{k}{2}+1\leq x_{i}\leq\frac{k}{2}\\
u(x_{1},\cdots,x_{i-1},x_{i}-1+k,x_{i+1},\cdots,x_{N}),  \quad-\frac{k}{2}\leq x_{i}\leq-\frac{k}{2}+1.\\
\end{array} \right.$$
Since $V(x)$ and $f(x,t)$ is $1-$periodic in $x_{k},\ k=1,\cdots,
N,$ we deduce that $J_{k}$ is invariant under the action of
$\mathbb{Z}_{k}^{N}.$

\begin{lemma}(\cite[Lemma 2]{P})\label{whoami}
There exist constants $C_{1}>0$ and $C_{2}>0$ which are independent
of $k$ such that for any $u\in E_{k},$
$$C_{1}||u||_{H^{1}(Q_{k})}\leq ||u||_{k}\leq C_{2}||u||_{H^{1}(Q_{k})}.$$
\end{lemma}

By the conditions ${\bf (f_{1})}-{\bf (f_{3})}$, we have the
following Lemma (one can see \cite{Rabin} for reference):

\begin{lemma}\label{dalianmao}
For any $k\in \mathbb{N}$, $J_{k}$ satisfies $(PS)$ conditions.
\end{lemma}

\begin{lemma}\label{dadaoruqintain5tqddu} (\cite[Lemma 3.3]{Pan} or \cite[Lemma 4]{P})
There exists $\epsilon_{1}>0$ not depending on $k$ such that
$||u_{k}||_{k}\geq\epsilon_{1},\ ||u||\geq \epsilon_{1}$ holds for
any nontrivial critical point $u_{k}$ of $J_{k}$ and $u$ of $J.$
In addition, there exists $\epsilon_{2}>0$ not depending on $k$
such that $J_{k}(u_{k})\geq\epsilon_{2},\ J(u)\geq \epsilon_{2}$
holds for any nontrivial critical point $u_{k}$ of $J_{k}$ and $u$
of $J.$
\end{lemma}

\begin{lemma}(\cite[Lemma 8]{P})\label{dahuoe}
There exist real number $\delta>0$ and $ r>0$ which are independent
of $k$ such that $\displaystyle\inf_{u\in N_{k}}J_{k}(u)\geq\delta$,
where $N_{k}=\{z\in Z_{k}  : ||z||_{k}=r\}.$
\end{lemma}

Now for each $k$, we fix a function $z_{k}\in Z_{k}$ with
$||z_{k}||_{k}=1.$ For $\rho>0,$ we define the sets
$$M_{k}=\{y+tz_{k} : ||y+tz_{k}||_{k}\leq\rho,\ t\geq 0,\ y\in Y_{k}\}.$$

\begin{lemma}(\cite[Lemma 9]{P})\label{cactv}
There exists a $\rho>r$ which is independent of $k$ such that
$$\displaystyle\sup_{u\in\partial M_{k}}J_{k}(u)=0.$$
\end{lemma}

\begin{lemma}(\cite[Theorem 3.4]{Pan} or \cite[Lemma 10]{P})\label{cctv}
The number $$c_{k}=\inf_{h\in \Gamma_{k}}\sup_{u\in
M_{k}}J(h(u))$$ is a critical value of $J_{k}$ and there exists
positive number $M$ which is independent of $k$ such that
$0<\delta\leq c_{k}\leq M<\infty$, where
$$\Gamma_{k}=\{h\in C(E_{k},E_{k}) : h|_{\partial M_{k}}=id\}.$$
\end{lemma}

Let $\chi_{k}$ be cut-off functions such that $0\leq\chi_{k}\leq 1,\
\chi_{k}\equiv 1$ on $Q_{k-1}, \ \chi_{k}\equiv0$ outside of $Q_{k}$
and $|\nabla \chi_{k}|\leq C,\ k=1,2,\cdots.$

\medskip

\begin{lemma}(\cite[Theorem 5.1]{Pan} or \cite[Theorem 11]{P})\label{cotgaadsewe}
Under the assumptions $\bf(f_{1})-\bf{(f_{3})}$. Let $v_{k}\in
E_{k}$ be a uniformly bounded sequence which satisfies
$J'_{k}(v_{k})\rightarrow 0$ and
$\widetilde{c}_{k}=J_{k}(v_{k})\rightarrow \widetilde{c}>0.$ Then
there exist critical points $v^{i}$, $i=1,2,\cdots,\nu$ of $J$ and
sequences $d^{i}_{k}\in\mathbb{Z}^{N}$ such that as
$k\rightarrow\infty,$ $|d^{i}_{k}-d^{j}_{k}|\rightarrow\infty,\
i\neq j,$
$$||v_{k}-\sum^{\nu}_{i=1}v^{i}(\cdot+d^{i}_{k})||_{k}\rightarrow 0$$
and $\displaystyle\sum^{\nu}_{i=1}J(v^{i})=\widetilde{c}.$
\end{lemma}

 Let $K$ and $K_{k}$ be the sets of critical points
of $J$ and $J_{k},\ k=1,2,\cdots$ respectively. For
$a,b\in\mathbb{R},$ let $K^{a}:= K\cap J^{a},\ K_{a}:=K\cap
J_{a}$, $ K^{b}_{a}:=K\cap J^{a}\cap J_{b}$ and
$K_{k}^{a}=J^{a}_{k}\cap K_{k}.$

Let
\begin{equation}\label{yhgtff45}
\displaystyle c_{0}=\sup_{k}c_{k},
\end{equation} where $c_{k}$ is
the minimax value defined in Lemma \ref{cctv}.    Now we impose
the following condition on Equation (\ref{1.1}):

\bigskip

\noindent{\large\bf (*).  There exists $\alpha_{0}>0$ such that
$K^{c_{0}+\alpha_{0}}/ \mathbb{Z}^{N}$ is finite.}

\bigskip

\begin{lemma}\label{shanshan}
If the condition {\bf(*)} holds, then the following three statements
hold:
\begin{description}
\item{(1).} There exists $\beta_{0}\in(0,\alpha_{0})$
 such that \begin{equation}\label{wumenpanbozouniw}
 \inf\{||\nabla J(u)|| : u\in X,\ J(u)=c_{0}+\beta_{0}\}>0
 \end{equation}
 and there exists constant $\epsilon_{3}>0$ not depending on $k$
 such that
 \begin{equation}\label{saniaodonglixia8j}
 \inf\{||\nabla J_{k}( u)||_{k} : u\in E_{k},\ J_{k}(u)=c_{0}+\beta_{0}\}>\epsilon_{3}
 \end{equation}
 if $k$ large enough.

\item{(2).} There exist $\delta_{0}>0$ and $k_{0}\in\mathbb{N}$ such
that if $k\geq k_{0},$ then for any $u_{k}\in
K^{c_{0}+\beta_{0}}_{k},$
$$ B_{E_{k}}(u_{k},\delta_{0})\cap
K^{c_{0}+\beta_{0}}_{k}=B_{E_{k}}(u_{k},\delta_{0}/2)\cap
K^{c_{0}+\beta_{0}}_{k}.$$ It follows that when $k\geq k_{0},$  for
any $u^{1}_{k}, u^{2}_{k}\in K^{c_{0}+\beta_{0}}_{k},$ either
$||u^{1}_{k}-u^{2}_{k}||_{k}\leq\delta_{0}/2$ or
$||u^{1}_{k}-u^{2}_{k}||_{k}\geq\delta_{0}.$

\item{(3).} For any $\epsilon>0,$ there exists
$k_{\epsilon}\in\mathbb{N}$ such that if $k\geq k_{\epsilon}$, then
for any $u_{k}\in K^{c_{0}+\beta_{0}}_{k},$
$$B_{E_{k}}(u_{k},\delta_{0})\cap
K^{c_{0}+\beta_{0}}_{k}=B_{E_{k}}(u_{k},\epsilon)\cap
K^{c_{0}+\beta_{0}}_{k},$$ where $\delta_{0}$ is the constant
appeared in $(2).$ It follows that when $k\geq k_{\epsilon},$ for
any $u^{1}_{k}, u^{2}_{k}\in K^{c_{0}+\beta_{0}}_{k},$ either
$||u^{1}_{k}-u^{2}_{k}||_{k}\leq\epsilon$ or
$||u^{1}_{k}-u^{2}_{k}||_{k}\geq\delta_{0}.$
\end{description}
\end{lemma}
\noindent{\bf Proof.}
 By Lemma \ref{dadaoruqintain5tqddu}, we know
that there exists $\epsilon_{2}>0$ such that for any $u\in K,$
$J(u)\geq\epsilon_{2}$. Let
$\overline{l}=[\frac{\alpha_{0}+c_{0}}{\epsilon_{2}}]+1$ and let
$$\mathcal{F}_{\overline{l}}(K^{c_{0}+\alpha_{0}}):=
\{\sum^{j}_{i=1}v_{i}(\cdot-b_{i}) : 1\leq i\leq j,\ \ 1\leq j\leq
\overline{l},\ v_{i}\in K^{c_{0}+\alpha_{0}},\
b_{i}\in\mathbb{Z}^{N}\}.$$ If $c\in (c_{0},c_{0}+\alpha_{0})$
satisfies that there exists a sequence $\{u_{m}\}$ such that as
$m\rightarrow\infty,$
\begin{equation}\label{huonanxiong7hppplm}
J(v_{m})\rightarrow c,\  ||\nabla J(v_{m})||\rightarrow 0,
\end{equation}
then by Proposition 1.24 of \cite{Rabi}, we deduce that there exist
at most $\overline{l}$ nontrivial solutions $v_{i}\in
K^{c_{0}+\alpha_{0}},\ i=1,\cdots,\overline{l}$ and $\overline{l}$
sequence $\{d^{i}_{m}\}\subset\mathbb{Z}^{N}$ such that as
$m\rightarrow\infty,$ $|d^{i}_{m}-d^{j}_{m}|\rightarrow\infty,\
i\neq j$, $c=\sum^{\overline{l}}_{i=1}J(v^{i})$ and
\begin{equation}\label{nb8hahags55dfe7}
||u_{m}-\sum^{\overline{l}}_{i=1}v^{i}(\cdot-d^{i}_{m})||\rightarrow
0.
\end{equation}
By the condition {\bf (*)}, we know that
$$\mathcal{A}=\{\sum^{j}_{i=1}J(u_{i})\ :\ 1\leq i\leq j,\ 1\leq j\leq \overline{l},\
u_{i}\in K^{c_{0}+\alpha_{0}}\}$$
 is a finite set. It follows that the possible $c\in (c_{0},c_{0}+\alpha_{0})$
which satisfies (\ref{huonanxiong7hppplm}) is finite. If we choose
$\beta_{0}\in (0,\alpha_{0})$ such that $c_{0}+\beta_{0}\in
(c_{0},c_{0}+\alpha_{0})\setminus\mathcal{A}$, then
(\ref{wumenpanbozouniw}) holds.

We are ready to prove that (\ref{saniaodonglixia8j}) holds. If not,
then there exists a sequence $\{u_{k}\}$ such that $u_{k}\in E_{k},\
J_{k}(u_{k})=c_{0}+\beta_{0},\ k=1,2,\cdots$ and $||\nabla
J_{k}(u_{k})||_{k}\rightarrow 0$ as $k\rightarrow\infty.$ Then by
Lemma \ref{cotgaadsewe} we deduce that there exist at most
$\overline{l}$ nontrivial solutions $v_{i}\in K^{c_{0}+\alpha_{0}},\
i=1,\cdots,\overline{l}$ and $\overline{l}$ sequence
$\{d^{i}_{m}\}\subset\mathbb{Z}^{N}$ such that as
$m\rightarrow\infty,$ $|d^{i}_{m}-d^{j}_{m}|\rightarrow\infty,\
i\neq j$, $c_{0}+\beta_{0}=\sum^{\overline{l}}_{i=1}J_{k}(v^{i})$
and
$$||u_{m}-\sum^{\overline{l}}_{i=1}v^{i}(\cdot-d^{i}_{m})||_{k}\rightarrow 0.$$
It follows that $c_{0}+\beta_{0}\in \mathcal{A}.$ It is a
contradiction. Thus (\ref{saniaodonglixia8j}) holds.

 By the condition {\bf{(*)}} and Proposition 1.55 of
\cite{Rabi}, we know that
\begin{equation}\label{contizelatiandrobi}
\mu=\mu(\mathcal{F}_{\overline{l}}(K^{c_{0}+\alpha_{0}})):=\inf\{||x-y||\
:\ x\neq y \in \mathcal{F}_{\overline{l}}(K^{c_{0}+\alpha_{0}})\}>0.
\end{equation}
Choose $\delta_{0}=\mu/2.$ If there exist  two sequences
$\{u^{1}_{k}\}$ and $\{u^{2}_{k}\}$ such that $u^{i}_{k}\in
K^{c_{0}+\beta_{0}}_{k},\ k=1,2,\cdots,\ i=1,2$ and
$\frac{\delta_{0}}{2}<||u^{1}_{k}-u^{2}_{k}||<\delta_{0}$, then by
Lemma \ref{cotgaadsewe}, we deduce that there exist $v^{1}_{k},
v^{2}_{k}\in \mathcal{F}_{\overline{l}}(K^{c_{0}+\alpha_{0}})$ such
that as $k\rightarrow\infty$
$$||u^{i}_{k}-v^{i}_{k}||_{k}\rightarrow 0,\ i=1,2.$$
It follows that $||v^{1}_{k}-v^{2}_{k}||\leq\delta_{0}<\mu$ when $k$
large enough. It contradicts the definition of $\mu$. Thus the
result of $(2)$ holds.

The proof of  result $(3)$ is similar. \hfill$\Box$

\begin{remark}\label{jaijianjianhian6tga5ffa5fag}In fact,
by Lemma \ref{cotgaadsewe}, we  get that as $k\rightarrow\infty,$
$$\mbox{dist}(K^{c_{0}+\beta_{0}}_{k},\mathcal{F}_{\overline{l}}(K^{c_{0}+\alpha_{0}}))\rightarrow 0.$$
By (\ref{contizelatiandrobi}), we know that
$\mathcal{F}_{\overline{l}}(K^{c_{0}+\alpha_{0}})$ is a discrete
set. Thus $K^{c_{0}+\beta_{0}}_{k}$ can be decomposed into a union
of its subsets which are disjoint each other.
\end{remark}

\medskip

By Lemma \ref{shanshan}, we have the following Lemma
\begin{lemma}\label{llkaddas}
Suppose that the condition {\bf{(*)}} holds. There exist positive
integer $m_{k}>0$ and $m_{k}$ subsets $K^{(i)}_{k},\ i=1,2\cdots,
m_{k}$ of $K^{c_{0}+\beta_{0}}_{k}$ such that if $k\geq k_{0}$, then
$$\displaystyle K_{k}^{c_{0}+\beta_{0}}=\displaystyle\bigcup^{m_{k}}_{i=1}K^{(i)}_{k},$$
$\mbox{dist}(K^{(i)}_{k},K^{(j)}_{k})\geq\delta_{0},\ i\neq j$ and $
\mbox{diam}(K^{(i)}_{k})\leq\delta_{0}/2,\ i=1,2\cdots,m_{k},$ where
$\delta_{0},$ $k_{0}$ and $\beta_{0}$ are the constants appeared in
Lemma \ref{shanshan}. Furthermore, for any $1\leq i\leq m_{k}$, as
$k\rightarrow\infty, $
$$\mbox{diam}(K^{(i)}_{k})\rightarrow 0.$$
\end{lemma}
\noindent{\bf Proof.} By Lemma \ref{shanshan}, we know that if
$k\geq k_{0},$ then $K^{c_{0}+\beta_{0}}_{k}$ can be decomposed into
 a union of its subsets $K^{(i)}_{k}$ which are disjoint each other.
 We show that the number of these subsets is finite. If not, by Lemma \ref{shanshan} there
 exist $K^{(i)}_{k},\ i=1,2,\cdots$ such that
 $K^{c_{0}+\beta_{0}}_{k}=\displaystyle\bigcup^{\infty}_{i=1}K^{(i)}_{k}$
 and $\mbox{dist}(K^{(i)}_{k},k^{(j)}_{k})\geq\delta_{0},\ i\neq j.$
 Choose $u_{i,k}\in K^{(i)}_{k}$. Then $u_{i,k}$ satisfies that
 \begin{equation}\label{635tgd8jnyc90008j}
 ||u_{i,k}-u_{j,k}||\geq\delta_{0},\ i\neq j.\end{equation}
 By Lemma \ref{dalianmao}, we get that there exists  a subsequence $\{u_{i_{m},k}\}$ of
 $\{u_{i,k}\}$ and $u\in E_{k}$ such that $||u_{i_{m},k}-u||_{k}\rightarrow
 0$ as $m\rightarrow\infty.$ It contradicts
 (\ref{635tgd8jnyc90008j}). Therefore, the number  of $K^{(i)}_{k}$
 is finite. We denote it by $m_{k}.$ Finally, by the result $(3)$ of Lemma
 \ref{shanshan}, we get that $\mbox{diam}(K^{(i)}_{k})\rightarrow 0$
 as $k\rightarrow\infty.$\hfill $\Box$

\medskip

Let $H_{*}(A,B)$ be the $*-$th singular homology group with
coefficient $\mathbb{Z}_{2}$. By the definition of
$c_{k}=\inf_{h\in\Gamma_{k}}\max_{u\in M_{k}}J_{k}(h(u))$, the
Linking Theorem (one can see \cite{Rabin} or \cite{Will} for
reference) and the proof of Theorem 7.5 of \cite{B}, we have the
following Lemma:
\begin{lemma}\label{renjianzhendaoqushi}
Let $\widehat{\delta}=\min\{\delta/2, \epsilon_{2}\}$ where $\delta$
and $\epsilon_{2}$ are the constants appeared in Lemma \ref{dahuoe}
and Lemma \ref{dadaoruqintain5tqddu} respectively, then
\begin{equation}\label{hapnaner}
H_{j(k)+1}(J^{c_{0}+\beta_{0}}_{k},J^{\widehat{\delta}}_{k})\neq 0.
\end{equation}
\end{lemma}

By Lemma \ref{shanshan} we know that $K^{(i)}_{k},\
i=1,2,\cdots,m_{k}$ are isolated critical sets of $J_{k}$. In
\cite{CG}, Chang and Ghoussoub provided a definition of critical
group for isolated critical set. In \cite{gm} and \cite{cha}, the
authors defined the Gromoll-Meyer pair (for short GM-pair) for an
isolated critical point for  $C^{1}-$functional. And in \cite{CG},
Chang and Ghoussoub extended the definition of GM-pair into a
dynamically isolated critical set (see Definition I.10 of
\cite{CG}).

\medskip

Let $f$ be a $C^{1}$ functional on a Finsler manifold $M$ (Banach
space is a special case of Finsler manifold) with critical set
$K_{f}$. And let $V$ be a pesudo-gradient vector field $V$ with
respct to $df$ on $M$. A pesudo-gradient flow associated with $V$ is
the unique solution of the following ordinary differential equation
in $M:$
$$\dot{\eta}= V_{1}(\eta(x,t)),\ \eta(x,0)=x,$$
where $V_{1}(x)=g(x)\frac{V(x)}{||V(x)||}$ and
$g(x)=\min\{\mbox{dist}(x,K_{f}), 1\}$. A subset $W$ of $M$ is said
to have the mean value property (for short (MVP)) if for any $x\in
M$ and any $t_{0}<t_{1}$ we have $\eta(x,[t_{0},t_{1}])\subset W$
whenever  $\eta(x,t_{i})\in W,\ i=1,2.$

\begin{definition}\label{afsd556ggbdvf}(Definition I.10 of \cite{CG})
Let $f$ be a $C^{1}$ functional on a Finsler manifold $M$. A subset
$S$ of the critical set $K$ of $f$ is said to be a dynamically
isolated critical set if there exist a closed neighborhood
$\mathcal{O}$ of $S$ and regular values $a<b$ of $f$ such that
\begin{equation}\label{homninchunenggqw545ff}
\mathcal{O}\subset f^{-1}[a,b]
\end{equation}
and \begin{equation}\label{mnnavddf99kh}
cl(\widetilde{\mathcal{O}})\cap K\cap f^{-1}[a,b]=S,
\end{equation}
where
$\widetilde{\mathcal{O}}=\bigcup_{t\in\mathbb{R}}\eta(\mathcal{O},t)$.
$(\mathcal{O}, a, b)$ is called an isolating triplet for $S.$
\end{definition}
After providing the definition of dynamically isolated critical set,
the authors of \cite{CG} give the definition of critical group for
dynamically isolated critical set as follows:
\begin{definition}\label{bbbnhaaamhj88jujull}
Let $S$ be a dynamically isolated critical set of a $C^{1}$
functional $f$ and let $(\mathcal{O},a,b)$ be any isolating triplet
for $S$. For each integer $q$, we shall call the $q$th homology goup
$$C_{q}(f,S)=H_{q}(f^{b}\cap\widetilde{\mathcal{O}}^{+},f^{a}\cap\widetilde{\mathcal{O}}^{+})$$
the $q$th critical group for $S$, where
$\widetilde{\mathcal{O}}^{+}=\bigcup_{t\geq 0}\eta(\mathcal{O},t).$
\end{definition}
\begin{remark}
In \cite{CG}, the  critical group is defined by the  cohomology
group of the topology pair
$(f^{b}\cap\widetilde{\mathcal{O}}^{+},f^{a}\cap\widetilde{\mathcal{O}}^{+})$.
Here we use the homology group instead. All results in \cite{CG}
still holds for homology group since the properties of cohomology
the authors used in \cite{CG} are excision property and homotopy
property.
\end{remark}
\begin{definition}\label{haoyouduoooouagfdf}(Definition III.1
of \cite{CG})\label{zheqingsabainiaotefr66yh} Let $f$ be a $C^{1}$
functional on a Finsler manifold $M$ and let $S$ be a subset of the
critical set $K_{f}$ for $f$. A pair $(W,W_{-})$ of subset is said
to be a GM-pair for $S$ associated with a pesudo-gradient  vector
field $V$, if the following conditions hold:
\begin{description}
\item{(1).} $W$ is a closed (MVP) neighborhood of $S$ satisfying
$W\cap K=S$ and $W\cap f_{\alpha}=\emptyset$ for some $\alpha.$

\item{(2).} $W_{-}$ is an exit set for $W,$ i.e., for each $x_{0}\in W$
and $t_{1}>0$ such that $\eta(x_{0},t_{1})\not\in W,$ there exists
$t_{0}\in [0,t_{1})$ such that $\eta(x_{0},[0,t_{0}])\subset W$ and
$\eta(x_{0},t_{0})\in W_{-}.$

\item{(3).} $W_{-}$ is closed and is a union of a finite number of
sub-manifolds that transversal to the flow $\eta.$
\end{description}
\end{definition}

In \cite{CG}, the authors proved the following theorem which can be
seen as another definition of the critical group for a dynamically
isolated critical set.
\begin{lemma}(Theorem III.3 of \cite{CG})\label{wauxianfenguanzai}
Let $f$ be a $C^{1}$ functional on a $C^{1}$ Finsler manifold $M$
and let $S$ be a dynamically isolated critical set for $f.$ Then for
any GM-pair $(W,W_{-})$ for $S,$ we have
$$H_{*}(W,W_{-})\cong H_{*}(f^{b}\cap\widetilde{\mathcal{O}}^{+},f^{a}\cap\widetilde{\mathcal{O}}^{+})
=C_{*}(f,S),$$ where $(\mathcal{O}, a,b)$ is an isolating triplet
for $S.$
\end{lemma}

By \cite{CG} and Lemma \ref{llkaddas}, we have the following
 Lemma
 \begin{lemma}\label{hhhnayooahj88765kk}
There exists an index $i_{k}$ satisfying $1\leq i_{k}\leq m_{k}$
such that
\begin{equation}\label{hunanren}
C_{j(k)+1}(J_{k},K^{(i_{k})}_{k})\neq 0,\ k=1,2,\cdots.
\end{equation}
\end{lemma}
\noindent{\bf Proof.} It is easy to verify that
$((J_{k})^{c_{0}+\beta_{0}}_{\widehat{\delta}},
c_{0}+\beta_{0},\widehat{\delta})$ is an isolating triplet for the
isolated critical set $K^{c_{0}+\beta_{0}}_{k}$, where
$\widehat{\delta}$ is the constant appeared in Lemma
\ref{renjianzhendaoqushi} and
$$(J_{k})^{c_{0}+\beta_{0}}_{\widehat{\delta}}=\{u\in E_{k}\
:\ \widehat{\delta}\leq J_{k}(u)\leq c_{0}+\beta_{0}\}.$$ Thus by
the Definition \ref{bbbnhaaamhj88jujull}  and Lemma
\ref{renjianzhendaoqushi}, we have
\begin{equation}\label{jnh7hytqq4f8i6x}
C_{j(k)+1}(J_{k},K^{c_{0}+\beta_{0}}_{k})
=H_{j(k)+1}(J^{c_{0}+\beta_{0}}_{k},J^{\widehat{\delta}}_{k})\neq
0.
\end{equation}
Let
$$\mathcal{C}=\{\sum^{j}_{i=1}J(v_{i})
\ |\ 1\leq i\leq j,\ 1\leq j\leq\overline{l},\ v_{i}\in
K^{c_{0}+\beta_{0}}\}\cap(-\infty, c_{0}+\beta_{0}].$$ By the
condition $\bf {(*)}$ we know that $\mathcal{C}$ is a finite set.
Without loss of generality, we may assume that
$$\mathcal{C}=\{c_{1},\cdots,c_{n_{0}}\}.$$
and $c_{1}<c_{2}<\cdots<c_{n_{0}}.$ Thus by Lemma 2.7, we get that
\begin{equation}\label{vvvxcsdsds}
\mbox{dist}_{\mathbb{R}}(J_{k}(K^{c_{0}+\beta_{0}}_{k}),\mathcal{C})\rightarrow
0,\ k\rightarrow\infty.
\end{equation}
Choose $\epsilon_{0}=\frac{1}{2}\min\{c_{i}-c_{i-1}\ |\
i=2,\cdots,n_{0}\}$. For $0<\epsilon<\epsilon_{0},$ let
$K^{(i,j)}_{k}=K^{(i)}_{k}\cap
(J_{k})^{c_{j}+\epsilon}_{c_{j}-\epsilon},$ $ j=1,\cdots,d_{i}, $
$i=1,\cdots,m_{k}$. By Lemma 2.8 and Remark 2.9, we know that if $k$
large enough, $K^{(i,j)}_{k}$ does not dependent on $\epsilon.$

By Proposition 2.2 of \cite{Rabi}, we get that for any
$0<\delta<\mu$ (for the definition of $\mu$ see (2.11) in the proof
of Lemma 2.8), there exists constant $\varsigma_{\delta}>0$ such
that
\begin{equation}\label{yyyyyqweeasasa}
\inf\{||\nabla J(u)||\ | \ u\in J^{c_{0}+\alpha_{0}}\setminus
N_{\delta}(\mathcal{F}_{\overline{l}}(K^{c_{0}+\alpha_{0}}))\}>0.
\end{equation}
Then by Lemma 2.7, we get that there exists constant
$\varsigma_{\delta}>0$  which is independent of $k$, such that
\begin{equation}\label{yyyyyqweeasasa}
\inf\{||\nabla J_{k}(u)||\ | \ u\in
J^{c_{0}+\alpha_{0}}_{k}\setminus
N_{\delta}(K^{c_{0}+\alpha_{0}}_{k})\}>\varsigma_{\delta}>0.
\end{equation}

By Lemma \ref{shanshan}, we get that
$\mbox{diam}(K^{(i)}_{k})\rightarrow 0$ as $k\rightarrow\infty$ and
$\mbox{dist}(K^{(i)}_{k},K^{(j)}_{k})\geq\delta_{0},\ i\neq j$, when
$k$ large enough. Thus by (\ref{vvvxcsdsds}),
(\ref{yyyyyqweeasasa}), Section 2 of \cite{gm} or page 49 and page
50 of \cite{cha}, we know that  we can choose $\delta\in
(0,\delta_{0})$ and $\epsilon\in (0,\epsilon_{0})$ such that there
exist GM-pair $(W_{i,j},W^{-}_{i,j})$ of $K^{(i,j)}_{k}$ such that
$$W_{i}\subset \{u\in E_{k}\ :\ \mbox{dist}(u, K^{(i)}_{k})\leq\delta/4\},
\ j=1,\cdots,d_{i},\  i=1,2,\cdots, m_{k}.$$ Thus
$\displaystyle(\bigcup
^{m_{k}}_{i=1}W_{i},\bigcup^{m_{k}}_{i=1}W^{-}_{i})$ is a GM-pair of
$K^{c_{0}+\beta_{0}}_{k}$. By (\ref{jnh7hytqq4f8i6x}) and Lemma
\ref{wauxianfenguanzai}, we get that
\begin{equation}\label{honhonlalamnnansiao8jb}
H_{j(k)+1}\displaystyle(\bigcup
^{m_{k}}_{i=1}W_{i},\bigcup^{m_{k}}_{i=1}W^{-}_{i})\cong
C_{j(k)+1}(J_{k},K^{c_{0}+\beta_{0}}_{k})=
H_{j(k)+1}(J^{c_{0}+\beta_{0}}_{k},J^{\widehat{\delta}}_{k})\neq 0.
\end{equation}
Since $W_{i}$ and $W_{j}$ are disjoint if $i\neq j,$ we get that
\begin{equation}\label{zhichenshuiyure8b}
H_{j(k)+1}\displaystyle(\bigcup
^{m_{k}}_{i=1}W_{i},\bigcup^{m_{k}}_{i=1}W^{-}_{i})\cong\bigoplus^{m_{k}}_{i=1}H_{j(k)+1}(W_{i},W^{-}_{i}).
\end{equation}
By Lemma \ref{wauxianfenguanzai}, we have
\begin{equation}\label{yueshijiaxmin8hyt}
H_{j(k)+1}(W_{i},W^{-}_{i})\cong
C_{j(k)+1}(J_{k},K^{(i)}_{k}).\end{equation} By
$(\ref{jnh7hytqq4f8i6x})-(\ref{yueshijiaxmin8hyt})$, we get that
\begin{equation}\label{burezhaohuanj89maz8}
\bigoplus^{m_{k}}_{i=1}C_{j(k)+1}(J_{k},K^{(i)}_{k})\cong
C_{j(k)+1}(J_{k},K^{c_{0}+\beta_{0}}_{k})\neq 0.
\end{equation}
 Thus there exists an index $i_{k}$ satisfying $1\leq i_{k}\leq m_{k}$ such that
$C_{j(k)+1}(J_{k},K^{(i_{k})}_{k})\neq 0.$\hfill$\Box$

\section{A reduction method}

\hspace*{\parindent}Let $K^{(i_{k})}_{k}$ satisfy (\ref{hunanren})
and $u_{k}\in K^{(i_{k})}_{k},\ k=1,2\cdots.$ By Lemma
\ref{cotgaadsewe}, we know that there exist a positive integer $n$,
$n$ functions $u^{i}\in K^{c_{0}+\beta_{0}},\ i=1,\cdots,n$ and $n$
sequences $\{b^{i}_{k}\}\subset \mathbb{Z}^{N},\ i=1,\cdots,n$ such
that as $k\rightarrow\infty,$
\begin{equation}\label{jiajiajiajnh7yhqq}
|b^{i}_{k}-b^{j}_{k}|\rightarrow+\infty,\ \mbox{for any}\ i\neq j\
\mbox{and}\
||u_{k}-\sum^{n}_{i=1}u^{i}(\cdot+b^{i}_{k})||_{H^{1}(Q_{k})}\rightarrow
0.\end{equation} Without loss of generality, we may assume that
$|b^{i}_{k}|\rightarrow +\infty$ as $k\rightarrow\infty,\
i=1,\cdots,n.$
\medskip

The proofs of the  following four Lemmas shall be provided in the
appendix. For convenience, we denote $f'_{t}(x,t)$ by $f'(x,t).$

\medskip

\begin{lemma} \label{mnnnssshgdfr}The following limit
holds uniformly  for any $\psi_{k}, \varphi_{k}\in E_{k}$ which
satisfy $ ||\psi_{k}||_{k}\leq 1,\ ||\varphi_{k}||_{k}\leq 1$,
$$\lim_{k\rightarrow\infty}\int_{Q_{k}}|f'(x,u_{k})-f'(x,\sum^{n}_{i=1}u^{i}(\cdot+b^{i}_{k}))|\cdot
|\psi_{k}|\cdot |\varphi_{k}|dx=0.$$
\end{lemma}

\begin{lemma}\label{yadanhexiawa}
If $\{\widetilde{v}_{k}\}$ and $\{v_{k}\}$ are two bounded sequences
in $H^{1}(\mathbb{R}^{N})$ which satisfy that
$||\widetilde{v}_{k}-v_{k}||\rightarrow 0$ as $k\rightarrow\infty,$
then the following two results holds:
\begin{description}

\item{(1).} The Limit
$$\int_{\mathbb{R}^{N}}|f(x,\widetilde{v}_{k})-f(x,v_{k})|\cdot
|\varphi|\rightarrow 0,\ k\rightarrow\infty$$ holds uniformly for
 $\varphi\in H^{1}(\mathbb{R}^{N})$ which satisfies $||\varphi||\leq 1$.

\item{(2).} The Limit
$$\int_{\mathbb{R}^{N}}|f'(x,\widetilde{v}_{k})-f'(x,v_{k})|\cdot
|\varphi\cdot\psi|\rightarrow 0,\ k\rightarrow\infty$$ holds
uniformly for  $\varphi,\psi\in H^{1}(\mathbb{R}^{N})$ which satisfy
$||\varphi||\leq 1$, $||\psi||\leq 1$.
\end{description}
\end{lemma}

\begin{lemma}\label{lantaojin}
As $k\rightarrow\infty,$ the limit
$$\int_{Q_{k}}|f'(x,\sum^{n}_{i=1}u^{i}(\cdot+b^{i}_{k}))-
\sum^{n}_{i=1}f'(x,u^{i}(\cdot+b^{i}_{k}))|\cdot
|\varphi_{k}\cdot\psi_{k}|dx\rightarrow 0$$ holds uniformly for any
$\psi_{k}, \varphi_{k}\in E_{k}$ which satisfy that
$||\varphi_{k}||_{k}\leq 1,\ ||\psi_{k}||_{k}\leq 1.$
\end{lemma}

\begin{lemma}\label{tongrenttthqqqw}
\begin{description}
\item{(1).} Suppose $v^{i}_{k}\in E_{k},$ $a^{i}_{k}\in
\mathbb{Z}^{N},\ i=1,2,\cdots,n$ which satisfy that
 for any $i\neq j,$
$|a^{i}_{k}-a^{j}_{k}|\rightarrow\infty$  as $k\rightarrow\infty$.
Then as $k\rightarrow\infty,$ the limit
$$\int_{Q_{k}}|f'(x,\sum^{n}_{i=1}v^{i}_{k}(\cdot+a^{i}_{k})-
\sum^{n}_{i=1}f'(x,v^{i}_{k}(\cdot+a^{i}_{k})))|\cdot
|\varphi_{k}\cdot\psi_{k}|dx\rightarrow 0$$ holds uniformly for any
$\psi_{k}, \varphi_{k}\in E_{k}$ which satisfy
$||\varphi_{k}||_{k}\leq 1,\ ||\psi_{k}||_{k}\leq 1.$

\item{(2).} Suppose $v^{i}\in H^{1}(\mathbb{R}^{N}),$ $a^{i}_{k}\in
\mathbb{Z}^{N},\ i=1,2,\cdots,n$ which satisfy that
 for any $i\neq j,$
$|a^{i}_{k}-a^{j}_{k}|\rightarrow\infty$  as $k\rightarrow\infty$.
Then as $k\rightarrow\infty,$ the limit
$$\int_{\mathbb{R}^{N}}|f'(x,\sum^{n}_{i=1}v^{i}(\cdot+a^{i}_{k})-
\sum^{n}_{i=1}f'(x,v^{i}(\cdot+a^{i}_{k})))|\cdot
|\varphi\cdot\psi|dx\rightarrow 0$$ holds uniformly for any  $\psi,\
\varphi\in H^{1}(\mathbb{R}^{N})$ which satisfy  $||\varphi||\leq
1,\ ||\psi||\leq 1.$
\end{description}
\end{lemma}
\begin{lemma}\label{fujianren}
Suppose that $v_{k}\in\ker \nabla^{2}J_{k}(u_{k})$ and
$||v_{k}||_{k}=1,\ k=1,2,\cdots.$ Then there exist $v^{i}\in \ker
\nabla^{2}J(u^{i}),\ i=1,\cdots,n$ such that as
$k\rightarrow\infty,$
$$||v_{k}-\sum^{n}_{i=1}v^{i}(\cdot+b^{i}_{k})||_{H^{1}(Q_{k})}\rightarrow 0.$$
\end{lemma}
\noindent{\bf Proof.} Since $v_{k}\in\ker \nabla^{2}J_{k}(u_{k}),$
we have \begin{equation}\label{chulaidaoqiulaidao} -\triangle
v_{k}+V(x)v_{k}=f'(x,u_{k})v_{k},\ \mbox{in}\ E_{k}.\end{equation}
By Lemma \ref{mnnnssshgdfr} and Lemma \ref{lantaojin}, we know that
the limit
$$\lim_{k\rightarrow\infty}\int_{Q_{k}}|f'(x,u_{k})-\sum^{n}_{i=1}f'(x,u^{i}(\cdot+b^{i}_{k}))|
\cdot |\varphi_{k}|\cdot |\psi_{k}|dx=0$$ holds  uniformly for
$\varphi_{k}, \phi_{k}\in E_{k}$ which satisfy $
||\varphi_{k}||_{k}\leq 1,\ ||\psi_{k}||_{k}\leq 1.$ Assume that
$v_{k}(\cdot-b^{i}_{k})\rightharpoonup v^{i}$ in
$H^{1}_{loc}(\mathbb{R}^{N})$ as $k\rightarrow\infty.$ By
$||v_{k}||_{k}=1$, we deduce that $v^{i}\in H^{1}(\mathbb{R}^{N})$.

By
$$-\triangle
v_{k}(\cdot-b^{i}_{k})+V(x)v_{k}(\cdot-b^{i}_{k})=f'(x,u_{k}(\cdot-b^{i}_{k}))v_{k}(\cdot-b^{i}_{k})\
\mbox{in}\ E_{k},$$ $u_{k}(\cdot-b^{i}_{k})\rightharpoonup u^{i}$
and $v_{k}(\cdot-b^{i}_{k})\rightharpoonup v^{i}$ in
$H^{1}_{loc}(\mathbb{R}^{N}),$ we get that
\begin{equation}\label{hundandgt5fr0olk9m}
-\triangle v^{i}+V(x)v^{i}=f'(x,u^{i})v^{i}\ \mbox{in}\
H^{1}(\mathbb{R}^{N}),
\end{equation}
i.e. $v^{i}\in \ker \nabla^{2}J(u^{i}),\ i=1,\cdots,n.$

Recall that $\chi_{k}$ be  cut-off functions satisfying that
$0\leq\chi_{k}\leq 1,\ \chi_{k}\equiv 1$ on $Q_{k-1}, \
\chi_{k}\equiv0$ outside of $Q_{k}$ and $|\nabla \chi_{k}|\leq C.$
Set $v^{i}_{k}=\chi_{k}v^{i},$ we have $v^{i}_{k}\in E_{k}.$ If
$\varphi\in E_{k}$, then by (\ref{hundandgt5fr0olk9m}), we have
\begin{eqnarray}\label{xiaoshanshan} &&\int_{Q_{k}}\nabla
v^{i}_{k}\cdot\nabla\varphi+\int_{Q_{k}}V(x)v^{i}_{k}\varphi\nonumber\\
&=&\int_{Q_{k}}\nabla v^{i}\cdot\nabla
(\chi_{k}\varphi)+\int_{Q_{k}}V(x)v^{i}\cdot(\chi_{k}\varphi)
-\int_{Q_{k}}\varphi\nabla
v^{i}\nabla\chi_{k}+\int_{Q_{k}}v^{i}\nabla
\chi_{k}\nabla\varphi\nonumber\\
&=&\int_{Q_{k}}f'(x,u^{i})\cdot
v^{i}\cdot(\chi_{k}\varphi)-\int_{Q_{k}}\varphi\nabla
v^{i}\nabla\chi_{k}+\int_{Q_{k}}v^{i}\nabla\chi_{k}\nabla\varphi\nonumber\\
&=&\int_{Q_{k}}f'(x,u^{i})\cdot
(\chi_{k}v^{i})\cdot\varphi-\int_{Q_{k}}\varphi\nabla
v^{i}\nabla\chi_{k}+\int_{Q_{k}}v^{i}\nabla\chi_{k}\nabla\varphi\nonumber\\
&=&\int_{Q_{k}}f'(x,u^{i})\cdot
v^{i}_{k}\cdot\varphi-\int_{\mathbb{R}^{N}\setminus
Q_{k}}\varphi\nabla
v^{i}\nabla\chi_{k}+\int_{\mathbb{R}^{N}\setminus
Q_{k}}v^{i}\nabla\chi_{k}\nabla\varphi.\nonumber
\end{eqnarray}
It follows that for any $\varphi\in E_{k},$
\begin{eqnarray}\label{zhiqiangbuxi}
&&\int_{Q_{k}}\nabla
v^{i}_{k}(\cdot+b^{i}_{k})\nabla\varphi+\int_{Q_{k}}V(x)v^{i}_{k}(\cdot+b^{i}_{k})\varphi\nonumber\\
&=&\int_{Q_{k}}f'(x,u^{i}(\cdot+b^{i}_{k}))\cdot
v^{i}_{k}(\cdot+b^{i}_{k})\cdot\varphi -\int_{Q_{k}}(\nabla
v^{i}(\cdot+b^{i}_{k})\cdot\nabla\chi_{k}(\cdot+b^{i}_{k}))\cdot\varphi\nonumber\\
&&
+\int_{Q_{k}}v^{i}(\cdot+b^{i}_{k})\nabla\chi_{k}(\cdot+b^{i}_{k})\cdot\nabla\varphi.
\end{eqnarray}
By (\ref{zhiqiangbuxi}), we get that
\begin{eqnarray}\label{chuanshuozhong} &&\int_{Q_{k}}\nabla
v^{i}_{k}(\cdot+b^{i}_{k})\nabla\varphi+\int_{Q_{k}}V(x)v^{i}_{k}(\cdot+b^{i}_{k})\varphi\nonumber\\
&=&\int_{Q_{k}}f'(x,u^{i}(\cdot+b^{i}_{k}))\cdot
v_{k}\cdot\varphi\nonumber\\
&&+\int_{Q_{k}}f'(x,u^{i}(\cdot+b^{i}_{k}))\cdot
v^{i}_{k}(\cdot+b^{i}_{k})\cdot\varphi-\int_{Q_{k}}f'(x,u^{i}(\cdot+b^{i}_{k}))\cdot
v_{k}\cdot\varphi\nonumber\\
&&+\int_{Q_{k}}v^{i}(\cdot+b^{i}_{k})\nabla\chi_{k}(\cdot+b^{i}_{k})\cdot\nabla\varphi
 -\int_{Q_{k}}(\nabla
v^{i}(\cdot+b^{i}_{k})\cdot\nabla\chi_{k}(\cdot+b^{i}_{k}))\cdot\varphi.\nonumber\\
\end{eqnarray}
By (\ref{chulaidaoqiulaidao}) and (\ref{chuanshuozhong}), we get
that for any $\varphi\in E_{k},$
\begin{eqnarray}\label{renshengrumena}
&&\int_{Q_{k}}\nabla
(v_{k}-\sum^{n}_{i=1}v^{i}_{k}(\cdot+b^{i}_{k}))\cdot\nabla\varphi
+\int_{Q_{k}}V(x)\cdot(v_{k}-\sum^{n}_{i=1}v^{i}_{k}(\cdot+b^{i}_{k}))\cdot\varphi\nonumber\\
&=&\int_{Q_{k}}(f'(x,u_{k})-\sum^{n}_{i=1}f'(x,u^{i}(\cdot+b^{i}_{k})))\cdot
v_{k}\cdot\varphi\nonumber\\
&&-\sum^{n}_{i=1}\int_{Q_{k}}f'(x,u^{i}(\cdot+b^{i}_{k}))
\cdot(v^{i}_{k}(\cdot+b^{i}_{k})-v_{k})\cdot\varphi\nonumber\\
&&-\sum^{n}_{i=1}\int_{Q_{k}}v^{i}(\cdot+b^{i}_{k})
\nabla\chi_{k}(\cdot+b^{i}_{k})\cdot\nabla\varphi
+\sum^{n}_{i=1}\int_{Q_{k}}(\nabla
v^{i}(\cdot+b^{i}_{k})\cdot\nabla\chi_{k}(\cdot+b^{i}_{k}))\cdot\varphi.\nonumber\\
\end{eqnarray}
Since $v_{k}(\cdot-b^{i}_{k})\rightharpoonup v^{i}$ in
$H^{1}_{loc}(\mathbb{R}^{N})$ as $k\rightarrow\infty,$ we get that
the limit
\begin{eqnarray}\label{akakbbbgtsewdes}
\lim_{k\rightarrow\infty}\int_{Q_{k}}f'(x,u^{i}(\cdot+b^{i}_{k}))\cdot
(v^{i}_{k}(\cdot+b^{i}_{k})-v_{k})\cdot\varphi= 0
\end{eqnarray}
holds uniformly for $\varphi\in E_{k}$ which satisfies
$||\varphi||_{k}\leq 1.$ By Lemma \ref{mnnnssshgdfr} and Lemma
\ref{lantaojin}, we know that that the following limit holds
uniformly for $\varphi\in E_{k} $ which satisfies $
||\varphi||_{k}\leq 1,$
\begin{equation}\label{juju7hy13246gf}
\lim_{k\rightarrow\infty}\int_{Q_{k}}(f'(x,u_{k})
-\sum^{n}_{i=1}f'(x,u^{i}(\cdot+b^{i}_{k})))\cdot
v_{k}\cdot\varphi=0.
\end{equation}
Furthermore, the following limits hold uniformly for $\varphi\in
E_{k} $ which satisfies $ ||\varphi||_{k}\leq 1,$
\begin{equation}\label{7653100087sd9}
\lim_{k\rightarrow\infty}\int_{Q_{k}}v^{i}(\cdot+b^{i}_{k})
\nabla\chi_{k}(\cdot+b^{i}_{k})\cdot\nabla\varphi=\lim_{k\rightarrow\infty}
\int_{\mathbb{R}^{N}\setminus Q_{k}}v^{i}(\cdot+b^{i}_{k})
\nabla\chi_{k}(\cdot+b^{i}_{k})\cdot\nabla\varphi=0,\end{equation}

\begin{equation}\label{nb5f5f5yg88}
\lim_{k\rightarrow\infty}\int_{Q_{k}}(\nabla
v^{i}(\cdot+b^{i}_{k})\cdot\nabla\chi_{k}(\cdot+b^{i}_{k}))\cdot\varphi
=\lim_{k\rightarrow\infty} \int_{\mathbb{R}^{N}\setminus
Q_{k}}(\nabla
v^{i}(\cdot+b^{i}_{k})\cdot\nabla\chi_{k}(\cdot+b^{i}_{k}))\cdot\varphi
=0.\end{equation} By $(\ref{renshengrumena})-(\ref{nb5f5f5yg88})$,
we deduce  that the following equality holds uniformly for
$\varphi\in E_{k} $ which satisfies $ ||\varphi||_{k}\leq 1,$
\begin{equation}\label{bianshuiliusishuiliu}
\int_{Q_{k}}\nabla
(v_{k}-\sum^{n}_{i=1}v^{i}_{k}(\cdot+b^{i}_{k}))\cdot\nabla\varphi
+\int_{Q_{k}}V(x)\cdot(v_{k}-\sum^{n}_{i=1}v^{i}_{k}(\cdot+b^{i}_{k}))\cdot\varphi
=o(1),\ \mbox{as}\ k\rightarrow\infty.
\end{equation}

Recall that $P_{k}:E_{k}\rightarrow Y_{k}$ and
$T_{k}:E_{k}\rightarrow Z_{k}$ are the orthogonal projections.
Choose
$\varphi=T_{k}(v_{k}-\sum^{n}_{i=1}v^{i}_{k}(\cdot+b^{i}_{k}))$ and
$\varphi=P_{k}(v_{k}-\sum^{n}_{i=1}v^{i}_{k}(\cdot+b^{i}_{k}))$
respectively in (\ref{bianshuiliusishuiliu}), we get that  as
$k\rightarrow\infty,$
$$||T_{k}(v_{k}-\sum^{n}_{i=1}v^{i}_{k}(\cdot+b^{i}_{k}))||_{k}\rightarrow 0,\
||P_{k}(v_{k}-\sum^{n}_{i=1}v^{i}_{k}(\cdot+b^{i}_{k}))||_{k}\rightarrow
0.$$ Therefore, as $k\rightarrow\infty,$
$||v_{k}-\sum^{n}_{i=1}v^{i}_{k}(\cdot+b^{i}_{k})||_{k}\rightarrow
0.$ \hfill$\Box$

\bigskip
The following notations will be used in this   and the next
Sections.

\begin{itemize}
 \item Let $N_{i}=\ker
\nabla^{2}J(u^{i})=\mbox{span}\{e_{i,1},\cdots,e_{i,l_{i}}\}$, where
$l_{i}=\dim N_{i},\ i=1,2,\cdots,n$ and
$(e_{i,s},e_{i,j})_{X}=\delta_{s,j}.$

\item Let $N^{k}_{i}=\mbox{span}\{(\chi_{k}e_{i,1})(\cdot+b^{i}_{k}),
\cdots,(\chi_{k}e_{i,l_{i}})(\cdot+b^{i}_{k})\}\subset E_{k},\
i=1,2,\cdots,n,$
$\Lambda_{k}=\mbox{span}\{\displaystyle\bigcup^{n}_{i=1}N^{k}_{i}\}\subset
E_{k}$ and $\Pi_{k}=\Lambda^{\bot}_{k}$, the orthogonal complement
space of $\Lambda_{k}$ in $E_{k}.$

 \item Let $\widetilde{N}^{k}_{i}
=\mbox{span}\{e_{i,1}(\cdot+b^{i}_{k}),\cdots,e_{i,l_{i}}(\cdot+b^{i}_{k})\},\
i=1,\cdots,n,$
$\widetilde{\Lambda}_{k}=\mbox{span}\{\displaystyle\bigcup^{n}_{i=1}\widetilde{N}^{k}_{i}\}\subset
X$ and $\widetilde{\Pi}_{k}=(\widetilde{\Lambda}_{k})^{\bot}$, the
orthogonal complement space of $\widetilde{\Lambda}_{k}$ in $X$.

 \item Let
$\widehat{N}^{k}_{i}=\mbox{span}\{\chi_{k}e_{i,1},\cdots,\chi_{k}e_{i,l_{i}}\}\subset
E_{k},$ $ i=1,2,\cdots,n.$

 \item For convenience, denote $(\chi_{k}u^{i})(\cdot+b^{i}_{k}),\
\chi_{k}u^{i},\ (\chi_{k}e_{i,j})(\cdot+b^{i}_{k})$ and
$e_{i,j}(\cdot+b^{i}_{k})$ by $u^{i}_{k}, \ \widehat{u}^{i}_{k},$
$e^{k}_{i,j}$ and $\widetilde{e}^{k}_{i,j}$ respectively, $
j=1,2,\cdots,l_{i},\ i=1,2,\cdots,n.$
\end{itemize}

\begin{lemma}\label{xiaojiazhi}
There exists $k_{0}\in\mathbb{N},\ \delta_{0}>0$ and $\eta>0$ which
are independent of $k$ such that
\begin{description}
\item{(1).} if $k\geq k_{0},$ then for any $u\in
B_{E_{k}}(u_{k},\delta_{0}),$ the operator
$$\widetilde{T}_{k}\circ\nabla^{2}J_{k}(u)|_{\Pi_{k}}:
\Pi_{k}\rightarrow \Pi_{k}$$ is invertible and
$$||(\widetilde{T}_{k}\circ\nabla^{2}J_{k}(u)|_{\Pi_{k}})^{-1}||\leq\eta,\ k=1,2,\cdots,$$
where $\widetilde{T}_{k}:E_{k}\rightarrow\Pi_{k}$ is the orthogonal
projection.

\item{(2).} if $k\geq k_{0},$ then for any $u\in
B_{E_{k}}(\widehat{u}^{i}_{k},\delta_{0}),$ the operator
$$P_{(\widehat{N}^{k}_{i})^{\bot}}\circ\nabla^{2}J_{k}(u)|_{(\widehat{N}^{k}_{i})^{\bot}}:
(\widehat{N}^{k}_{i})^{\bot}\rightarrow
(\widehat{N}^{k}_{i})^{\bot}$$ is invertible and
$$||(P_{(\widehat{N}^{k}_{i})^{\bot}}\circ\nabla^{2}J_{k}(u)
|_{(\widehat{N}^{k}_{i})^{\bot}})^{-1}||\leq\eta,\ k=1,2,\cdots,$$
where
$P_{(\widehat{N}^{k}_{i})^{\bot}}:E_{k}\rightarrow(\widehat{N}^{k}_{i})^{\bot}$
is the orthogonal projection.
\end{description}
\end{lemma}
\noindent{\bf Proof.} We only give the proof of the result (1),
since the proof of the result (2) is similar. If the result (1) is
not true, then there exists a sequence $\{\widetilde{u}_{k}\}$ such
that $\widetilde{u}_{k}\in E_{k}$ and as $k\rightarrow\infty$,
\begin{equation}\label{qeqeqeq}
||u_{k}-\widetilde{u}_{k}||_{k}\rightarrow 0,\
||\widetilde{T}_{k}\circ\nabla^{2}J_{k}(\widetilde{u}_{k})|_{\Pi_{k}}||\rightarrow
0.\end{equation} By (\ref{qeqeqeq}), we know that there exists
$v_{k}\in\Pi_{k}$ satisfying $||v_{k}||_{k}=1$ and as
$k\rightarrow\infty,$
\begin{equation}\label{saomiaommm}
||\widetilde{T}_{k}(\nabla^{2}J_{k}(\widetilde{u}_{k})v_{k})||_{k}\rightarrow
0.
\end{equation}

\noindent{\bf Step 1.} We shall prove that if
$v_{k}(\cdot-b^{i}_{k})\rightharpoonup v^{i}$ in
$H^{1}_{loc}(\mathbb{R}^{N})$, then
$v^{i}\in\ker\nabla^{2}J(u^{i}),\ i=1,\cdots,n.$

Choose $\varphi\in (\ker\nabla^{2}J(u^{i}))^{\bot}.$ Let
$\varphi_{k}=\chi_{k}\varphi.$ Then $\varphi_{k}\in E_{k}.$ Assume
that
$$\varphi_{k}(\cdot+b^{i}_{k})=\widetilde{\varphi}_{k}+\widehat{\varphi}_{k},$$
where $\widetilde{\varphi}_{k}\in\Pi_{k}$ and
$\widehat{\varphi}_{k}\in\Lambda_{k}.$ Since
$\lim_{k\rightarrow\infty}|b^{i}_{k}|=+\infty$, we deduce that for
any $u\in\ker\nabla^{2}J(u^{i})$, as $k\rightarrow\infty,$ the limit
$$(\varphi_{k}(\cdot+b^{i}_{k}),\chi_{k}u)_{k}\rightarrow 0.$$
It follows that $||\widehat{\varphi}_{k}||_{k}\rightarrow 0$ as
$k\rightarrow\infty$.
 Thus  as $k\rightarrow\infty,$
$$||\widetilde{\varphi}_{k}(\cdot-b^{i}_{k})
-\varphi_{k}||_{H^{1}(Q_{k})}\rightarrow 0.$$ By (\ref{saomiaommm}),
we get that as $k\rightarrow\infty,$
\begin{eqnarray}\label{dabaimmmnn}
&&(\nabla^{2}J_{k}(\widetilde{u}_{k})v_{k},\widetilde{\varphi}_{k})_{k}\nonumber\\
&=&\int_{Q_{k}}\nabla v_{k}\cdot\nabla
\widetilde{\varphi}_{k}+\int_{Q_{k}}V(x)v_{k}\cdot\widetilde{\varphi}_{k}-
\int_{Q_{k}}f'(x,\widetilde{u}_{k})\cdot
v_{k}\cdot\widetilde{\varphi}_{k}\rightarrow 0.
\end{eqnarray}
Note that $v_{k}(\cdot-b^{i}_{k})\rightharpoonup v^{i}$ in
$H^{1}_{loc}(\mathbb{R}^{N})$,
$\widetilde{u}_{k}(\cdot-b^{i}_{k})\rightharpoonup u^{i}$ in
$H^{1}_{loc}(\mathbb{R}^{N})$ and
$\widetilde{\varphi}_{k}(\cdot-b^{i}_{k})\rightarrow \varphi$ in
$H^{1}_{loc}(\mathbb{R}^{N})$, by (\ref{dabaimmmnn}) and
$$\lim_{k\rightarrow\infty}\int_{Q_{k}}f'(x,\widetilde{u}_{k}(\cdot-b^{i}_{k}))\cdot v_{k}(\cdot-b^{i}_{k})
\cdot\widetilde{\varphi}_{k}(\cdot-b^{i}_{k})=
\int_{\mathbb{R}^{N}}f'(x,u^{i})\cdot v^{i}\cdot \varphi,$$ we get
that
\begin{equation}\label{vvvghasssdwe}
\int_{\mathbb{R}^{N}}\nabla
v^{i}\cdot\nabla\varphi+\int_{\mathbb{R}^{N}}V(x)v^{i}\cdot\varphi
-\int_{\mathbb{R}^{N}}f'(x,u^{i})v^{i}\cdot\varphi=0.\end{equation}
Since    $\varphi$ is an arbitrary function in $
(\ker\nabla^{2}J(u^{i}))^{\bot}$, by (\ref{vvvghasssdwe}), we have
$v^{i}\in\ker\nabla^{2}J(u^{i}).$

\bigskip

\noindent{\bf Step 2.} We shall prove that $v^{i}=0,\ i=1,\cdots,n.$

By the result of {\bf Step 1} and the definition of $\Lambda_{k}$,
we know that $\chi_{k}v^{i}(\cdot+b^{i}_{k})\in\Lambda_{k}$. Since
$v_{k}\in\Pi_{k},$ we get that $(v_{k},
\chi_{k}v^{i}(\cdot+b^{i}_{k}))_{k}=0,$ i.e.,
$(v_{k}(\cdot-b^{i}_{k}), \chi_{k}v^{i})_{k}= 0,$ $ k=1,2\cdots$. By
$v_{k}(\cdot-b^{i}_{k})\rightharpoonup v^{i}$ in
$H^{1}_{loc}(\mathbb{R}^{N})$ and $\chi_{k}v^{i}\rightarrow v^{i}$
in $H^{1}_{loc}(\mathbb{R}^{N})$, we get that as
$k\rightarrow\infty$,
$$(v_{k}(\cdot-b^{i}_{k}), \chi_{k}v^{i})_{k}\rightarrow
||v^{i}||^{2}_{H^{1}(\mathbb{R}^{N})}.$$ Thus $v^{i}=0,\
i=1,\cdots,n.$

\bigskip

\noindent{\bf Step 3.} We shall prove that as $k\rightarrow\infty,$
$||\nabla^{2}J_{k}(\widetilde{u}_{k})v_{k}||_{k}\rightarrow 0.$

Since
$||\widetilde{T}_{k}(\nabla^{2}J_{k}(\widetilde{u}_{k})v_{k})||_{k}\rightarrow
0$ as $k\rightarrow\infty$, by the definition of
$\widetilde{T}_{k},$ to prove this claim we only need to prove that
for any $u\in \displaystyle\bigcup^{n}_{i=1}N^{k}_{i},$ as
$k\rightarrow\infty,$
$$(\nabla^{2}J_{k}(\widetilde{u}_{k})v_{k},u)_{k}\rightarrow 0.$$
Without loss of generality, we may assume that
$u=\chi_{k}\varphi(\cdot+b^{i}_{k}),$ where $\varphi\in N_{i}$ for
some $i.$ Then
\begin{eqnarray}\label{kaolulu}
&&(\nabla^{2}J_{k}(\widetilde{u}_{k})v_{k},\chi_{k}\varphi(\cdot+b^{i}_{k}))_{k}
\nonumber\\
&=&\int_{Q_{k}}\nabla v_{k}\nabla
(\chi_{k}\varphi(\cdot+b^{i}_{k}))+\int_{Q_{k}}V(x)v_{k}\cdot
\chi_{k}\varphi(\cdot+b^{i}_{k})-\int_{Q_{k}}f'(x,\widetilde{u}_{k})
v_{k}\cdot \chi_{k}\varphi(\cdot+b^{i}_{k})\nonumber\\
&=&\int_{Q_{k}}\nabla (v_{k}(\cdot-b^{i}_{k}))\cdot\nabla
(\chi_{k}\varphi)+\int_{Q_{k}}V(x)v_{k}(\cdot-b^{i}_{k})\cdot
\chi_{k}\varphi\nonumber\\
&&-\int_{Q_{k}}f'(x,\widetilde{u}_{k}(\cdot-b^{i}_{k}))\cdot
v_{k}(\cdot-b^{i}_{k})\cdot \chi_{k}\varphi.\nonumber
\end{eqnarray}
Since $$v_{k}(\cdot-b^{i}_{k})\rightharpoonup 0\ \mbox{in}\
H^{1}_{loc}(\mathbb{R}^{N}),\ (\mbox{by}\ {\bf Step 2})$$
$$\chi_{k}\varphi\rightarrow\varphi \ \mbox{in}\
H^{1}_{loc}(\mathbb{R}^{N}),$$
$$\widetilde{u}_{k}(\cdot-b^{i}_{k})\rightharpoonup u^{i}\
\mbox{in}\ H^{1}_{loc}(\mathbb{R}^{N}),\ (\mbox{by}\
(\ref{qeqeqeq})\ \mbox{and}\ (\ref{jiajiajiajnh7yhqq}))$$ by
(\ref{vvvghasssdwe}), we get that as $k\rightarrow\infty,$
\begin{eqnarray}\label{77715rffadlllm09855}
&&\lim_{k\rightarrow\infty}(\nabla^{2}J_{k}(\widetilde{u}_{k})v_{k},\chi_{k}\varphi(\cdot+b^{i}_{k}))_{k}
\nonumber\\
&=&\int_{\mathbb{R}^{N}}\nabla
v^{i}\cdot\nabla\varphi+\int_{\mathbb{R}^{N}}V(x)v^{i}\cdot\varphi
-\int_{\mathbb{R}^{N}}f'(x,u^{i})v^{i}\cdot\varphi=
0.\nonumber\end{eqnarray}This proves the result of this step.

\bigskip

\noindent{\bf Step 4.} We are ready to prove that
$||v_{k}||_{k}\rightarrow 0$ as $k\rightarrow\infty.$ Then it
induces a contradiction since $||v_{k}||_{k}=1$ for any $ k.$ We
prove $||P_{k}v_{k}||_{k}\rightarrow 0$ as $k\rightarrow\infty$
firstly.

Since
\begin{equation}\label{n888888876}
(\nabla^{2}J_{k}(\widetilde{u}_{k})v_{k},
P_{k}v_{k})_{k}=(v_{k},P_{k}v_{k})-\int_{Q_{k}}f'(x,\widetilde{u}_{k})v_{k}\cdot
(P_{k}v_{k}),\end{equation} by (\ref{qeqeqeq}),
(\ref{jiajiajiajnh7yhqq}) and Lemma \ref{lantaojin}, we get that as
$k\rightarrow\infty,$
\begin{equation}\label{quququad}
\int_{Q_{k}}f'(x,\widetilde{u}_{k})v_{k}\cdot
(P_{k}v_{k})=\sum^{n}_{i=1}\int_{Q_{k}}f'(x,u^{i}(\cdot+b^{i}_{k}))v_{k}\cdot
(P_{k}v_{k})+o(1).\end{equation} By (\ref{quququad}) and the fact
that $v_{k}(\cdot-b^{i}_{k})\rightharpoonup 0$ in
$H^{1}_{loc}(\mathbb{R}^{N})$ (see {\bf Step 2}), we get that as
$k\rightarrow\infty,$ for $i=1,\cdots,n,$
$$\int_{Q_{k}}f'(x,u^{i}(\cdot+b^{i}_{k}))v_{k}\cdot
(P_{k}v_{k})=\int_{Q_{k}}f'(x,u^{i})v_{k}(\cdot-b^{i}_{k})\cdot
(P_{k}v_{k})(\cdot-b^{i}_{k})\rightarrow 0.$$ Thus as
$k\rightarrow\infty,$
\begin{equation}\label{777tgfad34f}
\int_{Q_{k}}f'(x,\widetilde{u}_{k})v_{k}\cdot
(P_{k}v_{k})=o(1).\end{equation} By {\bf Step 3},
 we get that as $k\rightarrow\infty,$ \begin{equation}\label{bnnzx5r}
  (\nabla^{2}J_{k}(\widetilde{u}_{k})v_{k},
P_{k}v_{k})_{k}\rightarrow 0.\end{equation} By (\ref{n888888876}),
(\ref{777tgfad34f}) and (\ref{bnnzx5r}), we get that as
$k\rightarrow\infty,$
$$-||P_{k}v_{k}||^{2}_{k}=(v_{k},P_{k}v_{k})_{k}\rightarrow 0.$$
In the same way, we can prove that
$||T_{k}v_{k}||^{2}_{k}\rightarrow 0$ as $k\rightarrow\infty.$  Thus
$||v_{k}||_{k}=\sqrt{||T_{k}v_{k}||^{2}_{k}+||P_{k}v_{k}||^{2}_{k}}\rightarrow
0$ as $k\rightarrow\infty.$ It is a contradiction. This completes
the proof of this Lemma. \hfill$\Box$

\bigskip

\begin{lemma}\label{chengdudud}
There exists $k_{0}>0$ such that for $k\geq k_{0}$, $\dim
\Lambda_{k}=\displaystyle\sum^{n}_{i=1}l_{i}.$
\end{lemma}
\noindent{\bf Proof.} By the definition of $\Lambda_{k}$, we know
that $\dim \Lambda_{k}\leq\sum^{n}_{i=1}l_{i}.$ Note that if $i\neq
j$, then  for any $u^{k}_{i}\in N^{k}_{i},\ v_{i}^{k}\in N^{k}_{j}$
which satisfy that $ ||u^{k}_{i}||=1$ and $ ||v_{i}^{k}||=1,$
 $(u_{i}^{k},v_{i}^{k})_{k}\rightarrow 0$ as
$k\rightarrow\infty,$ since $|b^{i}_{k}-b^{j}_{k}|\rightarrow
\infty,\ i\neq j.$ Thus there exists $k_{0}>0$ such that when $k\geq
k_{0}$, $\chi_{k}e_{i,j}(\cdot+b^{i}_{k}),\ j=1,\cdots,l_{i},\
i=1,\cdots,n$ are linear independence. Hence
$\dim\Lambda_{k}\geq\sum^{n}_{i=1}l_{i}$ when $k\geq k_{0}.$ It
follows that $\dim\Lambda_{k}=\sum^{n}_{i=1}l_{i}$ when $k\geq
k_{0}.$\hfill$\Box$

\bigskip

By Lemma \ref{chengdudud}, we can define an equivalent norm on
$\Lambda_{k}$ (and $N^{k}_{i},\ \widetilde{N}^{k}_{i},\
\widehat{N}^{k}_{i},\ \widetilde{\Lambda}_{k}$ as well)   by
$$|||h|||=\sqrt{\sum^{n}_{i=1}\sum_{j=1}^{l_{i}}
x^{2}_{i,j}},$$ where  $h=\sum^{n}_{i=1}\sum_{j=1}^{l_{i}}
x_{i,j}(\chi_{k}e_{i,j})(\cdot+b^{i}_{k})\in \Lambda_{k}.$ In the
left part of this section and the next section, we use $|||\cdot|||$
as the norm of $\Lambda_{k}$ (and  $N^{k}_{i},\
\widetilde{N}^{k}_{i},\ \widehat{N}^{k}_{i},\
\widetilde{\Lambda}_{k}$ as well).

\begin{lemma}\label{taolichunfen}
There exist $k_{0}>0,\ \delta_{0}>0$ and $\tau_{0}>0$ which are
independent of $k$ such that when $k\geq k_{0},$ the following five
statements hold:
\begin{description}
\item{(1).} There exists a $C^{1}-$mapping
$$w_{k}:\overline{B_{\Lambda_{k}}(0,\delta_{0})}\rightarrow\overline{B_{\Pi_{k}}(0,\tau_{0})}$$
such that $w_{k}(0)=0$ and for any $h\in
\overline{B_{\Lambda_{k}}(0,\delta_{0})},$ $\widetilde{T}_{k}\nabla
J_{k}(u_{k}+w_{k}(h)+h)=0$, where
$\widetilde{T}_{k}:E_{k}\rightarrow \Pi_{k}$ is the orthogonal
projection.

\item{(2).} There exists a $C^{1}-$mapping
$$\widehat{w}_{k}^{i}:\overline{B_{\widehat{N}^{k}_{i}}(0,\delta_{0})}
\rightarrow\overline{B_{(\widehat{N}^{k}_{i})^{\bot}}(0,\tau_{0})}$$
such that $\widehat{w}_{k}^{i}(0)=0$ and $
P_{(\widehat{N}^{k}_{i})^{\bot}}\nabla
J_{k}(\widehat{u}^{i}_{k}+\widehat{w}_{k}^{i}(h)+h)=0$ for any $h\in
\overline{B_{\widehat{N}^{k}_{i}}(0,\delta_{0})},\ i=1,2,\cdots,n,$
where $P_{(\widehat{N}^{k}_{i})^{\bot}}:E_{k}\rightarrow
(\widehat{N}^{k}_{i})^{\bot}$ is the orthogonal projection.

\item{(3).} There exists a $C^{1}-$mapping
$$\omega^{i}:\overline{B_{N_{i}}(0,\delta_{0})}
\rightarrow\overline{B_{(N_{i})^{\bot}}(0,\tau_{0})}$$ such that
$\omega^{i}(0)=0$ and $P_{(N_{i})^{\bot}} \nabla
J(u^{i}+\omega^{i}(h)+h)=0$ for any $h\in
\overline{B_{N_{i}}(0,\delta_{0})},\ i=1,2,\cdots,n,$ where
$P_{(N_{i})^{\bot}}:X\rightarrow (N_{i})^{\bot}$ is the orthogonal
projection.

\item{(4).} There exists a $C^{1}-$mapping
$$w_{k}^{i}:\overline{B_{N^{k}_{i}}(0,\delta_{0})}
\rightarrow\overline{B_{(N^{k}_{i})^{\bot}}(0,\tau_{0})}$$ such that
$w_{k}^{i}(0)=0$ and $ P_{(N^{k}_{i})^{\bot}}\nabla
J_{k}(u^{i}_{k}+w_{k}^{i}(h)+h)=0$ for any $h\in
\overline{B_{N^{k}_{i}}(0,\delta_{0})},\ i=1,2,\cdots,n,$ where
$P_{(N^{k}_{i})^{\bot}}:E_{k}\rightarrow (N^{k}_{i})^{\bot}$ is the
orthogonal projection.

\item{(5).} There exists a $C^{1}-$mapping
$$\widetilde{w}_{k}:\overline{B_{\widetilde{\Lambda}_{k}}(0,\delta_{0})}
\rightarrow\overline{B_{\widetilde{\Pi}_{k}}(0,\tau_{0})}$$ such
that $\widetilde{w}_{k}(0)=0$ and for any $h\in
\overline{B_{\widetilde{\Lambda}_{k}}(0,\delta_{0})}$,
$$P_{\widetilde{\Pi}_{k}}\nabla
J(\sum^{n}_{i=1}u^{i}(\cdot+b^{i}_{k})+\widetilde{w}_{k}(h)+h)=0,$$
where $P_{\widetilde{\Pi}_{k}}:X\rightarrow \widetilde{\Pi}_{k}$ is
the orthogonal projection.
\end{description}
\end{lemma}
\noindent{\bf Proof.} We only give the prove of result $(1),$ since
the proofs of the other results are similar.  Set
$I_{k}(w+h)=\widetilde{T}_{k}\nabla J_{k}(u_{k}+w+h),$ then
$I_{k}(0+0)=0.$ By Lemma \ref{xiaojiazhi}, we know that there exist
$k_{0}>0$ and $\delta_{0}>0$ such that when $k\geq k_{0},$
$\widetilde{T}_{k}\nabla
I_{k}|_{\Pi_{k}}=\widetilde{T}_{k}\nabla^{2}J_{k}(u_{k}+w+h)|_{\Pi_{k}}$
is invertible if $||w+h||_{k}\leq\delta_{0}$
 and there exits $\eta>0$ such that $||(\widetilde{T}_{k}\nabla I_{k}|_{\Pi_{k}})^{-1}||\leq\eta$ for any $k\geq
 k_{0}$. Then by the implicit functional theorem, we get that there
 exist $\tau_{0}>0$ and a $C^{1}-$mapping $$w_{k}:
 \overline{B_{\Lambda_{k}}(0,\delta_{0})}\rightarrow\overline{B_{\Pi_{k}}(0,\tau_{0})}$$
such that $w_{k}(0)=0$ and $I_{k}(w_{k}(h)+h)=0.$\hfill$\Box$

\begin{remark}
Since $J_{k}$ is invariant under  the  action of
$\mathbb{Z}^{N}_{k}$, we get that
\begin{equation}\label{tapainqinshan}
(\widehat{w}^{i}_{k}(\sum^{l_{i}}_{j=1}x_{i,j}\chi_{k}e_{i,j}))(\cdot+b^{i}_{k})=
w^{i}_{k}(\sum^{l_{i}}_{j=1}x_{i,j}e^{k}_{i,j}).\end{equation}
\end{remark}

If $f:E\rightarrow F$ is a $C^{1}$ map between two Banach spaces $E$
and $F,$ we denote the derivative operator of $f$ at $u$ by $f'(u)$
and the action of $f'(u)$ on $v\in E$ is denoted by $f'(u)v$. The
proofs of the  following two Lemmas will be given in appendix.
\begin{lemma}\label{nbbzvcccfrtyujhgppp}
 For any $1\leq i\leq n,$  the following two
 statements hold:
\begin{description}
\item{(1).} As $k\rightarrow\infty,$
$$\sup\{||\widehat{w}^{i}_{k}(\sum^{l_{i}}_{j=1}x_{i,j}\chi_{k}e_{i,j})
-\omega^{i}(\sum^{l_{i}}_{j=1}x_{i,j}e_{i,j})||_{H^{1}(Q_{k})}\ :\
\sqrt{\sum^{l_{i}}_{j=1}x^{2}_{i,j}}\leq\delta_{0}\}\rightarrow 0.$$

\item{(2).} As $k\rightarrow\infty,$ for any $1\leq s\leq l_{i}$ and
$1\leq i\leq n,$
$$\sup\{||(\widehat{w}^{i}_{k})'(\sum^{l_{i}}_{j=1}x_{i,j}\chi_{k}e_{i,j})(\chi_{k}e_{i,s})
-(\omega^{i})'(\sum^{l_{i}}_{j=1}x_{i,j}e_{i,j})e_{i,s}||_{H^{1}(Q_{k})}
: \sqrt{\sum^{l_{i}}_{j=1}x^{2}_{i,j}}\leq\delta_{0}\}\rightarrow
0.$$
\end{description}
\end{lemma}

  For $h\in
\Lambda_{k},\ h=\sum^{n}_{i=1}\sum^{l_{i}}_{j=1}x_{i,j}e^{k}_{i,j},$
we denote $h_{i}=\sum^{l_{i}}_{j=1}x_{i,j}e^{k}_{i,j},\
i=1,\cdots,n.$

\begin{lemma}\label{tiandiyisa}

\begin{description}
\item{(1).} As $k\rightarrow\infty,$
$$\sup
\{||w_{k}(\sum^{n}_{i=1}h_{i})-\sum^{n}_{i=1}w^{i}_{k}(h_{i})||_{k}
\ : \ h=\sum^{n}_{i=1}h_{i}\in
\overline{B_{\Lambda_{k}}(0,\delta_{0})} \}\rightarrow 0.$$

\item{(2).} As $k\rightarrow\infty,$ for any $1\leq j\leq l_{i},\
1\leq i\leq n,$
$$\sup\{|| w'_{k}(\sum^{n}_{s=1}h_{s}) e^{k}_{i,j}
-\sum^{n}_{s=1}(w^{s}_{k})'(h_{s}) e^{k}_{i,j}||_{k}\ :\
h=\sum^{n}_{s=1}h_{s}\in\overline{B_{\Lambda_{k}}(0,\delta_{0})}\}\rightarrow
0.$$
\end{description}
\end{lemma}

\section{Critical groups of  reduction functions}

For
$x=(x_{1,1},\cdots,x_{1,l_{1}},\cdots,x_{n,1},\cdots,x_{n,l_{n}})
\in \overline{B_{\mathbb{R}^{\sum^{n}_{i=1}l_{i}}}(0,\delta_{0})}$,
we define
\begin{equation}\label{wulinweiyitenxilian7uj}
I_{k}(x)=J_{k}(u_{k}+w_{k}(\sum^{n}_{i=1}\sum^{l_{i}}_{j=1}x_{i,j}e^{k}_{i,j})+
\sum^{n}_{i=1}\sum^{l_{i}}_{j=1}x_{i,j}e^{k}_{i,j}),
\end{equation}
where
$B_{\mathbb{R}^{\sum^{n}_{i=1}l_{i}}}(0,\delta_{0})=\{x\in\mathbb{R}^{\sum^{n}_{i=1}l_{i}}
: \sqrt{\sum^{n}_{i=1}\sum^{l_{i}}_{j=1}x^{2}_{i,j}}<\delta_{0}\}.$

Denote $x_{i}=(x_{i,1},\cdots,x_{i,l_{i}})\in
\overline{B_{\mathbb{R}^{l_{i}}}(0,\delta_{0})}$ and
$$I^{i}(x_{i})=J(u^{i}+\omega^{i}(\sum^{l_{i}}_{j=1}x_{i,j}e_{i,j})+\sum^{l_{i}}_{j=1}x_{i,j}e_{i,j}),\
i=1,2,\cdots,n.$$ By (\ref{jiajiajiajnh7yhqq}),
(\ref{tapainqinshan}), Lemma \ref{tiandiyisa} and Lemma
\ref{nbbzvcccfrtyujhgppp}, we have the following Lemma

\begin{lemma}\label{cheonfqingshenlin}
As $k\rightarrow\infty,$
$||I_{k}(x)-\displaystyle\sum^{n}_{i=1}I^{i}(x_{i})||
_{C^{1}(\overline{B_{\mathbb{R}^{\sum^{n}_{i=1}l_{i}}}(0,\delta_{0}/2)})}\rightarrow
0.$
\end{lemma}
\noindent{\bf Proof.} Let
$$\mathcal{I}^{i}_{k}(x_{i})=J_{k}(u^{i}_{k}+w_{k}^{i}
(\sum^{l_{i}}_{j=1}x_{i,j}e^{k}_{i,j})+\sum^{l_{i}}_{j=1}x_{i,j}e^{k}_{i,j}),\
i=1,2,\cdots,n.$$ Then by Lemma \ref{tiandiyisa}, we have
\begin{equation}\label{unjha5556qwe7765}
||I_{k}(x)-\displaystyle\sum^{n}_{i=1}\mathcal{I}^{i}_{k}(x_{i})||
_{C^{1}(\overline{B_{\mathbb{R}^{\sum^{n}_{i=1}l_{i}}}(0,\delta_{0}/2)})}\rightarrow
0,\ k\rightarrow\infty.\end{equation} Let
$$\widehat{I}^{i}_{k}(x_{i})=J_{k}(\widehat{u}^{i}_{k}+\widehat{w}_{k}^{i}
(\sum^{l_{i}}_{j=1}x_{i,j}\chi_{k}e_{i,j})+\sum^{l_{i}}_{j=1}x_{i,j}\chi_{k}e_{i,j}),\
i=1,2,\cdots,n.$$ Then by Lemma \ref{nbbzvcccfrtyujhgppp}, we get
\begin{equation}\label{unjha5556qwe7765tgf}
||\widehat{I}^{i}_{k}-I^{i}||
_{C^{1}(\overline{B_{\mathbb{R}^{l_{i}}}(0,\delta_{0}/2)})}\rightarrow
0,\ k\rightarrow\infty,\ i=1,\cdots,n.\end{equation} By the
invariance of the $\mathbb{Z}_{k}^{N}$ action on $J_{k} $ and
(\ref{tapainqinshan}), we have
\begin{equation}\label{uuuuu555130008}
\mathcal{I}^{i}_{k}(x_{i})\equiv \widehat{I}^{i}_{k}(x_{i}),\
i=1,2,\cdots,n.
\end{equation}
By $(\ref{unjha5556qwe7765})-(\ref{uuuuu555130008})$ we get the
result of this Lemma. \hfill$\Box$

\bigskip

By the properties of $w_{k}$ and $\omega^{i},$ we have the following
Lemma
\begin{lemma}\label{chenjinchanhai}
\begin{description}
\item{(1).} If $x^{0}\in
B_{\mathbb{R}^{\sum^{n}_{i=1}l_{i}}}(0,\delta_{0}/2)$ is a critical
point of $I_{k},$ then
$$u_{k}+w_{k}(\sum^{n}_{i=1}\sum^{l_{i}}_{j=1}x^{0}_{i,j}e^{k}_{i,j})
+\sum^{n}_{i=1}\sum^{l_{i}}_{j=1}x^{0}_{i,j}e^{k}_{i,j}$$ is a
critical point of $J_{k}.$

\item{(2).} If $x^{0}\in
B_{\mathbb{R}^{\sum^{n}_{i=1}l_{i}}}(0,\delta_{0}/2)$ is a critical
point of $I^{i},$ then
$$u^{i}+\omega^{i}(\sum^{l_{i}}_{j=1}x^{0}_{i,j}e_{i,j})
+\sum^{l_{i}}_{j=1}x^{0}_{i,j}e_{i,j}$$ is a critical point of $J.$
\end{description}
\end{lemma}
\noindent{\bf Proof.} We only give the proof of result $(1)$, since
the proof of result $(2)$ is similar. If $x^{0}\in
B_{\mathbb{R}^{\sum^{n}_{i=1}l_{i}}}(0,\delta_{0}/2)$ is a critical
point of $I_{k},$ then for any $x_{s,t},$
\begin{eqnarray}\label{jja89u99kqn8jiiu44444t}
0&=&\frac{\partial I_{k}(x^{0})}{\partial x_{s,t}}\nonumber\\
&=&\left(\nabla
J_{k}(u_{k}+w_{k}(\sum^{n}_{i=1}\sum^{l_{i}}_{j=1}x^{0}_{i,j}e^{k}_{i,j})
+\sum^{n}_{i=1}\sum^{l_{i}}_{j=1}x^{0}_{i,j}e^{k}_{i,j}),\
w'_{k}(\sum^{n}_{i=1}\sum^{l_{i}}_{j=1}x^{0}_{i,j}e^{k}_{i,j})e^{k}_{s,t}\right)_{k}\nonumber\\
&&+\left(\nabla
J_{k}(u_{k}+w_{k}(\sum^{n}_{i=1}\sum^{l_{i}}_{j=1}x^{0}_{i,j}e^{k}_{i,j})
+\sum^{n}_{i=1}\sum^{l_{i}}_{j=1}x^{0}_{i,j}e^{k}_{i,j}),\
e^{k}_{s,t}\right)_{k}
\end{eqnarray}
Since
$w'_{k}(\sum^{n}_{i=1}\sum^{l_{i}}_{j=1}x^{0}_{i,j}e^{k}_{i,j})e^{k}_{s,t}\in
\Pi_{k},$ by the result $(1)$ of Lemma \ref{taolichunfen}, we get
that
\begin{equation}\label{kkmag5555555}
\left(\nabla
J_{k}(u_{k}+w_{k}(\sum^{n}_{i=1}\sum^{l_{i}}_{j=1}x^{0}_{i,j}e^{k}_{i,j})
+\sum^{n}_{i=1}\sum^{l_{i}}_{j=1}x^{0}_{i,j}e^{k}_{i,j}),\
w'_{k}(\sum^{n}_{i=1}\sum^{l_{i}}_{j=1}x^{0}_{i,j}e^{k}_{i,j})e^{k}_{s,t}\right)_{k}=0.
\end{equation}
By (\ref{jja89u99kqn8jiiu44444t}) and (\ref{kkmag5555555}), we get
that for  any $1\leq t\leq l_{s},$ $1\leq s\leq n,$
\begin{equation}\label{jiaxianrenmin7y6}
\left(\nabla
J_{k}(u_{k}+w_{k}(\sum^{n}_{i=1}\sum^{l_{i}}_{j=1}x^{0}_{i,j}e^{k}_{i,j})
+\sum^{n}_{i=1}\sum^{l_{i}}_{j=1}x^{0}_{i,j}e^{k}_{i,j}),\
e^{k}_{s,t}\right)_{k}=0.
\end{equation}
By (\ref{jiaxianrenmin7y6}) and the result (1) of Lemma
\ref{taolichunfen}, we know that
$$u_{k}+w_{k}(\sum^{n}_{i=1}\sum^{l_{i}}_{j=1}x^{0}_{i,j}e^{k}_{i,j})
+\sum^{n}_{i=1}\sum^{l_{i}}_{j=1}x^{0}_{i,j}e^{k}_{i,j}$$ is a
critical point of $J_{k}.$ \hfill$\Box$

\begin{remark}\label{xiaojiji58888uh}
By the condition {\bf (*)} and Lemma \ref{chenjinchanhai}, we know
that $0$ is the unique critical point of $I^{i}(x_{i})$ in
$\overline{B_{\mathbb{R}^{l_{i}}}(0,\delta_{0})},\ i=1,\cdots,n.$
\end{remark}
By Lemma \ref{xiaojiazhi}, we know that $\widetilde{T}_{k}\circ
\nabla^{2}J_{k}(u_{k})$ is a bounded, invertible and self-adjoint
operator in Hilbert space $\Pi_{k}.$ Let $P^{+}_{k}$  (resp. $
P^{-}_{k}$) be the orthogonal projection from $E_{k}$ into the
positive (resp. negative) subspace $\Pi_{k}^{+}$ (resp. $
\Pi^{-}_{k}$) with respect to the spectral decomposition of
$\widetilde{T}_{k}\circ \nabla^{2}J_{k}(u_{k}).$ By Lemma 1 of
\cite{gm}, we have the following Lemma:

\begin{lemma}\label{chuxiaoshengyan67ytgfff}
 For any $u\in  \overline{B_{E_{k}}(u_{k},\delta_{0}/2)}$,
 $u$ has the unique decomposition $u=u_{k}+w+h$ where
 $w\in \Pi_{k}$ and $h\in \Lambda_{k}$.
  There exist a diffeomorphism $\Psi_{k}: \overline{B_{E_{k}}(u_{k},\delta_{0}/2)}
 \rightarrow E_{k}$ which satisfies that $ \Psi(u_{k})=u_{k}$  such that for any
 $u=u_{k}+w+\sum^{n}_{i=1}\sum^{l_{i}}_{j=1}x_{i,j}e^{k}_{i,j}\in \overline{B_{E_{k}}(u_{k},\delta_{0}/2)},$
 $$J_{k}(\Psi_{k}(u))=||P^{+}_{k}w||^{2}_{k}-||P^{-}_{k}w||^{2}_{k}
+J_{k}(u_{k}+w_{k}(\sum^{n}_{i=1}\sum^{l_{i}}_{j=1}x_{i,j}e^{k}_{i,j})+
\sum^{n}_{i=1}\sum^{l_{i}}_{j=1}x_{i,j}e^{k}_{i,j}).$$
\end{lemma}
\begin{remark}\label{haenkanchenhlinche9i}
By Lemma \ref{shanshan}, we know that if $k$ large enough, then
$K^{(i_{k})}_{k}\subset B_{E_{k}}(u_{k}, \delta_{0}/2).$ Let
$\mathcal{K}_{k}:=\{x\in\overline{B_{\mathbb{R}^{\sum^{n}_{i=1}l_{i}}}(0,\delta_{0}/2)}\
:\ \mbox{there exists}\ u\in K^{(i_{k})}_{k},\ \mbox{such that}\
P_{\Lambda_{k}}(u-u_{k})=\sum^{n}_{i=1}\sum^{l_{i}}_{j=1}x_{i,j}e^{k}_{i,j}\},$
where $P_{\Lambda_{k}}:E_{k}\rightarrow \Lambda_{k}$ is the
orthogonal projection. Then by Lemma \ref{chenjinchanhai},
$\mathcal{K}_{k}$ is the critical set of $I_{k}$ in
$\overline{B_{\mathbb{R}^{\sum^{n}_{i=1}l_{i}}}(0,\delta_{0}/2)}$.
By Lemma \ref{shanshan}, we know that
$\mbox{diam}(\mathcal{K}_{k})\rightarrow 0$ as $k\rightarrow\infty.$
\end{remark}

\begin{lemma}\label{dadaoqianxianquiz}
If $k$ large enough, then there exists integer $\widehat{m}_{k}$
which satisfies that $0\leq \widehat{m}_{k}\leq \sum^{n}_{i=1}l_{i}$
such that $C_{\widehat{m}_{k}}(I_{k}, \mathcal{K}_{k})\neq 0.$
\end{lemma}
\noindent{\bf Proof.}
 Let $(W_{1},W^{-}_{1})$ and $(W_{2},W_{2}^{-})$
be the GM-pairs for the functionals
$\mathcal{J}_{k}(w+h)=||w||^{2}_{k}-||h||^{2}_{k}$ in the unique
critical point $0$ and $I_{k}$ in the critical set $\mathcal{K}_{k}$
respectively. Then by \cite[Lemma 5.1]{cha}, we know that
$(W_{1}\times W_{2},(W_{1}^{-}\times W_{2})\cup(W_{1}\times
W^{-}_{2}))$ is a GM-pair for the isolated critical set
$\Psi^{-1}_{k}(K^{(i_{k})}_{k})$ of the functional $J_{k}\circ
\Psi_{k}.$ Then by Theorem 5.5 of \cite{cha}, we get that
\begin{equation}\label{huanghedahechangju87y}
C_{*}(J_{k}\circ\Psi_{k},\Psi^{-1}_{k}(K^{(i_{k})}_{k}))=C_{*}(\mathcal{J}_{k},0)\otimes
C_{*}(I_{k},\mathcal{K}_{k}).
\end{equation}
Since
$C_{q}(\mathcal{J}_{k},0)=\delta_{q,\dim\widetilde{\Pi}^{-}_{k}}\mathbb{Z}_{2}$,
by (\ref{huanghedahechangju87y}), we get that
\begin{equation}\label{dadaorunndhfgrt77herfmjk}
C_{q}(I_{k},\mathcal{K}_{k})
=C_{q+\dim\widetilde{\Pi}^{-}_{k}}(J_{k}\circ\Psi_{k},\Psi^{-1}_{k}(K^{(i_{k})}_{k})).
\end{equation}
Since $C_{*}(J_{k}\circ\Psi_{k},\Psi^{-1}_{k}(K^{(i_{k})}_{k})\cong
C_{*}(J_{k}, K^{(i_{k})}_{k})$, by (\ref{dadaorunndhfgrt77herfmjk}),
we get that
\begin{equation}\label{yisoufuyinzuocajnh9i}
C_{q}(I_{k},\mathcal{K}_{k})\cong
C_{q+\dim\widetilde{\Pi}^{-}_{k}}(J_{k}, K^{(i_{k})}_{k}).
\end{equation}
By Lemma \ref{hhhnayooahj88765kk} and (\ref{yisoufuyinzuocajnh9i}),
we have
$C_{j(k)-\dim\widetilde{\Pi}^{-}_{k}+1}(I_{k},\mathcal{K}_{k})\cong
C_{j(k)+1}(J_{k}, K^{(i_{k})}_{k})\neq 0.$ Let
$\widehat{m}_{k}=j(k)-\dim\widetilde{\Pi}^{-}_{k}+1.$ Notice that
$I_{k}$ is a function defined in some subset of
$\mathbb{R}^{\sum^{n}_{i=1}l_{i}}$, we get that
$C_{q}(I_{k},\mathcal{K}_{k})=0$ if $q>\sum^{n}_{i=1}l_{i}.$ Thus
$0\leq \widehat{m}_{k}\leq \sum^{n}_{i=1}l_{i}.$ \hfill$\Box$

\begin{lemma}\label{aijianshanbuai}
Let $(W,W_{-})$ be a GM-pair of the isolated critical point $0$  of
$\displaystyle\sum^{n}_{i=1}I^{i}(x_{i})$ with respect to
$-d(\sum^{n}_{i=1}I^{i}(x_{i}))$ and $W\subset
B_{\mathbb{R}^{\sum^{n}_{i=1}l_{i}}}(0,\delta_{0}/2). $ If $k$ large
enough, then  $(W,W_{-})$ is also a GM-pair of the isolated critical
set $\mathcal{K}_{k}$ of $I_{k}$ with respect to certain
pseudo-gradient vector field  of $I_{k}$.
\end{lemma}
\noindent{\bf Proof.} Since $\mbox{diam}(\mathcal{K}_{k})\rightarrow
0,$ as $k\rightarrow\infty,$ we know that there exists $ r>0$ such
that if $k$ large enough, then
$B_{\mathbb{R}^{\sum^{n}_{i=1}l_{i}}}(0,r)\subset \mbox{int}(W)$,
the interior of $W$,  and $\mathcal{K}_{k}\subset
B_{\mathbb{R}^{\sum^{n}_{i=1}l_{i}}}(0,r/4)$. Note that
$$\beta:=\inf\{||d(\sum^{n}_{i=1}I^{i}(x_{i}))||\ :\ x\in W\setminus
B_{\mathbb{R}^{\sum^{n}_{i=1}l_{i}}}(0,r/2)\}>0.$$ Define $\rho\in
C^{2}(\mathbb{R}^{\sum^{n}_{i=1}l_{i}}, \mathbb{R})$ satisfying
$$\rho(x)=\left\{
\begin{array}{l}
1, \quad x\in B_{\mathbb{R}^{\sum^{n}_{i=1}l_{i}}}(0,r/2)\\
0, \quad x\not\in B_{\mathbb{R}^{\sum^{n}_{i=1}l_{i}}}(0,r),\\
\end{array} \right.$$
with $0\leq\rho(x)\leq 1$ and a vector field
$$V(x)=\frac{3}{2}(\rho(x)dI_{k}(x)+(1-\rho(x))d(\sum^{n}_{i=1}I^{i}(x_{i}))).$$
Choosing $0<\epsilon<\beta/4$, by Lemma \ref{cheonfqingshenlin}, we
know that if $k$ large enough, then
\begin{equation}\label{7771776156hhyt}
||I_{k}(x)-\sum^{k}_{i=1}I^{i}(x_{i})||
_{C^{1}(\overline{B_{\mathbb{R}^{\sum^{n}_{i=1}l_{i}}}(0,\delta_{0})})}<\epsilon.
\end{equation}
We shall prove that $||V(x)||\leq 2||dI_{k}(x)||$ and
$(V(x),dI_{k}(x))\geq ||d(I_{k}(x))||^{2}$, where $(\cdot,\cdot)$
denotes the inner product in $\mathbb{R}^{\sum^{n}_{i=1}l_{i}}.$

By (\ref{7771776156hhyt}), we know that for any  $ x\not\in
B_{\mathbb{R}^{\sum^{n}_{i=1}l_{i}}}(0,r/2),$
\begin{equation}\label{7hbsv5refs4f}||dI_{k}(x)||\geq
||d(\sum^{n}_{i=1}I^{i}(x_{i}))||-\epsilon\geq\beta-\epsilon>3\epsilon.\end{equation}

Thus we have
\begin{eqnarray}\label{damenxianjue}
&&(V(x),dI_{k}(x))\nonumber\\
&=&
\left(\frac{3}{2}(\rho(x)dI_{k}(x)+(1-\rho(x))d(\sum^{n}_{i=1}I^{i}(x_{i}))),\
dI_{k}(x)\right)\nonumber\\
&=&\left(\frac{3}{2}(dI_{k}(x)
+(1-\rho(x))(d(\sum^{n}_{i=1}I^{i}(x_{i}))-dI_{k}(x))),\
dI_{k}(x)\right)\nonumber\\
&\geq&\frac{3}{2}(||dI_{k}(x)||^{2}-\epsilon||dI_{k}(x)||)\quad \quad (\mbox{by}\ (\ref{7771776156hhyt}))\nonumber\\
&\geq&\frac{3}{2}(||dI_{k}(x)||^{2}-\frac{1}{3}||dI_{k}(x)||^{2})\nonumber\\
&=&||dI_{k}(x)||^{2}\nonumber
\end{eqnarray}
and
\begin{eqnarray}\label{qianjianmencao}
||V(x)||&=&||\frac{3}{2}((1-\rho(x))d(\sum^{n}_{i=1}I^{i}(x_{i}))+\rho(x)dI_{k}(x))||\nonumber\\
&\leq& \frac{3}{2}(||dI_{k}(x)||+\epsilon)\nonumber\\
&\leq& 2||dI_{k}(x)||\quad\quad (\mbox{by}\
(\ref{7hbsv5refs4f})).\nonumber
\end{eqnarray}
Since for any $ x\in B_{\mathbb{R}^{\sum^{n}_{i=1}l_{i}}}(0,r/2),\
V(x)=dI_{k}(x)$, the verification is trivial.

Notice that $V(x)=-d(\sum^{n}_{i=1}I^{i}(x_{i}))$ outside a ball
$B_{\mathbb{R}^{\sum^{n}_{i=1}l_{i}}}(0,r)\subset \mbox{int}(W).$ It
is not difficult to verify that $(W,W_{-})$ is a GM-pair of
$I_{k}(x)$ with respect to $V$.\hfill $\Box$

\begin{lemma}\label{wanjinrenzhongmie}
There exist an index $i_{0}$ satisfying  $1\leq i_{0}\leq n$ and a
nonnegative integer $m_{0}$  satisfying  $0\leq
m_{0}\leq\sum^{n}_{i=1}l_{i}$ such that $C_{m_{0}}(I^{i_{0}},0)\neq
0.$
\end{lemma}
 \noindent{\bf Proof.} Let $(W,W^{-})$ be  a GM-pair of $\sum^{n}_{i=1}I^{i}(x_{i})$ for
 the isolated critical
 point $0$ which satisfies
 that $W\subset
 B_{\mathbb{R}^{\sum^{n}_{i=1}l_{i}}}(0,\delta_{0}/4).$ Then by
 Lemma \ref{aijianshanbuai}, we know that $(W,W^{-})$ is also a
 GM-pair of $I_{k}$ for the  isolated critical set $\mathcal{K}_{k}$
 if $k$ large enough.
  Thus by Lemma \ref{wauxianfenguanzai} and Lemma \ref{dadaoqianxianquiz}, we get that
 \begin{equation}\label{tapianqingshanerenewihgs5frk}
 C_{\widehat{m}_{k}}(\sum^{n}_{i=1}I^{i}(x_{i}),0)=H_{\widehat{m}_{k}}(W,W_{-})
 =C_{\widehat{m}_{k}}(I_{k},\mathcal{K}_{k})\neq 0.
 \end{equation}
By  \cite[Theorem 5.5]{cha}, we get that
\begin{equation}\label{jiannanhao6tgpqbbc}
C_{*}(\sum^{n}_{i=1}I^{i}(x_{i}),0)=\bigotimes^{n}_{i=1}C_{*}(I^{i},0).
\end{equation}
  The  result of this Lemma follows from
(\ref{tapianqingshanerenewihgs5frk}) and (\ref{jiannanhao6tgpqbbc}).
 \hfill$\Box$

\section{Proof of Theorem \ref{theorem1}}

\hspace*{\parindent}If the condition {\bf (*)} does not hold, then
Equation (\ref{1.1}) has infinitely many geometrically different
solutions and the proof terminates. In the following, we always
assume that the condition {\bf (*)} holds.

We denote $u^{i_{0}}$ by $u_{0}$, where $i_{0}$ is the index
appeared in Lemma \ref{wanjinrenzhongmie}. By the condition
{\bf{(*)}}, we know that $u_{0}$ is an isolated critical point of
$J.$ We will use $u_{0}$ as a basic ``one-bump" solution to
construct multi-bump solutions for equation (\ref{1.1}).

For positive integer $m\geq 2$ and $b_{i}\in \mathbb{Z}^{N},\
i=1,2,\cdots,m$, denote $u_{b_{i}}=u_{0}(x-b_{i}),\ i=1,2,\cdots,m.$
 For $\mathbf{k}=(b_{1},\cdots,b_{m})\in (\mathbb{Z}^{N})^{m}$,
 denote $u_{\mathbf{k}}=\sum^{m}_{i=1}u_{b_{i}}$
 and let
$l_{\mathbf{k}}=\min\{|b_{i}-b_{j}|\ :\ i\neq j\}$. Let
$\mathcal{N}$ be the kernel of $\nabla^{2}J(u_{0})$ and let
$\mathcal{N}=\mbox{span}\{e_{1},\cdots,e_{l}\}$, where  $l=\dim
\mathcal{N}$ and $e_{i},\ i=1,\cdots,l$ satisfy
$(e_{i},e_{j})_{X}=\delta_{i,j}$.   Let $\mathcal{N}_{b_{i}}$ be the
kernel of $\nabla^{2}J(u_{b_{i}}),\ i=1,2,\cdots,m,$ then
$\mathcal{N}_{b_{i}}=\mbox{span}\{e_{1}(\cdot-b_{i}),\cdots,e_{l}(\cdot-b_{i})\}.$
Let $\mathcal{N}_{\mathbf{k}}=\mbox{span}\{\mathcal{N}_{b_{i}}\ :\
i=1,2,\cdots,m\}$ and
$Z_{\mathbf{k}}=(\mathcal{N}_{\mathbf{k}})^{\bot}\subset X$.

\medskip

As the same argument as Lemma \ref{taolichunfen}
 and Lemma \ref{chenjinchanhai}, we have the
following two Lemmas:
\begin{lemma}\label{dongfensongmuanru} There exist
$\delta>0$ and a $C^{1}-$mapping $\omega:
\overline{B_{\mathbb{R}^{l}}(0,\delta)}\rightarrow
\mathcal{N}^{\bot}$ such that
\begin{description}
\item{(1).} $P_{\mathcal{N}^{\bot}}\nabla
J(u_{0}+\omega(x_{1},\cdots,x_{l})+\displaystyle\sum^{l}_{i=1}x_{i}e_{i})=0$
and $\omega(0)=0,$ where $P_{\mathcal{N}^{\bot}}:X\rightarrow
\mathcal{N}^{\bot}$ is the orthogonal projection.

\item{(2).} If $(x^{0}_{1},\cdots,x^{0}_{l})$ is a critical point of
$J(u_{0}+\omega(x_{1},\cdots,x_{l})+\sum^{l}_{i=1}x_{i}e_{i}),$ then
$u_{0}+\omega(x_{1},\cdots,x_{l})+\sum^{l}_{i=1}x_{i}e_{i}$ is a
critical point of $J.$
\end{description}
\end{lemma}
\begin{lemma}\label{chuangzhangdadadase} There exist
$L>0$ and $\delta>0$ such that if $l_{\mathbf{k}}>L,$ then there is
a $C^{1}-$mapping
$w_{\mathbf{k}}:\overline{B_{\mathbb{R}^{ml}}(0,\delta)}\rightarrow
Z_{\mathbf{k}}$ satisfying that
\begin{description}
\item{(1).} $P_{Z_{\mathbf{k}}}\nabla J(u_{\mathbf{k}}+w_{\mathbf{k}}
(x_{1,1},\cdots,x_{1,l},\cdots,x_{m,1},\cdots,x_{m,l})
+\displaystyle\sum^{m}_{i=1}\sum^{l}_{j=1}x_{i,j}e_{j}(\cdot-b_{i}))=0.$

\item{(2).} if $(x^{0}_{1,1},\cdots,x^{0}_{m,l})$ is a critical point of
$$J(u_{\mathbf{k}}+w_{\mathbf{k}}
(x_{1,1},\cdots,x_{1,l},\cdots,x_{m,1},\cdots,x_{m,l})
+\sum^{m}_{i=1}\sum^{l}_{j=1}x_{i,j}e_{j}(\cdot-b_{i}))$$ then
$u_{\mathbf{k}}+w_{\mathbf{k}}
(x^{0}_{1,1},\cdots,x^{0}_{1,l},\cdots,x^{0}_{m,1},\cdots,x^{0}_{m,l})
+\sum^{m}_{i=1}\sum^{l}_{j=1}x^{0}_{i,j}e_{j}(\cdot-b_{i})$ is a
critical point of $J$.
\end{description}
\end{lemma}

Let
$H_{\mathbf{k}}(x_{1,1},\cdots,x_{m,l})=J(u_{\mathbf{k}}+w_{\mathbf{k}}
(x_{1,1},\cdots,x_{m,l})
+\sum^{m}_{i=1}\sum^{l}_{j=1}x_{i,j}e_{j}(\cdot-b_{i}))$ and
$H(x_{1},\cdots,x_{l})=J(u_{0}+\omega(x_{1},\cdots,x_{l})+\sum^{l}_{j=1}x_{j}e_{j}).$
Let $x=(x_{1,1},\cdots,x_{m,l})$ and $
x_{i}=(x_{i,1},\cdots,x_{i,l}),\ i=1,2,\cdots,m.$

\begin{remark}\label{hanb6235999999}
By the condition {\bf (*)}, we know that $u_{0}$ is an isolated
critical point of $J$. Then by Lemma \ref{dongfensongmuanru}, we
know that $0$ is the unique critical point of $H$ in
$\overline{B_{\mathbb{R}^{l}}(0,\delta)}$.
\end{remark}

As the same argument as Lemma \ref{cheonfqingshenlin}, we have the
following Lemma:
\begin{lemma}\label{beifenyoujiaren} There exists
$\delta>0$ such that as $l_{\mathbf{k}}\rightarrow\infty,$
$$||H_{\mathbf{k}}(x)-\sum^{m}_{i=1} H(x_{i})
||_{C^{1}(\overline{B_{\mathbb{R}^{ml}}(0,\delta)})} \rightarrow
0.$$
\end{lemma}

Note that $H(x)=I^{i_{0}}(x)$ for $x\in
\overline{B_{\mathbb{R}^{l}}(0,\delta)}$. By Lemma
\ref{wanjinrenzhongmie}, we know that $C_{m_{0}}(I^{i_{0}},0)\neq
0.$ Thus we have the following Lemma:
\begin{lemma}\label{qianjiaqnchaoruleillk}
$C_{m_{0}}(H,0)\neq 0.$
\end{lemma}

By Lemma \ref{qianjiaqnchaoruleillk} and $\displaystyle
C_{mm_{0}}(\sum^{m}_{i=1}H(x_{i}),0)=\displaystyle\bigotimes^{m}_{i=1}C_{m_{0}}(H(x_{i}),0)$
(see (\ref{jiannanhao6tgpqbbc}) for reference), we get that the
following Lemma:

\begin{lemma}\label{yejinchunsakykonf}
$\displaystyle C_{mm_{0}}(\sum^{m}_{i=1}H(x_{i}),0)\neq 0.$
\end{lemma}

By Lemma \ref{yejinchunsakykonf} and Lemma \ref{beifenyoujiaren}, we
have the following Lemma
\begin{lemma}\label{vxcxsdsdsdsdsdsdd} If
$l_{\mathbf{k}}$ large enough,  then $H_{\mathbf{k}}$ has at least a
critical point $x^{\mathbf{k}}\in B_{\mathbb{R}^{ml}}(0,\delta)$
which satisfies that $x^{\mathbf{k}}\rightarrow 0$ as
$l_{\mathbf{k}}\rightarrow\infty.$
\end{lemma}
\noindent{\bf Proof.} Let $(W,W_{-})$ be a GM-pair of the isolated
critical point $0$ of $\sum^{m}_{i=1}H(x_{i})$ with respective to
the gradient field $-d(\sum^{m}_{i=1}H(x_{i}))$. By Lemma
\ref{beifenyoujiaren} and the proof of Lemma \ref{aijianshanbuai},
we know that when $l_{\mathbf{k}}$ large enough, if $\eta$ is the
flow generalized by the following ordinary differential equation
$$\dot{\eta}= V_{1}(\eta(x,t)),\ \eta(x,0)=x,$$
where $V_{1}(x)=g(x)\frac{V(x)}{||V(x)||}$,
$V(x)=\frac{3}{2}(\rho(x)dH_{\mathbf{k}}(x)+(1-\rho(x))d(\sum^{n}_{i=1}H(x_{i}))),$
 $\rho\in C^{2}(\mathbb{R}^{\sum^{n}_{i=1}l_{i}},
\mathbb{R})$ satisfying $0\leq\rho(x)\leq 1$ for any $x$ and
$$\rho(x)=\left\{
\begin{array}{l}
1, \quad x\in B_{\mathbb{R}^{ml}}(0,\delta/4)\\
0, \quad x\not\in B_{\mathbb{R}^{ml}}(0,\delta/2),\\
\end{array} \right.$$ and
$g(x)=\left\{
\begin{array}{l}
\min\{\mbox{dist}(x,K_{H_{\mathbf{k}}}), 1\}, \ \mbox{if}\ K_{H_{\mathbf{k}}}\neq\emptyset\\
1, \quad\quad\quad\quad\quad\quad\quad\quad\quad \mbox{if}\ K_{H_{\mathbf{k}}}=\emptyset,\\
\end{array} \right.$ then $(W,W_{-})$ satisfies
the following conditions:
\begin{description}
\item{(1).}  $W$ has the (MVP) property with respect to the flow $\eta.$

\item{(2).} $W_{-}$ is an exit set for $W,$ i.e., for each $x_{0}\in W$
and $t_{1}>0$ such that $\eta(x_{0},t_{1})\not\in W,$ there exists
$t_{0}\in [0,t_{1})$ such that $\eta(x_{0},[0,t_{0}])\subset W$ and
$\eta(x_{0},t_{0})\in W_{-}.$

\item{(3).} $W_{-}$ is closed and is a union of a finite number of
sub-manifolds that transversal to the flow $\eta.$
\end{description}

If $H_{\mathbf{k}}$ has no critical point in $W,$ then for any $x\in
B_{\mathbb{R}^{ml}}(0,\delta)$
\begin{equation}\label{nn7hhqp098}
g(x)\geq\iota\ \mbox{and}\ ||V(x)||\geq\iota \end{equation} for some
$\iota>0$. Thus for any $x\in W$, there exists $t\geq 0$ such that
$\eta(x,t)\not\in W.$ In fact, if there exists $x\in W$ such that
for any $t\geq 0$, $\eta(x,t)\in W,$ then by the fact that $W$ is a
bounded closed set in finite dimension space $\mathbb{R}^{ml}$, we
deduce that there exists a sequence $\{t_{n}\}$ such that
$t_{n}\rightarrow +\infty$ as $n\rightarrow\infty$ and
$\eta(x,t_{n})$ converges to some point $x_{0}\in W.$ Then $x_{0}$
must  satisfy that $V_{1}(x_{0})=0.$ It contradicts to
(\ref{nn7hhqp098}). Thus for any $x\in W,$ $t_{x}=\inf\{t'\geq 0\ :\
\eta(x,t')\in W_{-}\}<+\infty.$ It is easy to verify that $t_{x}=0$
for any $x\in W_{-}.$ Define $\mathcal{H}:W\times [0,1]\rightarrow
W,$ $\mathcal{H}(x,s)=\eta(x,st_{x}).$ It follows that $W_{-}$ is a
deformation retract of $W.$ Thus $H_{q}(W,W_{-})=0,\ \forall q.$ But
$(W,W_{-})$ is a GM-pair of the isolated critical point $0$ of
$\sum^{m}_{i=1}H(x_{i})$, by Lemma \ref{wauxianfenguanzai} and Lemma
\ref{yejinchunsakykonf}, we get that
$$H_{mm_{0}}(W,W_{-})\cong C_{mm_{0}}(\sum^{m}_{i=1}H(x_{i}),0)\neq 0.$$
It is a contradiction. Thus if $l_{\mathbf{k}}$ large enough, then
$H_{\mathbf{k}}$ has at least a critical point $x^{\mathbf{k}}$ in
$B_{\mathbb{R}^{ml}}(0,\delta)$. Finally, by Lemma
\ref{beifenyoujiaren} and Remark \ref{hanb6235999999}, we get that
$x^{\mathbf{k}}\rightarrow 0$ as $l_{\mathbf{k}}\rightarrow\infty.$
\hfill$\Box$

\bigskip

 By Lemma \ref{vxcxsdsdsdsdsdsdd} and Lemma
\ref{chuangzhangdadadase}, we get  the following result:
\begin{Theorem}\label{dagonggaochenggg}
$u_{\mathbf{k}}+w_{\mathbf{k}}(x^{\mathbf{k}}_{1,1},\cdots,x^{\mathbf{k}}
_{m,l})+\sum^{m}_{i=1}\sum^{l}_{j=1}x^{\mathbf{k}}_{i,j}e_{j}(\cdot-b_{i})$
is a critical point of $J.$ Furthermore, as
$l_{\mathbf{k}}\rightarrow +\infty$,
$w_{\mathbf{k}}(x^{\mathbf{k}}_{1,1},\cdots,x^{\mathbf{k}}_{m,l})\rightarrow
0$ and
$\displaystyle\sum^{m}_{i=1}\sum^{l}_{j=1}x^{\mathbf{k}}_{i,j}e_{j}(\cdot-b_{i})\rightarrow
0.$
\end{Theorem}

\section{Appendix}
\hspace*{\parindent} In this Section, we shall give the proofs of
Lemma \ref{mnnnssshgdfr}, Lemma \ref{lantaojin}, Lemma
\ref{nbbzvcccfrtyujhgppp} and Lemma \ref{tiandiyisa}. The proof of
Lemma \ref{yadanhexiawa} is similar to the proof of Lemma
\ref{mnnnssshgdfr} and the proof of Lemma \ref{tongrenttthqqqw} is
similar to the proof of  Lemma \ref{lantaojin}.

\bigskip

\noindent{\bf \large Proof of Lemma \ref{mnnnssshgdfr}:}

\medskip

 For convenience, we set
$$v_{k}=\sum^{n}_{i=1}u^{i}(\cdot+b^{i}_{k})-u_{k},$$
by (\ref{jiajiajiajnh7yhqq}), we have
\begin{equation}\label{xiyanwuxianhao765}
||v_{k}||_{H^{1}(Q_{k})}\rightarrow 0\ \mbox{ as}\
k\rightarrow\infty.
\end{equation}
Let $$ \Omega^{0}_{\epsilon,k}=\{x\in Q_{k}\ :\ |u_{k}(x)|\leq
1/\epsilon\},\ \Omega^{1}_{\epsilon,k}=\{x\in Q_{k}\ :\
|u_{k}(x)|> 1/\epsilon\},$$ and
$$U^{0}_{\epsilon,k}=\{x\in Q_{k}\ : \
|v_{k}(x)|< \epsilon\},\ U^{1}_{\epsilon,k}=\{x\in Q_{k}\ : \
|v_{k}(x)|\geq\epsilon\}.$$  By
$||v_{k}||_{H^{1}(Q_{k})}\rightarrow 0$ as $k\rightarrow\infty$,
we get that for every $\epsilon>0$
\begin{equation}\label{nnna7yhhab444}
\mbox{mes}(U^{1}_{\epsilon,k})\rightarrow 0\ \mbox{as}\
k\rightarrow\infty\end{equation} and by the fact that $\{u_{k}\}$
is bounded in $H^{1}(Q_{k})$, we get that as $ \epsilon\rightarrow
0$,
\begin{equation}\label{hhhhh777y76}
\mbox{mes}(\Omega^{1}_{\epsilon,k})\rightarrow
0,\end{equation}holds uniformly for $k\in\mathbb{N},$ here
$\mbox{mes}(A)$ denotes the Lebesgue measure of the set $A$.

By the definition of $v_{k}$ and $\Omega^{i}_{\epsilon,k}$ and
$U^{i}_{\epsilon,k},\ i=0,1$, we have
\begin{eqnarray}\label{jianjianren}
&&\int_{Q_{k}}|f'(x,\sum^{n}_{i=1}u^{i}(\cdot+b^{i}_{k}))-f'(x,u_{k})|\cdot
|\psi_{k}|\cdot|\varphi_{k}|dx\nonumber\\
&=&\int_{Q_{k}}|f'(x,u_{k}+v_{k})-f'(x,u_{k})|\cdot
|\psi_{k}|\cdot|\varphi_{k}|dx\nonumber\\
&\leq&\int_{
U^{0}_{\epsilon,k}\cap\Omega^{0}_{\epsilon,k}}+\int_{\Omega^{1}_{\epsilon,k}}
+\int_{U^{1}_{\epsilon,k}}.
\end{eqnarray}By the condition ${\bf(f_{2})}$ and H\"older inequality, we get
\begin{eqnarray}\label{xiaopp}
&&\int_{U^{1}_{\epsilon,k}}|f'(x,u_{k}+v_{k})-f'(x,u_{k})|\cdot
|\varphi_{k}|\cdot|\psi_{k}|dx\nonumber\\
&\leq& C\int_{U^{1}_{\epsilon,k}}|u_{k}+v_{k}|^{q-2}\cdot
|\varphi_{k}|\cdot|\psi_{k}|dx+C\int_{U^{1}_{\epsilon,k}}|u_{k}+v_{k}|^{p-2}\cdot
|\varphi_{k}|\cdot|\psi_{k}|dx\nonumber\\
&&+C\int_{U^{1}_{\epsilon,k}}|u_{k}|^{q-2}\cdot
|\varphi_{k}|\cdot|\psi_{k}|dx+C\int_{U^{1}_{\epsilon,k}}|u_{k}|^{p-2}\cdot
|\varphi_{k}|\cdot|\psi_{k}|dx\nonumber
\end{eqnarray}
and
\begin{eqnarray}\label{xiaojiji}
&&\int_{U^{1}_{\epsilon,k}}|u_{k}+v_{k}|^{q-2}\cdot
|\varphi_{k}|\cdot|\psi_{k}|dx\nonumber\\
&\leq& (\int_{U^{1}_{\epsilon}}1^{r}dx)^{\frac{1}{r}}
\cdot(\int_{U^{1}_{\epsilon,k}}|u_{k}+v_{k}|^{q+(q-2)\delta'})^{\frac{1}{\delta'+\frac{q}{q-2}}}
\cdot (\int_{U^{1}_{\epsilon,k}}|\varphi_{k}|^{q})^{\frac{1}{q}}
\cdot(\int_{U^{1}_{\epsilon,k}}|\psi_{k}|^{q})^{\frac{1}{q}},\nonumber
\end{eqnarray}
where $q+(q-2)\delta'\leq 2^{*}$ and
$\frac{1}{r}+\frac{1}{\delta'+\frac{q}{q-2}}+\frac{2}{q}=1.$ It
follows that there exists $C'_{1}>0$ which is independent of $k,\
\psi_{k}$ and $\varphi_{k}$ such that
$$\int_{U^{1}_{\epsilon,k}}|u_{k}+v_{k}|^{q-2}|\psi_{k}\cdot\varphi_{k}|dx
\leq C'_{1}(\mbox{mes}(U^{1}_{\epsilon,k}))^{\frac{1}{r}}.$$ Thus
there exists $C_{1}>0$ which is independent of $k,\ \psi_{k}$ and
$\varphi_{k}$ such that
\begin{equation}\label{huanabv5frthhhh}
\int_{U^{1}_{\epsilon,k}}|f'(x,u_{k}+v_{k})-f'(x,u_{k})|\cdot
|\psi_{k}\cdot\varphi_{k}|dx \leq
C_{1}(\mbox{mes}(U^{1}_{\epsilon,k}))^{\frac{1}{r}}.
\end{equation}
In the same way, there exists $C_{2}>0$ which is independent of $k,\
\psi_{k}$ and $\varphi_{k}$ such that
\begin{equation}\label{shenminchenke6gttgftfrr11111}
\int_{\Omega^{1}_{\epsilon,k}}|f'(x,u_{k}+v_{k})-f'(x,u_{k})|\cdot
|\psi_{k}\cdot\varphi_{k}|dx \leq
C_{2}(\mbox{mes}(\Omega^{1}_{\epsilon,k}))^{\frac{1}{r}}.
\end{equation}

Choose $\eta\in C^{\infty}_{0}(\mathbb{R})$ which satisfies that
$0\leq\eta\leq 1,$ $|\eta'(t)|\leq 2,$
$\int_{\mathbb{R}}\eta(t)dt=1$, $\eta\equiv 1$ in $(-\frac{1}{2},
\frac{1}{2})$ and $\eta\equiv 0$ in $\mathbb{R}\setminus (-1,1).$
Let $\eta_{\delta}(t)=\frac{1}{\delta}\eta(\frac{t}{\delta})$ and
$$g_{\delta}(x,t)=\int_{\mathbb{R}}f'_{s}(x,s)\eta_{\delta}(t-s)ds.$$
Since $f'(x,t)$ is a Caratheodory function, we deduce that for
almost all $x\in\mathbb{R}^{N}$ and for all $t\in\mathbb{R},$
\begin{equation}\label{jjhsf554233}
\lim_{\delta\rightarrow 0}g_{\delta}(x,t)=f'(x,t).
\end{equation}

We shall prove that for every $\epsilon\in (0,1)$, the following
limit holds uniformly for $k\in\mathbb{N},$
\begin{equation}\label{gdfdr442323rr}
\lim_{\delta\rightarrow 0}\int_{\Omega^{0}_{\epsilon,k}\cap
U^{0}_{\epsilon,k}}|f'(x,u_{k}+v_{k})-g_{\delta}(x,u_{k}+v_{k})|^{\frac{q}{q-2}}dx=0.
\end{equation}
If not, then there exist $\epsilon_{0}>0,\ \eta_{0}>0$ and
sequences $\{\delta_{m}\},\ \{k_{m}\}$ which satisfy that
$\delta_{m}\rightarrow 0$ and
\begin{equation}\label{ncbvrrrere552233}
\int_{\Omega^{0}_{\epsilon_{0},k_{m}}\cap
U^{0}_{\epsilon_{0},k_{m}}}|f'(x,u_{k_{m}}+v_{k_{m}})-g_{\delta_{m}}(x,u_{k_{m}}+v_{k_{m}})|^{\frac{q}{q-2}}dx>\eta_{0},\
m=1,2,\cdots.
\end{equation}
By the condition $\bf (f_{2})$ and the definition of
$g_{\delta_{m}}$, we know that there exist $C, C'>0$ which are
independent of $m$ such that
\begin{eqnarray}\label{vvvscxiii888}
&&|f'(x,u_{k_{m}}(x)+v_{k_{m}}(x))-g_{\delta_{m}}(x,u_{k_{m}}(x)+v_{k_{m}}(x))|^{\frac{q}{q-2}}\nonumber\\
&\leq& C'|u_{k_{m}}(x)+v_{k_{m}}(x)|^{q}\nonumber \\
&\leq& C(\frac{1}{\epsilon_{0}^{q}}+|v_{k}(x)|^{q}), \ \forall
x\in \Omega^{1}_{\epsilon_{0},k_{m}}\cap
U^{0}_{\epsilon_{0},k_{m}}.
\end{eqnarray}
By (\ref{jjhsf554233}), we get that
\begin{equation}\label{yyyyyetrrr43435ff}
\lim_{m\rightarrow\infty}
|f'(x,u_{k_{m}}(x)+v_{k_{m}}(x))-g_{\delta_{m}}(x,u_{k_{m}}(x)+v_{k_{m}}(x))|^{\frac{q}{q-2}}=0
\ a.e.
\end{equation}
Moreover, by the fact that there exists constant $C>0$ which is
independent of $m$ such that
\begin{equation}
\epsilon^{q}_{0}\mbox{mes}(\Omega^{0}_{\epsilon_{0},k_{m}}\cap
U^{0}_{\epsilon_{0},k_{m}})\leq\int_{\Omega^{0}_{\epsilon_{0},k_{m}}\cap
U^{0}_{\epsilon_{0},k_{m}}}|u_{k_{m}}|^{q}\leq\int_{Q_{k}}|u_{k_{m}}|^{q}\leq
C
\end{equation}
we get that there exists constant $C>0$ which is independent of
$m$ such that
\begin{equation}\label{kkkkkjjj9998i8i8u}
\mbox{mes}(\Omega^{0}_{\epsilon_{0},k_{m}}\cap
U^{0}_{\epsilon_{0},k_{m}})\leq C.
\end{equation}
By (\ref{kkkkkjjj9998i8i8u}), we may assume that the limit
$\lim_{m\rightarrow\infty}\mbox{mes}(\Omega^{0}_{\epsilon_{0},k_{m}}\cap
U^{0}_{\epsilon_{0},k_{m}})$ exists.
 By (\ref{vvvscxiii888}) and
Fatou theorem, we get that
\begin{eqnarray}\label{bcvuuuuji898uu}
&&\int_{\Omega^{0}_{\epsilon_{0},k_{m}}\cap
U^{0}_{\epsilon_{0},k_{m}}}\liminf_{m\rightarrow\infty}(C(\frac{1}{\epsilon_{0}^{q}}+|v_{k}|^{q})-
|f'(x,u_{k_{m}}(x)+v_{k_{m}}(x))-g_{\delta_{m}}(x,u_{k_{m}}(x)+v_{k_{m}}(x))|^{\frac{q}{q-2}})dx\nonumber\\
&\leq&\liminf_{m\rightarrow\infty}\int_{\Omega^{0}_{\epsilon_{0},k_{m}}\cap
U^{0}_{\epsilon_{0},k_{m}}}(C(\frac{1}{\epsilon_{0}^{q}}+|v_{k}|^{q})-
|f'(x,u_{k_{m}}(x)+v_{k_{m}}(x))-g_{\delta_{m}}(x,u_{k_{m}}(x)+v_{k_{m}}(x))|^{\frac{q}{q-2}})dx\nonumber\\
\end{eqnarray}
By the fact that
$\lim_{k\rightarrow\infty}\int_{Q_{k}}|v_{k}|^{q}=0$ and
(\ref{bcvuuuuji898uu}), we get that
\begin{eqnarray}\label{klmuiuyy777766}
&&C\frac{1}{\epsilon_{0}^{q}}\lim_{m\rightarrow\infty}\mbox{mes}(\Omega^{0}_{\epsilon_{0},k_{m}}\cap
U^{0}_{\epsilon_{0},k_{m}})-\int_{\Omega^{0}_{\epsilon_{0},k_{m}}\cap
U^{0}_{\epsilon_{0},k_{m}}}\lim_{m\rightarrow\infty}
|f'(x,u_{k_{m}}+v_{k_{m}})-g_{\delta_{m}}(x,u_{k_{m}}+v_{k_{m}})|^{\frac{q}{q-2}}dx\nonumber\\
&\leq&
C\frac{1}{\epsilon_{0}^{q}}\lim_{m\rightarrow\infty}\mbox{mes}(\Omega^{0}_{\epsilon_{0},k_{m}}\cap
U^{0}_{\epsilon_{0},k_{m}})-\limsup_{m\rightarrow\infty}\int_{\Omega^{0}_{\epsilon_{0},k_{m}}\cap
U^{0}_{\epsilon_{0},k_{m}}}
|f'(x,u_{k_{m}}+v_{k_{m}})-g_{\delta_{m}}(x,u_{k_{m}}+v_{k_{m}})|^{\frac{q}{q-2}}dx\nonumber\\
\end{eqnarray}
By (\ref{klmuiuyy777766}) and (\ref{yyyyyetrrr43435ff}), we get
that
\begin{equation}\label{iopklmi88766}
\limsup_{m\rightarrow\infty}\int_{\Omega^{0}_{\epsilon_{0},k_{m}}\cap
U^{0}_{\epsilon_{0},k_{m}}}
|f'(x,u_{k_{m}}+v_{k_{m}})-g_{\delta_{m}}(x,u_{k_{m}}+v_{k_{m}})|^{\frac{q}{q-2}}dx=0.
\end{equation}
It contradicts to (\ref{ncbvrrrere552233}). Thus
(\ref{gdfdr442323rr}) holds uniformly for $k\in\mathbb{N}$. Then
by (\ref{gdfdr442323rr}) and
\begin{eqnarray}
&&\int_{\Omega^{0}_{\epsilon,k}\cap
U^{0}_{\epsilon,k}}|f'(x,u_{k}+v_{k})-g_{\delta}(x,u_{k}+v_{k})|\cdot
|\psi_{k}\cdot\varphi_{k}|dx\nonumber\\
&\leq&( \int_{\Omega^{0}_{\epsilon,k}\cap
U^{0}_{\epsilon,k}}|f'(x,u_{k}+v_{k})-g_{\delta}(x,u_{k}+v_{k})|^{\frac{q}{q-2}}dx)^{\frac{q-2}{q}}
(\int_{\mathbb{R}^{N}}|\psi_{k}|^{q})^{\frac{1}{q}}(\int_{\mathbb{R}^{N}}|\varphi_{k}|^{q})^{\frac{1}{q}}\nonumber\\
&\leq& C( \int_{\Omega^{0}_{\epsilon,k}\cap
U^{0}_{\epsilon,k}}|f'(x,u_{k}+v_{k})-g_{\delta}(x,u_{k}+v_{k})|^{\frac{q}{q-2}}dx)^{\frac{q-2}{q}}
\end{eqnarray}
we get that the following limit holds uniformly for
$k\in\mathbb{N}$, $||\varphi_{k}||_{k}\leq 1$ and
$||\psi_{k}||_{k}\leq 1$,
\begin{equation}\label{kkkklio9999i8}
\lim_{\delta\rightarrow 0}\int_{\Omega^{0}_{\epsilon,k}\cap
U^{0}_{\epsilon,k}}|f'(x,u_{k}+v_{k})-g_{\delta}(x,u_{k}+v_{k})|\cdot
|\psi_{k}\cdot\varphi_{k}|dx=0.
\end{equation}
As the same argument as (\ref{kkkklio9999i8}), we get that the
following limit holds uniformly for $k\in\mathbb{N}$,
$||\varphi_{k}||_{k}\leq 1$ and $||\psi_{k}||_{k}\leq 1$,
\begin{equation}\label{kkkg999i8}
\lim_{\delta\rightarrow 0}\int_{\Omega^{0}_{\epsilon,k}\cap
U^{0}_{\epsilon,k}}|f'(x,u_{k})-g_{\delta}(x,u_{k})|\cdot
|\psi_{k}\cdot\varphi_{k}|dx=0.
\end{equation}
By the definition of $g_{\delta}(x,t)$ and the condition $\bf
(f_{2})$, we get that for all $t\in [-2/\epsilon, 2/\epsilon]$,
there exists constant $C_{\epsilon}>0$ such that
\begin{equation}\label{yyyyeterrre665rr}
|\frac{\partial g_{\delta}}{\partial
t}(x,t)|=|\int_{\mathbb{R}}f'(x,s)\frac{\partial\eta_{\delta}}{\partial
s}(t-s)ds|\leq C_{\epsilon}/\delta.
\end{equation}

By (\ref{yyyyeterrre665rr}), we get that for all  $x\in
\Omega^{0}_{\epsilon,k}\cap U^{0}_{\epsilon,k},$
\begin{equation}\label{hgggggg555540008}
|g_{\delta}(x,u_{k}(x)+v_{k}(x))-g_{\delta}(x,u_{k}(x))|\leq
\frac{C_{\epsilon}}{\delta}|v_{k}(x)|.\end{equation} If $q-2\leq
1,$ then there exists constant $C(\epsilon)>0$ such that
$$|g_{\delta}(x,u_{k}(x)+v_{k}(x))-g_{\delta}(x,u_{k}(x))|\leq \frac{C_{\epsilon}}{\delta}|v_{k}(x)|\leq
\frac{C(\epsilon)}{\delta}|v_{k}(x)|^{q-2},\ \forall x\in
\Omega^{0}_{\epsilon,k}\cap U^{0}_{\epsilon,k}.$$ Then
\begin{eqnarray}\label{hhhhgafse}&&\int_{\Omega^{0}_{\epsilon,k}\cap U^{0}_{\epsilon,k}}
|g_{\delta}(x,u_{k}+v_{k})-g_{\delta}(x,u_{k})|\cdot
|\psi_{k}\cdot\varphi_{k}|dx\nonumber\\
&\leq& \frac{C(\epsilon)}{\delta}\int_{\Omega^{0}_{\epsilon,k}\cap
U^{0}_{\epsilon,k}}|v_{k}|^{q-2}|\varphi_{k}\cdot\psi_{k}|dx\nonumber\\
&\leq&
\frac{C(\epsilon)}{\delta}(\int_{\Omega^{0}_{\epsilon,k}\cap
U^{0}_{\epsilon,k}}|v_{k}|^{q})^{\frac{q-2}{q}}(\int_{\Omega^{0}_{\epsilon,k}\cap
U^{0}_{\epsilon,k}}|\varphi_{k}|^{q})^{\frac{1}{q}}
(\int_{\Omega^{0}_{\epsilon,k}\cap
U^{0}_{\epsilon,k}}|\psi_{k}|^{q})^{\frac{1}{q}}\nonumber\\
&\leq&
\frac{C_{2}(\epsilon)}{\delta}||v_{k}||^{q-2}_{H^{1}(Q_{k})},\end{eqnarray}
where $C_{2}(\epsilon)>0$ is a constant which is independent of
$\psi_{k},\ \varphi_{k}$ and $k$.

If $q-2>1$, then
\begin{eqnarray}\label{aaabscdlllllju}
&&\int_{\Omega^{0}_{\epsilon,k}\cap U^{0}_{\epsilon,k}}
|g_{\delta}(x,u_{k}+v_{k})-g_{\delta}(x,u_{k})|\cdot |\psi_{k}\cdot\varphi_{k}|dx\nonumber\\
&\leq&
\frac{C_{\epsilon}}{\delta}\int_{\Omega^{0}_{\epsilon,k}\cap
U^{0}_{\epsilon,k}}|v_{k}|\cdot |\varphi_{k}\cdot\psi_{k}|dx\nonumber\\
&\leq&\frac{C_{\epsilon}}{\delta}(\int_{\Omega^{0}_{\epsilon,k}\cap
U^{0}_{\epsilon,k}}|v_{k}|^{3})^{\frac{1}{3}}(\int_{\Omega^{0}_{\epsilon,k}\cap
U^{0}_{\epsilon,k}}|\psi_{k}|^{3})^{\frac{1}{3}}(\int_{\Omega^{0}_{\epsilon,k}\cap
U^{0}_{\epsilon,k}}|\varphi_{k}|^{3})^{\frac{1}{3}}\nonumber\\
&\leq&\frac{C'_{2}(\epsilon)}{\delta}(\int_{\Omega^{0}_{\epsilon,k}\cap
U^{0}_{\epsilon,k}}|v_{k}|^{3})^{\frac{1}{3}}\nonumber\\
&\leq& \frac{C''_{2}(\epsilon)}{\delta}||v_{k}||_{H^{1}(Q_{k})},
\end{eqnarray}
where $C''_{2}(\epsilon)>0$ is a constant which is independent of
$\psi_{k},\ \varphi_{k}$ and $k$.

By (\ref{hhhhgafse}) and (\ref{aaabscdlllllju}), we get that the
following limit holds uniformly for $ ||\varphi_{k}||_{k}\leq 1$
and $||\psi_{k}||_{k}\leq 1$,
\begin{equation}\label{uuyetrffhhhvsdxc00}
\lim_{k\rightarrow\infty}\int_{\Omega^{0}_{\epsilon,k}\cap
U^{0}_{\epsilon,k}}
|g_{\delta}(x,u_{k}+v_{k})-g_{\delta}(x,u_{k})|\cdot
|\psi_{k}\cdot\varphi_{k}|dx=0.
\end{equation}
By (\ref{kkkklio9999i8}), (\ref{kkkg999i8}) and
(\ref{uuyetrffhhhvsdxc00}), we get that the following limit holds
uniformly for $ ||\varphi_{k}||_{k}\leq 1$ and
$||\psi_{k}||_{k}\leq 1$,
\begin{equation}\label{uuy665c00}
\lim_{k\rightarrow\infty}\int_{\Omega^{0}_{\epsilon,k}\cap
U^{0}_{\epsilon,k}} |f'(x,u_{k}+v_{k})-f'(x,u_{k})|\cdot
|\psi_{k}\cdot\varphi_{k}|dx=0.
\end{equation}
Finally, by (\ref{nnna7yhhab444}), (\ref{hhhhh777y76}),
(\ref{huanabv5frthhhh}), (\ref{shenminchenke6gttgftfrr11111}) and
(\ref{uuy665c00}), we get that the limit
$$\lim_{k\rightarrow\infty}\int_{Q_{k}} |f'(x,u_{k}+v_{k})-f'(x,u_{k})|\cdot
|\psi_{k}\cdot\varphi_{k}|dx=0$$ holds uniformly for $
||\varphi_{k}||_{k}\leq 1$ and $||\psi_{k}||_{k}\leq 1$.
\hfill$\Box$

\bigskip

\noindent{\bf \large Proof of Lemma \ref{lantaojin}:}

\medskip

Firstly, we shall prove that as $k\rightarrow\infty,$ the limit
\begin{equation}\label{bbbbbbbb7nh990}
\int_{Q_{k}}|f'(x,u^{1}-\sum^{n}_{i=2}u^{i}(\cdot+b^{i}_{k}-b^{1}_{k}))
-f'(x,\sum^{n}_{i=2}u^{i}(\cdot+b^{i}_{k}-b^{1}_{k}))-f'(x,u^{1})|
\cdot|\varphi_{k}\cdot\psi_{k}|\rightarrow 0\end{equation} holds
uniformly for $\varphi_{k}, \varphi_{k}\in E_{k}$ which satisfy
$||\varphi_{k}||_{k}\leq 1,\ ||\psi_{k}||_{k}\leq 1.$

For convenience, we set
$v_{k}=\displaystyle\sum^{n}_{i=2}u^{i}(\cdot+b^{i}_{k}-b^{1}_{k}).$

Let $\Omega^{1}_{\epsilon,k}=\{x\in Q_{k}\ :\ |v_{k}(x)|\leq
1/\epsilon\},$ $\Omega^{2}_{\epsilon,k}=\{x\in Q_{k}\ :\
|v_{k}(x)|> 1/\epsilon\},$
 $$U^{1}_{\epsilon}=\{x\in Q_{k}\
:\ |u^{1}(x)|\leq 1/\epsilon\},\ U^{2}_{\epsilon}=\{x\in Q_{k}\ :\
|u^{1}(x)|> 1/\epsilon\}.$$ Then the limit
\begin{equation}\label{chaiyuzhinan7h}
\lim_{\epsilon\rightarrow 0}\mbox{mes}(\Omega^{2}_{\epsilon,k})=
0,
\end{equation}
holds uniformly for $k\in\mathbb{N}$ and
\begin{equation}\label{hdgfere666633448}
 \mbox{mes}(U^{2}_{\epsilon})\rightarrow 0,\ \mbox{as}\
\epsilon\rightarrow 0.
\end{equation}
Note that
\begin{equation}\label{ujna666666ufae}
\int_{Q_{k}}|f'(x,u^{1}+v_{k})-f'(x,v_{k})-f'(x,u^{1})|\cdot
|\varphi_{k}\cdot\psi_{k}|dx \leq\int_{\Omega^{1}_{\epsilon}\cap
U^{1}_{\epsilon}}+\int_{\Omega^{2}_{\epsilon}}+\int_{U^{2}_{\epsilon}}.\end{equation}
As the same argument as (\ref{huanabv5frthhhh}) and
(\ref{shenminchenke6gttgftfrr11111}), we know that  there exist
constants $r>0,\ C_{1}>0$ and $C_{2}>0$ which are independent of
$k, \varphi_{k}$ and $\psi_{k}$ such that
\begin{equation}\label{njjmloo0980000}
\int_{\Omega^{2}_{\epsilon,k}}|f'(x,u^{1}+v_{k})-f'(x,v_{k})-f'(x,u^{1})|\cdot
|\varphi_{k}\cdot\psi_{k}|dx\leq
C_{1}(\mbox{mes}(\Omega^{2}_{\epsilon,k}))^{\frac{1}{r}},\end{equation}
\begin{equation}\label{njiokmlz109}
\int_{U^{2}_{\epsilon}}|f'(x,u^{1}+v_{k})-f'(x,v_{k})-f'(x,u^{1})|\cdot
|\varphi_{k}\cdot\psi_{k}|dx\leq
C_{2}(\mbox{mes}(U^{2}_{\epsilon}))^{\frac{1}{r}}.\end{equation}
As the same argument as (\ref{kkkklio9999i8}) and
(\ref{kkkg999i8}), we know that for every $\epsilon\in (0,1)$, the
following limits holds uniformly for $k\in\mathbb{N},$
$||\varphi_{k}||_{k}\leq 1$ and $||\psi_{k}||_{k}\leq 1$,
\begin{equation}\label{gdfdr442323rr12}
\lim_{\delta\rightarrow 0}\int_{\Omega^{1}_{\epsilon,k}\cap
U^{1}_{\epsilon}}|f'(x,u^{1}+v_{k})-g_{\delta}(x,u^{1}+v_{k})|\cdot
|\varphi_{k}\cdot\psi_{k}|dx=0,
\end{equation}
\begin{equation}\label{gdfdr442323rr123}
\lim_{\delta\rightarrow 0}\int_{\Omega^{1}_{\epsilon,k}\cap
U^{1}_{\epsilon}}|f'(x,v_{k})-g_{\delta}(x,v_{k})|\cdot
|\varphi_{k}\cdot\psi_{k}|dx=0,
\end{equation}
\begin{equation}\label{gdfdr442323rr124}
\lim_{\delta\rightarrow 0}\int_{\Omega^{1}_{\epsilon,k}\cap
U^{1}_{\epsilon}}|f'(x,u^{1})-g_{\delta}(x,u^{1})|\cdot
|\varphi_{k}\cdot\psi_{k}|dx=0.
\end{equation}

 For $x\in
\Omega^{1}_{\epsilon,k}\cap U^{1}_{\epsilon}$, as the same
argument as (\ref{hgggggg555540008}), we deduce that there exists
constant  $M_{\epsilon,\delta}>0$ which  depends only on
$\epsilon$ and $\delta$ such that
\begin{equation}\label{guizhouqian7h11d}
|g_{\delta}(x,v_{k}(x)+u^{1}(x))-g_{\delta}(x,v_{k}(x))-g_{\delta}(x,u^{1}(x))|
\leq
M_{\epsilon,\delta}|u^{1}(x)|+C(|u^{1}(x)|^{q-2}+|u^{1}(x)|^{p-2}).
\end{equation}

We shall prove that for any  $\epsilon>0$, as $k\rightarrow\infty,$
the limit
\begin{equation}\label{acnnnuyte5}
\int_{\Omega^{1}_{\epsilon,k}\cap U^{1}_{\epsilon}}
|g_{\delta}(x,v_{k}+u^{1})-g_{\delta}(x,v_{k})-g_{\delta}(x,u^{1})|\cdot
|\varphi_{k}|\cdot|\psi_{k}|\rightarrow 0
\end{equation}
holds uniformly for $\varphi_{k},\psi_{k}\in E_{k}$ which satisfy
$||\varphi_{k}||_{k}\leq 1$ and $||\psi_{k}||_{k}\leq 1$.

We distinguish it by two cases:

\noindent{\bf Case 1.} $0<q-2\leq 1.$ In this case, by
(\ref{guizhouqian7h11d}), there exists
$\widetilde{M}_{\epsilon,\delta}>0$ such that
$$|g_{\delta}(x,v_{k}+u^{1})-g_{\delta}(x,v_{k})-g_{\delta}(x,u^{1})|\leq\widetilde{M}_{\epsilon,\delta}|u^{1}|^{q-2},
\ \forall x\in \Omega^{1}_{\epsilon,k}\cap U^{1}_{\epsilon}.$$ Then
$$|g_{\delta}(x,v_{k}+u^{1})-g_{\delta}(x,v_{k})-g_{\delta}(x,u^{1})|^{\frac{q}{q-2}}\leq
\widetilde{M}^{\frac{q}{q-2}}_{\epsilon,\delta}|u^{1}|^{q},\ \forall
x\in \Omega^{1}_{\epsilon,k}\cap U^{1}_{\epsilon}.$$ Since as
$k\rightarrow\infty,$
$$|g_{\delta}(x,v_{k}+u^{1})-g_{\delta}(x,v_{k})-g_{\delta}(x,u^{1})|^{\frac{q}{q-2}}\rightarrow 0\ a.e.,$$
by Lebesgue convergence theorem, we get that as
$k\rightarrow\infty,$
$$\int_{\Omega^{1}_{\epsilon,k}\cap U^{1}_{\epsilon}}
|g_{\delta}(x,v_{k}+u^{1})-g_{\delta}(x,v_{k})-g_{\delta}(x,u^{1})|^{\frac{q}{q-2}}\rightarrow
0.$$ It follows that as $k\rightarrow\infty,$ the limit
\begin{eqnarray}\label{bbbbbgads}
&&\int_{\Omega^{1}_{\epsilon,k}\cap U^{1}_{\epsilon}}
|g_{\delta}(x,v_{k}+u^{1})-g_{\delta}(x,v_{k})-g_{\delta}(x,u^{1})|\cdot
|\varphi_{k}\cdot\psi_{k}|dx\nonumber\\
&\leq& (\int_{\Omega^{1}_{\epsilon,k}\cap U^{1}_{\epsilon}}
|g_{\delta}(x,v_{k}+u^{1})-g_{\delta}(x,v_{k})-g_{\delta}(x,u^{1})|^{\frac{q}{q-2}})^{\frac{q-2}{q}}
(\int_{Q_{k}}|\varphi_{k}|^{q})^{\frac{1}{q}}(\int_{Q_{k}}|\psi_{k}|^{q})^{\frac{1}{q}}\rightarrow
0\nonumber
\end{eqnarray}
holds uniformly for $\varphi_{k}, \psi_{k}\in E_{k}$ which satisfy
$||\varphi_{k}||_{k}\leq 1,\ ||\psi_{k}||_{k}\leq 1$.

\noindent{\bf Case 2.} $q-2>1.$ In this case, by
(\ref{guizhouqian7h11d}), we get that there exists
$\widehat{M}_{\epsilon,\delta}>0$ such that
$$|g_{\delta}(x,v_{k}+u^{1})-g_{\delta}(x,v_{k})-g_{\delta}(x,u^{1})|\leq\widehat{M}_{\epsilon,\delta}|u^{1}|.$$
By Lebesgue convergence theorem, we deduce that, as
$k\rightarrow\infty,$ the limit
$$\int_{\Omega^{1}_{\epsilon,k}\cap U^{1}_{\epsilon}}
|g_{\delta}(x,v_{k}+u^{1})-g_{\delta}(x,v_{k})-g_{\delta}(x,u^{1})|^{3}\rightarrow
0$$ holds uniformly for $\varphi_{k}, \psi_{k}\in E_{k}$ which
satisfy $||\varphi_{k}||_{k}\leq 1,\ ||\psi_{k}||_{k}\leq 1$.
Therefore, as $k\rightarrow\infty,$ the limit
\begin{eqnarray}\label{bbbbbgads}
&&\int_{\Omega^{1}_{\epsilon,k}\cap U^{1}_{\epsilon}}
|g_{\delta}(x,v_{k}+u^{1})-g_{\delta}(x,v_{k})-g_{\delta}(x,u^{1})|\cdot
|\varphi_{k}\cdot\psi_{k}|dx\nonumber\\
&\leq& (\int_{\Omega^{1}_{\epsilon,k}\cap U^{1}_{\epsilon}}
|g_{\delta}(x,v_{k}+u^{1})-g_{\delta}(x,v_{k})-g_{\delta}(x,u^{1})|^{3})^{\frac{1}{3}}
(\int_{Q_{k}}|\varphi_{k}|^{3})^{\frac{1}{3}}(\int_{Q_{k}}|\psi_{k}|^{3})^{\frac{1}{3}}\rightarrow
0\nonumber
\end{eqnarray}
holds uniformly for $\varphi_{k}, \psi_{k}\in E_{k}$ which satisfy
$||\varphi_{k}||_{k}\leq 1,\ ||\psi_{k}||_{k}\leq 1$.

Thus (\ref{acnnnuyte5}) holds. By
$(\ref{chaiyuzhinan7h})-(\ref{gdfdr442323rr124})$ and
(\ref{acnnnuyte5}), we get (\ref{bbbbbbbb7nh990}). Finally, by
inductive argument, we can get the desired result of this Lemma.
\hfill$\Box$

\bigskip

\noindent{\bf\large Proof of Lemma \ref{nbbzvcccfrtyujhgppp}:}

\medskip

\noindent{\bf (1).} By Lemma \ref{taolichunfen}, we get that for
any $\varphi\in (N^{k}_{i})^{\bot}$,
\begin{eqnarray}\label{zhongbaxintaohuan}
0&=&\int_{Q_{k}}\nabla
(\widehat{u}^{i}_{k}+\widehat{w}^{i}_{k}(\sum^{l_i}_{j=1}x_{i,j}\chi_{k}e_{i,j})
+\sum^{l_i}_{j=1}x_{i,j}\chi_{k}e_{i,j})\cdot\nabla\varphi\nonumber\\
&&+\int_{Q_{k}}V(x)\cdot
(\widehat{u}^{i}_{k}+\widehat{w}^{i}_{k}(\sum^{l_i}_{j=1}x_{i,j}\chi_{k}e_{i,j})
+\sum^{l_i}_{j=1}x_{i,j}\chi_{k}e_{i,j})\cdot\varphi\nonumber\\
&&-\int_{Q_{k}}f(x,\widehat{u}^{i}_{k}+\widehat{w}^{i}_{k}(\sum^{l_i}_{j=1}x_{i,j}\chi_{k}e_{i,j})
+\sum^{l_i}_{j=1}x_{i,j}\chi_{k}e_{i,j})\varphi,
\end{eqnarray}
and for any $\psi\in N_{i}^{\bot},$
\begin{eqnarray}\label{zhongbaxintaohuankkn}
0&=&\int_{\mathbb{R}^{N}}\nabla
(u^{i}+\omega^{i}(\sum^{l_i}_{j=1}x_{i,j}e_{i,j})
+\sum^{l_i}_{j=1}x_{i,j}e_{i,j})\cdot\nabla\psi\nonumber\\
&&+\int_{\mathbb{R}^{N}}V(x)\cdot
(u^{i}+\omega^{i}(\sum^{l_i}_{j=1}x_{i,j}e_{i,j})
+\sum^{l_i}_{j=1}x_{i,j}e_{i,j})\cdot\psi\nonumber\\
&&-\int_{\mathbb{R}^{N}}f(x,u^{i}+\omega^{i}(\sum^{l_i}_{j=1}x_{i,j}e_{i,j})+\sum^{l_i}_{j=1}x_{i,j}e_{i,j})\psi.
\end{eqnarray}
Note that the limit
\begin{equation}\label{huannnnamkkjuippk}
||\omega^{i}(\sum^{l_i}_{j=1}x_{i,j}e_{i,j})-\chi_{k}
\omega^{i}(\sum^{l_i}_{j=1}x_{i,j}e_{i,j})||_{H^{1}(\mathbb{R}^{N})}\rightarrow
0,\ k\rightarrow\infty
\end{equation}
holds uniformly for
$\sum^{l_{i}}_{j=1}x^{2}_{i,j}\leq\delta^{2}_{0}$ and
\begin{equation}\label{shuinuanyaxianziyg5}
||e_{i,j}-\chi_{k}e_{i,j}||_{H^{1}(\mathbb{R}^{N})}\rightarrow 0,\
||u^{i}-\widehat{u}^{i}_{k}||_{k}\rightarrow 0,\ k\rightarrow\infty.
\end{equation}
By $(\ref{zhongbaxintaohuankkn})- (\ref{shuinuanyaxianziyg5})$ and
Lemma \ref{yadanhexiawa}, we get that as $k\rightarrow\infty,$
\begin{eqnarray}\label{wubianluomuxiaoxiaoxiaoolk}
o(1)&=&\int_{Q_{k}}\nabla
(\widehat{u}^{i}+\chi_{k}\omega^{i}(\sum^{l_i}_{j=1}x_{i,j}e_{i,j})
+\chi_{k}\sum^{l_i}_{j=1}x_{i,j}e_{i,j})\cdot\nabla\psi\nonumber\\
&&+\int_{Q_{k}}V(x)\cdot
(\widehat{u}^{i}+\chi_{k}\omega^{i}(\sum^{l_i}_{j=1}x_{i,j}e_{i,j})
+\chi_{k}\sum^{l_i}_{j=1}x_{i,j}e_{i,j})\cdot\psi\nonumber\\
&&-\int_{Q_{k}}f(x,\widehat{u}^{i}+\chi_{k}\omega^{i}(\sum^{l_i}_{j=1}x_{i,j}e_{i,j})
+\chi_{k}\sum^{l_i}_{j=1}x_{i,j}e_{i,j})\psi
\end{eqnarray}
holds uniformly for  $\psi\in N^{\bot}_{i}$ satisfying $||\psi||\leq
1$ and $\sum^{l_{i}}_{j=1}x^{2}_{i,j}\leq\delta^{2}_{0}.$ It is easy
to verify that the limit
\begin{equation}\label{hhhjuuuyioooooopk778}
(\chi_{k}\omega^{i}(\sum^{l_i}_{j=1}x_{i,j}e_{i,j}),\
\chi_{k}e_{i,j})_{k}\rightarrow 0,\ k\rightarrow\infty
\end{equation}
holds uniformly for
$\sum^{l_{i}}_{j=1}x^{2}_{i,j}\leq\delta^{2}_{0},$ since
$(\omega^{i}(\sum^{l_i}_{j=1}x_{i,j}e_{i,j}),\ e_{i,j})=0.$ Recall
that $P_{(\widehat{N}^{k}_{i})^{\bot}}$ is the orthogonal projection
 from $E_{k}$ into
$(\widehat{N}^{k}_{i})^{\bot}$, by (\ref{hhhjuuuyioooooopk778}) and
(\ref{huannnnamkkjuippk}), we get that the limit
\begin{equation}\label{tonzhishannianmmmkiu}
||\omega^{i}(\sum^{l_i}_{j=1}x_{i,j}e_{i,j})-P_{(\widehat{N}^{k}_{i})^{\bot}}
(\chi_{k}\omega^{i}(\sum^{l_i}_{j=1}x_{i,j}e_{i,j}))||_{H^{1}(\mathbb{R}^{N})}\rightarrow
0,\ k\rightarrow\infty
\end{equation}
holds uniformly for
$\sum^{l_{i}}_{j=1}x^{2}_{i,j}\leq\delta^{2}_{0}.$ By
(\ref{tonzhishannianmmmkiu}) and (\ref{wubianluomuxiaoxiaoxiaoolk}),
we get that as $k\rightarrow\infty,$
\begin{eqnarray}\label{wubianluomuxiaoxiaoxiaoolkggg6781}
o(1)&=&\int_{Q_{k}}\nabla
(\widehat{u}^{i}+P_{(\widehat{N}^{k}_{i})^{\bot}}
(\chi_{k}\omega^{i}(\sum^{l_i}_{j=1}x_{i,j}e_{i,j}))
+\chi_{k}\sum^{l_i}_{j=1}x_{i,j}e_{i,j})\cdot\nabla\psi\nonumber\\
&&+\int_{Q_{k}}V(x)\cdot
(\widehat{u}^{i}+P_{(\widehat{N}^{k}_{i})^{\bot}}
(\chi_{k}\omega^{i}(\sum^{l_i}_{j=1}x_{i,j}e_{i,j}))
+\chi_{k}\sum^{l_i}_{j=1}x_{i,j}e_{i,j})\cdot\psi\nonumber\\
&&-\int_{Q_{k}}f(x,\widehat{u}^{i}+P_{(\widehat{N}^{k}_{i})^{\bot}}
(\chi_{k}\omega^{i}(\sum^{l_i}_{j=1}x_{i,j}e_{i,j}))
+\chi_{k}\sum^{l_i}_{j=1}x_{i,j}e_{i,j})\psi.
\end{eqnarray}

Choose
$\varphi=\widehat{w}^{i}_{k}(\sum^{l_i}_{j=1}x_{i,j}\chi_{k}e_{i,j})
-P_{(\widehat{N}^{k}_{i})^{\bot}}
(\chi_{k}\omega^{i}(\sum^{l_i}_{j=1}x_{i,j}e_{i,j}))$ and let
$\psi\in H^{1}(\mathbb{R}^{N})$ be an extension of $\varphi.$ We
have $\varphi\in(\widehat{N}^{k}_{i})^{\bot}$.  Then minus
(\ref{zhongbaxintaohuan}) by
(\ref{wubianluomuxiaoxiaoxiaoolkggg6781}), by mean value theorem, we
get that
\begin{eqnarray}\label{zhengshihetunyushangshi}
\left(P_{(\widehat{N}^{k}_{i})^{\bot}}(\nabla^{2}J_{k}
(u_{i,k,t}))\varphi,\varphi\right)_{k}\rightarrow 0,\
k\rightarrow\infty
\end{eqnarray}
holds uniformly for
$\sum^{l_{i}}_{j=1}x^{2}_{i,j}\leq\delta^{2}_{0},$ where
$$u_{i,k,t}=
\widehat{u}^{i}_{k}+(1-t)\widehat{w}^{i}_{k}(\sum^{l_i}_{j=1}x_{i,j}\chi_{k}e_{i,j})
+tP_{(\widehat{N}^{k}_{i})^{\bot}}
(\chi_{k}\omega^{i}(\sum^{l_i}_{j=1}x_{i,j}e_{i,j}))+\sum^{l_i}_{j=1}x_{i,j}\chi_{k}e_{i,j}$$
and $t\in [0,1].$ By Lemma \ref{xiaojiazhi}, we know that the
operator $P_{(\widehat{N}^{k}_{i})^{\bot}}(\nabla^{2}J_{k}
(u_{i,k,t})|_{(\widehat{N}^{k}_{i})^{\bot}})$ is invertible and
there exists constant $\eta>0$ such that
\begin{equation}\label{weiyuyansuanfei}
||(P_{(\widehat{N}^{k}_{i})^{\bot}}(\nabla^{2}J_{k}
(u_{i,k,t})|_{(\widehat{N}^{k}_{i})^{\bot}}))^{-1}||\leq \eta.
\end{equation}
By (\ref{zhengshihetunyushangshi}) and (\ref{weiyuyansuanfei}), we
get that \begin{equation}\label{honminlilllqushi}
||\varphi||_{k}\rightarrow 0,\ k\rightarrow\infty
\end{equation}
holds uniformly for
$\sum^{l_{i}}_{j=1}x^{2}_{i,j}\leq\delta^{2}_{0}.$ Thus the result
(1) of this Lemma follows from (\ref{honminlilllqushi}) and
(\ref{tonzhishannianmmmkiu}) directly.

\bigskip

\noindent{\bf (2).} Differentiating the equalities
(\ref{zhongbaxintaohuan}) and (\ref{zhongbaxintaohuankkn}) for the
variable $x_{i,s}$, we get that for any $\varphi\in
(N^{k}_{i})^{\bot}$,
\begin{eqnarray}\label{zhongbaxintaohuanggg6hh}
0&=&\int_{Q_{k}}\nabla
((\widehat{w}^{i}_{k})'(\sum^{l_i}_{j=1}x_{i,j}\chi_{k}e_{i,j})\chi_{k}e_{i,s}
+\chi_{k}e_{i,s})\cdot\nabla\varphi\nonumber\\
&&+\int_{Q_{k}}V(x)\cdot
((\widehat{w}^{i}_{k})'(\sum^{l_i}_{j=1}x_{i,j}\chi_{k}e_{i,j})\chi_{k}e_{i,s}
+\chi_{k}e_{i,s})\cdot\varphi\nonumber\\
&&-\int_{Q_{k}}f'(x,\widehat{u}^{i}_{k}+\widehat{w}^{i}_{k}(\sum^{l_i}_{j=1}x_{i,j}\chi_{k}e_{i,j})
+\sum^{l_i}_{j=1}x_{i,j}\chi_{k}e_{i,j})\cdot((\widehat{w}^{i}_{k})'(\sum^{l_i}_{j=1}x_{i,j}\chi_{k}e_{i,j})
\chi_{k}e_{i,s})
\cdot\varphi\nonumber\\
&&-\int_{Q_{k}}f'(x,\widehat{u}^{i}_{k}+\widehat{w}^{i}_{k}(\sum^{l_i}_{j=1}x_{i,j}\chi_{k}e_{i,j})
+\sum^{l_i}_{j=1}x_{i,j}\chi_{k}e_{i,j})\cdot(\chi_{k}e_{i,s})
\cdot\varphi
\end{eqnarray}
and for any $\psi\in N_{i}^{\bot},$
\begin{eqnarray}\label{zhongbaxintaohuankknuujhnb0000}
0&=&\int_{\mathbb{R}^{N}}\nabla
((\omega^{i})'(\sum^{l_i}_{j=1}x_{i,j}e_{i,j})e_{i,s}
+e_{i,s})\cdot\nabla\psi\nonumber\\
&&+\int_{\mathbb{R}^{N}}V(x)\cdot
((\omega^{i})'(\sum^{l_i}_{j=1}x_{i,j}e_{i,j})e_{i,s}
+e_{i,s})\cdot\psi\nonumber\\
&&-\int_{\mathbb{R}^{N}}f'(x,u^{i}+\omega^{i}(\sum^{l_i}_{j=1}x_{i,j}e_{i,j})+\sum^{l_i}_{j=1}x_{i,j}e_{i,j})
\cdot((\omega^{i})'(\sum^{l_i}_{j=1}x_{i,j}e_{i,j})e_{i,s})\cdot\psi\nonumber\\
&&-\int_{\mathbb{R}^{N}}f'(x,u^{i}+\omega^{i}(\sum^{l_i}_{j=1}x_{i,j}e_{i,j})+\sum^{l_i}_{j=1}x_{i,j}e_{i,j})
\cdot e_{i,s}\cdot\psi.
\end{eqnarray}
Note that the limit
\begin{equation}\label{huannnnamkkjuippknnbv}
||(\omega^{i})'(\sum^{l_i}_{j=1}x_{i,j}e_{i,j})e_{i,s}-\chi_{k}
(\omega^{i})'(\sum^{l_i}_{j=1}x_{i,j}e_{i,j})e_{i,s}||_{H^{1}(\mathbb{R}^{N})}\rightarrow
0,\ k\rightarrow\infty
\end{equation}
holds uniformly for
$\sum^{l_{i}}_{j=1}x^{2}_{i,j}\leq\delta^{2}_{0}.$ By
(\ref{huannnnamkkjuippk}), (\ref{shuinuanyaxianziyg5}),
(\ref{huannnnamkkjuippknnbv}), the result (1) of Lemma
\ref{nbbzvcccfrtyujhgppp} and Lemma \ref{yadanhexiawa}, we get that
as $k\rightarrow\infty,$
\begin{eqnarray}\label{tapainqinashanreyyyqhfsd}
&&\int_{\mathbb{R}^{N}}f'(x,u^{i}+\omega^{i}(\sum^{l_i}_{j=1}x_{i,j}e_{i,j})+\sum^{l_i}_{j=1}x_{i,j}e_{i,j})
\cdot((\omega^{i})'(\sum^{l_i}_{j=1}x_{i,j}e_{i,j})e_{i,s})\cdot\psi\nonumber\\
&=&\int_{Q_{k}}f'(x,\widehat{u}^{i}_{k}+\widehat{w}^{i}_{k}(\sum^{l_i}_{j=1}x_{i,j}\chi_{k}e_{i,j})
+\sum^{l_i}_{j=1}x_{i,j}\chi_{k}e_{i,j})
\cdot(\chi_{k}(\omega^{i})'(\sum^{l_i}_{j=1}x_{i,j}e_{i,j})e_{i,s})\cdot\psi\nonumber\\
&&+o(1)
\end{eqnarray}
and \begin{eqnarray}\label{bushilianzhinanwab776g}
&&\int_{\mathbb{R}^{N}}f'(x,u^{i}+\omega^{i}(\sum^{l_i}_{j=1}x_{i,j}e_{i,j})+\sum^{l_i}_{j=1}x_{i,j}e_{i,j})
\cdot e_{i,s}\cdot\psi\nonumber\\
&=&\int_{Q_{k}}f'(x,\widehat{u}^{i}_{k}+\widehat{w}^{i}_{k}(\sum^{l_i}_{j=1}x_{i,j}\chi_{k}e_{i,j})
+\sum^{l_i}_{j=1}x_{i,j}\chi_{k}e_{i,j})\cdot(\chi_{k}e_{i,s})
\cdot\psi+o(1)\nonumber\\
\end{eqnarray}
By (\ref{huannnnamkkjuippk}), (\ref{shuinuanyaxianziyg5}) and
$(\ref{zhongbaxintaohuanggg6hh})-(\ref{bushilianzhinanwab776g})$, we
can get that
\begin{eqnarray}\label{zhongbaxintaohuanggg6hhuuuh78j}
&&\int_{Q_{k}}\nabla
((\widehat{w}^{i}_{k})'(\sum^{l_i}_{j=1}x_{i,j}\chi_{k}e_{i,j})\chi_{k}e_{i,s}
+\chi_{k}e_{i,s})\cdot\nabla\varphi\nonumber\\
&&+\int_{Q_{k}}V(x)\cdot
((\widehat{w}^{i}_{k})'(\sum^{l_i}_{j=1}x_{i,j}\chi_{k}e_{i,j})\chi_{k}e_{i,s}
+\chi_{k}e_{i,s})\cdot\varphi\nonumber\\
&&-\int_{Q_{k}}f'(x,\widehat{u}^{i}_{k}+\widehat{w}^{i}_{k}(\sum^{l_i}_{j=1}x_{i,j}\chi_{k}e_{i,j})
+\sum^{l_i}_{j=1}x_{i,j}\chi_{k}e_{i,j})
\cdot((\widehat{w}^{i}_{k})'(\sum^{l_i}_{j=1}x_{i,j}\chi_{k}e_{i,j})\chi_{k}e_{i,s})
\cdot\varphi\nonumber\\
&=&\int_{Q_{k}}f'(x,\widehat{u}^{i}_{k}+\widehat{w}^{i}_{k}(\sum^{l_i}_{j=1}x_{i,j}\chi_{k}e_{i,j})
+\sum^{l_i}_{j=1}x_{i,j}\chi_{k}e_{i,j})\cdot(\chi_{k}e_{i,s})
\cdot\varphi
\end{eqnarray}
and
\begin{eqnarray}\label{zhongbaxintaohuankknuujhnb77y0000}
&&\int_{Q_{k}}\nabla
(\chi_{k}(\omega^{i})'(\sum^{l_i}_{j=1}x_{i,j}e_{i,j})e_{i,s}
+\chi_{k}e_{i,s})\cdot\nabla\psi\nonumber\\
&&+\int_{Q_{k}}V(x)\cdot
(\chi_{k}(\omega^{i})'(\sum^{l_i}_{j=1}x_{i,j}e_{i,j})e_{i,s}
+\chi_{k}e_{i,s})\cdot\psi\nonumber\\
&&-\int_{Q_{k}}f'(x,\widehat{u}^{i}+\widehat{w}_{k}^{i}(\sum^{l_i}_{j=1}x_{i,j}\chi_{k}e_{i,j})
+\sum^{l_i}_{j=1}x_{i,j}\chi_{k}e_{i,j})
\cdot(\chi_{k}(\omega^{i})'(\sum^{l_i}_{j=1}x_{i,j}e_{i,j})e_{i,s})\cdot\psi\nonumber\\
&=&\int_{Q_{k}}f'(x,\widehat{u}^{i}+\widehat{w}_{k}^{i}(\sum^{l_i}_{j=1}x_{i,j}\chi_{k}e_{i,j})
+\chi_{k}\sum^{l_i}_{j=1}x_{i,j}e_{i,j}) \cdot
\chi_{k}e_{i,s}\cdot\psi+o(1).
\end{eqnarray}
Since $((\omega^{i})'(\sum^{l_i}_{j=1}x_{i,j}e_{i,j})e_{i,s},\
e_{i,t})=0,$ we get that
\begin{equation}\label{nnnhafff4f56h}
((\chi_{k}\omega^{i})'(\sum^{l_i}_{j=1}x_{i,j}e_{i,j})e_{i,s},\
\chi_{k}e_{i,t})=o(1),\ k\rightarrow\infty
\end{equation}
holds uniformly for
$\sum^{l_{i}}_{j=1}x^{2}_{i,j}\leq\delta^{2}_{0}$. By
(\ref{nnnhafff4f56h}) and (\ref{huannnnamkkjuippknnbv}), we get that
the limit \begin{equation}\label{tonzhishannianmmmkiuttt}
||(\omega^{i})'(\sum^{l_i}_{j=1}x_{i,j}e_{i,j})e_{i,s}-P_{(\widehat{N}^{k}_{i})^{\bot}}
(\chi_{k}(\omega^{i})'(\sum^{l_i}_{j=1}x_{i,j}e_{i,j})e_{i,s})||_{H^{1}(\mathbb{R}^{N})}\rightarrow
0,\ k\rightarrow\infty
\end{equation}
holds uniformly for
$\sum^{l_{i}}_{j=1}x^{2}_{i,j}\leq\delta^{2}_{0}.$ By
(\ref{zhongbaxintaohuankknuujhnb77y0000}) and
(\ref{tonzhishannianmmmkiuttt}), we get that
\begin{eqnarray}\label{zhongbaxintaohuankkwwwwnb77y0000}
&&\int_{Q_{k}}\nabla (P_{(\widehat{N}^{k}_{i})^{\bot}}
(\chi_{k}(\omega^{i})'(\sum^{l_i}_{j=1}x_{i,j}e_{i,j})e_{i,s})
+\chi_{k}e_{i,s})\cdot\nabla\psi\nonumber\\
&&+\int_{Q_{k}}V(x)\cdot (P_{(\widehat{N}^{k}_{i})^{\bot}}
(\chi_{k}(\omega^{i})'(\sum^{l_i}_{j=1}x_{i,j}e_{i,j})e_{i,s})
+\chi_{k}e_{i,s})\cdot\psi\nonumber\\
&&-\int_{Q_{k}}f'(x,\widehat{u}^{i}+\widehat{w}_{k}^{i}(\sum^{l_i}_{j=1}x_{i,j}\chi_{k}e_{i,j})
+\sum^{l_i}_{j=1}x_{i,j}\chi_{k}e_{i,j})
\cdot(P_{(\widehat{N}^{k}_{i})^{\bot}}
(\chi_{k}(\omega^{i})'(\sum^{l_i}_{j=1}x_{i,j}e_{i,j})e_{i,s}))\cdot\psi\nonumber\\
&=&\int_{Q_{k}}f'(x,\widehat{u}^{i}+\widehat{w}_{k}^{i}(\sum^{l_i}_{j=1}x_{i,j}\chi_{k}e_{i,j})
+\chi_{k}\sum^{l_i}_{j=1}x_{i,j}e_{i,j}) \cdot
\chi_{k}e_{i,s}\cdot\psi+o(1).
\end{eqnarray}
Choose
$\varphi=(\widehat{w}^{i}_{k})'(\sum^{l_i}_{j=1}x_{i,j}\chi_{k}e_{i,j})\chi_{k}e_{i,s}
-P_{(\widehat{N}^{k}_{i})^{\bot}}
(\chi_{k}(\omega^{i})'(\sum^{l_i}_{j=1}x_{i,j}e_{i,j})e_{i,s})$ and
let $\psi\in H^{1}(\mathbb{R}^{N})$ be an extension of $\varphi.$ We
have $\varphi\in(\widehat{N}^{k}_{i})^{\bot}$. Then minus
(\ref{zhongbaxintaohuanggg6hhuuuh78j}) by
(\ref{zhongbaxintaohuankkwwwwnb77y0000}), we get that
\begin{eqnarray}\label{zhengshihetunyushangshiuuuhg}
\left(P_{(\widehat{N}^{k}_{i})^{\bot}}(\nabla^{2}J_{k}
(u_{i,k}))\varphi,\ \varphi\right)_{k}\rightarrow 0,\
k\rightarrow\infty
\end{eqnarray}
holds uniformly for
$\sum^{l_{i}}_{j=1}x^{2}_{i,j}\leq\delta^{2}_{0},$ where
$u_{i,k}=\widehat{u}^{i}+\widehat{w}_{k}^{i}(\sum^{l_i}_{j=1}x_{i,j}\chi_{k}e_{i,j})
+\chi_{k}\sum^{l_i}_{j=1}x_{i,j}e_{i,j}.$ By Lemma \ref{xiaojiazhi},
we know that the operator
$P_{(\widehat{N}^{k}_{i})^{\bot}}(\nabla^{2}J_{k}
(u_{i,k})|_{(\widehat{N}^{k}_{i})^{\bot}})$ is invertible and there
exists constant $M>0$ such that
\begin{equation}\label{weiyuyansuanfeiyyyga4}
||(P_{(\widehat{N}^{k}_{i})^{\bot}}(\nabla^{2}J_{k}
(u_{i,k})|_{(\widehat{N}^{k}_{i})^{\bot}}))^{-1}||\leq M.
\end{equation}
By (\ref{zhengshihetunyushangshiuuuhg}) and
(\ref{weiyuyansuanfeiyyyga4}), we get that
\begin{equation}\label{honminlilllqushiff} ||\varphi||_{k}\rightarrow
0,\ k\rightarrow\infty
\end{equation}
holds uniformly for
$\sum^{l_{i}}_{j=1}x^{2}_{i,j}\leq\delta^{2}_{0}.$ Thus the result
(2) of this Lemma follows from (\ref{honminlilllqushiff}) and
(\ref{tonzhishannianmmmkiuttt}) directly.\hfill$\Box$

\bigskip

\noindent{\bf\large   Proof of Lemma \ref{tiandiyisa}:}

\bigskip

\noindent{\bf (1).} Let
$\theta^{i}_{k}(h_{i})=w^{i}_{k}(h_{i})-\widetilde{T}_{k}(w^{i}_{k}(h_{i})),\
i=1,2,\cdots,n.$ By (\ref{tapainqinshan}) and the fact that
$|b^{i}_{k}-b^{j}_{k}|\rightarrow\infty$ as $k\rightarrow\infty$
for $i\neq j,$ we get that
\begin{eqnarray}\label{shuimen7hytrqq93}
\left(w^{i}_{k}(\sum^{l_{i}}_{s=1}x_{i,s}e^{k}_{i,s}),\
e^{k}_{j,t}\right)_{k}
&=&\left((\widehat{w}^{i}_{k}(\sum^{l_{i}}_{s=1}x_{i,s}\chi_{k}e_{i,j}))
(\cdot+b^{i}_{k}),\ e^{k}_{j,t}\right)_{k}\nonumber\\
&=&\left((\widehat{w}^{i}_{k}(\sum^{l_{i}}_{s=1}x_{i,s}\chi_{k}e_{i,j}))
,\ (\chi_{k}e_{j,t})(\cdot+b^{j}_{k}-b^{i}_{k})\right)_{k}\nonumber\\
&=&\left\{
\begin{array}{l}
0,  \quad\quad i=j\\
o(1), \ \mbox{as}\ k\rightarrow\infty,\  i\neq j.\\
\end{array} \right.
\end{eqnarray}
By (\ref{shuimen7hytrqq93}), we get that the limit
\begin{equation}\label{haojie66666eqwc}
||\theta^{i}_{k}(h_{i})||_{k}\rightarrow 0\ \mbox{ as}\
k\rightarrow\infty\end{equation} holds uniformly for $h_{i}$
satisfying $|||h_{i}|||\leq\delta_{0}.$ Thus by the result (1) of
Lemma \ref{taolichunfen}, we know that as $k\rightarrow\infty,$ the
limit
$$\widetilde{T}_{k}(\nabla J_{k}(u^{i}_{k}+\widetilde{T}_{k}(w^{i}_{k}(h_{i}))+h_{i}))
\rightarrow 0$$ holds uniformly for $h_{i}$ satisfying
$|||h_{i}|||\leq\delta_{0}.$ By Lemma \ref{tongrenttthqqqw} and
the fact that
$\lim_{k\rightarrow\infty}||u_{k}-\sum^{n}_{i=1}u^{i}_{k}||_{H^{1}(Q_{k})}=0$,
we get that as $k\rightarrow\infty,$ the limit
$$\int_{Q_{k}}|f(x,u_{k}+\sum^{n}_{i=1}\widetilde{T}_{k}(w^{i}_{k}(h_{i}))+\sum^{n}_{i=1}
h_{i})-\sum^{n}_{i=1}f(x,u^{i}_{k}+\widetilde{T}_{k}(w^{i}_{k}(h_{i}))+h_{i})|\cdot
|\varphi_{k}|\rightarrow 0$$ holds uniformly for
$h=\sum^{n}_{i=1}h_{i}$  satisfying $|||h|||\leq\delta_{0}$ and
$\varphi_{k}\in E_{k}$ satisfying $||\varphi_{k}||_{k}\leq 1$.
Thus we deduce that as $k\rightarrow\infty,$ the limit
\begin{equation}\label{i99i88iiii76ii}
\widetilde{T}_{k}(\nabla
J_{k}(u_{k}+\sum^{n}_{i=1}\widetilde{T}_{k}(w^{i}_{k}(h_{i}))+\sum^{n}_{i=1}h_{i}))\rightarrow
0\end{equation} holds uniformly for $h=\sum^{n}_{i=1}h_{i}$
satisfying $|||h|||\leq\delta_{0}.$ Since $\widetilde{T}_{k}(\nabla
J_{k}(u_{k}+w_{k}(\sum^{n}_{i=1}h_{i})+\sum^{n}_{i=1}h_{i}))=0,$ by
(\ref{i99i88iiii76ii}), we get that as $k\rightarrow\infty,$
\begin{eqnarray}\label{changjianhoulan}
&&\widetilde{T}_{k}(\nabla
J_{k}(u_{k}+w_{k}(\sum^{n}_{i=1}h_{i})+\sum^{n}_{i=1}h_{i})) -
\widetilde{T}_{k}(\nabla
J_{k}(u_{k}+\sum^{n}_{i=1}\widetilde{T}_{k}(w^{i}_{k}(h_{i}))+\sum^{n}_{i=1}h_{i}))\nonumber\\
&=&\left\{\int^{1}_{0}\widetilde{T}_{k}(\nabla^{2}
J_{k}(u_{k}+(1-t)w_{k}(\sum^{n}_{i=1}h_{i})+t\sum^{n}_{i=1}\widetilde{T}_{k}(w^{i}_{k}(h_{i}))
+\sum^{n}_{i=1}h_{i}))dt\right\}\nonumber\\
&&\times(\sum^{n}_{i=1}\widetilde{T}_{k}(w^{i}_{k}(h_{i}))
-w_{k}(\sum^{n}_{i=1}h_{i}))\nonumber\\
&=&o(1)
\end{eqnarray}
holds uniformly for $h\in\overline{B_{\Lambda_{k}}(0,\delta_{0})}.$
By Lemma \ref{xiaojiazhi}, we know that when $k\geq k_{0},$
$$||(\int^{1}_{0}\widetilde{T}_{k}(\nabla
J_{k}(u_{k}+(1-t)w_{k}(\sum^{n}_{i=1}h_{i})+t\sum^{n}_{i=1}\widetilde{T}_{k}(\omega^{i}_{k}(h_{i}))
+\sum^{n}_{i=1}h_{i}))dt)^{-1}||\leq\eta$$ holds uniformly for
$h=\sum^{n}_{i=1}\overline{B_{\Lambda_{k}}(0,\delta_{0})}.$ Thus we
get that as $k\rightarrow\infty,$
\begin{equation}\label{uuuuu6666717777787}
\sup
\{||w_{k}(\sum^{n}_{i=1}h_{i})-\sum^{n}_{i=1}\widetilde{T}_{k}(w^{i}_{k}(h_{i}))||_{k}
\ : \ h=\sum^{n}_{i=1}h_{i}\in
\overline{B_{\Lambda_{k}}(0,\delta_{0})} \}\rightarrow
0.\end{equation} The result (1) of this lemma follows from
(\ref{uuuuu6666717777787}) and (\ref{haojie66666eqwc}) directly.

\bigskip

\noindent{\bf (2).} By Lemma \ref{taolichunfen}, we know that for
any $\varphi\in\Pi_{k},$
\begin{eqnarray}\label{fenqiyunyong}
0&=&(\widetilde{T}_{k}\nabla
J_{k}(u_{k}+w_{k}(\sum^{n}_{s=1}h_{s})+\sum^{n}_{s=1}h_{s}),\
\varphi)_{k}\nonumber\\
&=&\int_{Q_{k}}\nabla(u_{k}+w_{k}(\sum^{n}_{s=1}h_{s})+\sum^{n}_{s=1}h_{s})\cdot\nabla\varphi+
\int_{Q_{k}}V(x)\cdot(u_{k}+w_{k}(\sum^{n}_{s=1}h_{s})+\sum^{n}_{s=1}h_{s})\cdot\varphi\nonumber\\
&&-\int_{Q_{k}}f(x,u_{k}+w_{k}(\sum^{n}_{s=1}h_{s})+\sum^{n}_{s=1}h_{s})\cdot\varphi.\nonumber
\end{eqnarray}
Differentiating the above equality for the variable $x_{i,j}$, we
get that
\begin{eqnarray}\label{qinqinhebiancao}
0&=&\int_{Q_{k}}\nabla(w'_{k}
(\sum^{n}_{s=1}h_{s})e^{k}_{i,j}+e^{k}_{i,j})\cdot\nabla\varphi+
\int_{Q_{k}}V(x)\cdot(w'_{k}
(\sum^{n}_{s=1}h_{s})e^{k}_{i,j}+e^{k}_{i,j})\cdot\varphi\nonumber\\
&&-\int_{Q_{k}}f'(x,u_{k}+w_{k}(\sum^{n}_{s=1}h_{s})+\sum^{n}_{s=1}h_{s})
\cdot(w'_{k}
(\sum^{n}_{s=1}h_{s})e^{k}_{i,j})\cdot\varphi\nonumber\\
&&-\int_{Q_{k}}f'(x,u_{k}+w_{k}(\sum^{n}_{s=1}h_{s})+\sum^{n}_{s=1}h_{s})\cdot
e^{k}_{i,j}\cdot\varphi.
\end{eqnarray}
Since  $\varphi\in\Pi_{k},$ we have $(e^{k}_{i,j},\varphi)_{k}=0.$
Thus by (\ref{qinqinhebiancao}), we get that
\begin{eqnarray}\label{chunmiambujiuexiao}
&&\int_{Q_{k}}\nabla(w'_{k}
(\sum^{n}_{s=1}h_{s})e^{k}_{i,j})\cdot\nabla\varphi+
\int_{Q_{k}}V(x)\cdot(w'_{k}
(\sum^{n}_{s=1}h_{s})e^{k}_{i,j})\cdot\varphi\nonumber\\
&&-\int_{Q_{k}}f'(x,u_{k}+w_{k}(\sum^{n}_{s=1}h_{s})+\sum^{n}_{s=1}h_{s})
\cdot(w'_{k}
(\sum^{n}_{s=1}h_{s})e^{k}_{i,j})\cdot\varphi\nonumber\\
&=&\int_{Q_{k}}f'(x,u_{k}+w_{k}(\sum^{n}_{s=1}h_{s})+\sum^{n}_{s=1}h_{s})\cdot
e^{k}_{i,j}\cdot\varphi.
\end{eqnarray}
By the same argument and noting that
$\displaystyle(w_{k}^{s})'(\sum^{l_{s}}_{t=1}x_{s,t}e^{k}_{s,t})
e^{k}_{i,j}=0$ if $s\neq i$, we know that for any
$\varphi\in\Pi_{k},$
\begin{eqnarray}\label{dadifeiying}
&&\int_{Q_{k}}\nabla((w_{k}^{s})'(\sum^{l_{s}}_{t=1}x_{s,t}e^{k}_{s,t})
e^{k}_{i,j}))\cdot\nabla\varphi +\int_{Q_{k}}V(x)\cdot
((w_{k}^{s})'(\sum^{l_{s}}_{t=1}x_{s,t}e^{k}_{s,t})
e^{k}_{i,j}) \cdot\varphi\nonumber\\
&&-\int_{Q_{k}}f'(x,u^{s}_{k}+w_{k}^{s}(\sum^{l_{s}}_{t=1}x_{s,t}e^{k}_{s,t})
+\sum^{l_{s}}_{t=1}x_{s,t}e^{k}_{s,t})\cdot((w_{k}^{s})'(\sum^{l_{s}}_{t=1}x_{s,t}e^{k}_{s,t})
e^{k}_{i,j})\cdot\varphi\nonumber\\
&=&\delta_{s,i}\int_{Q_{k}}f'(x,u^{s}_{k}+w_{k}^{s}(\sum^{l_{s}}_{t=1}x_{s,t}e^{k}_{s,t})
+\sum^{l_{s}}_{t=1}x_{s,t}e^{k}_{s,t})\cdot e^{k}_{i,j}\cdot\varphi
\end{eqnarray}
By Lemma \ref{tongrenttthqqqw}, we get that as $k\rightarrow\infty,$
the following two equalities
\begin{eqnarray}\label{vvvacdseeewee}
&&\int_{Q_{k}}f'(x,\sum^{n}_{s=1}u^{s}_{k}+\sum^{n}_{s=1}w_{k}^{s}(\sum^{l_{s}}_{t=1}x_{s,t}e^{k}_{s,t})
+\sum^{l_{s}}_{t=1}x_{s,t}e^{k}_{s,t})\cdot
e^{k}_{i,j}\cdot\varphi\nonumber\\
&=&\sum^{n}_{s=1}\int_{Q_{k}}f'(x,u^{s}_{k}+w_{k}^{s}(\sum^{l_{s}}_{t=1}x_{s,t}e^{k}_{s,t})
+\sum^{l_{s}}_{t=1}x_{s,t}e^{k}_{s,t})\cdot
e^{k}_{i,j}\cdot\varphi+o(1)
\end{eqnarray}
and \begin{eqnarray}\label{chhagsfpppli}
&&\int_{Q_{k}}f'(x,\sum^{n}_{s=1}u^{s}_{k}+\sum^{n}_{s=1}w_{k}^{s}(\sum^{l_{s}}_{t=1}x_{s,t}e^{k}_{s,t})
+\sum^{l_{s}}_{t=1}x_{s,t}e^{k}_{s,t})\cdot
((w_{k}^{s})'(\sum^{l_{s}}_{t=1}x_{s,t}e^{k}_{s,t})e^{k}_{i,j})\cdot\varphi\nonumber\\
&=&\sum^{n}_{s=1}\int_{Q_{k}}f'(x,u^{s}_{k}+w_{k}^{s}(\sum^{l_{s}}_{t=1}x_{s,t}e^{k}_{s,t})
+\sum^{l_{s}}_{t=1}x_{s,t}e^{k}_{s,t})\cdot
((w_{k}^{s})'(\sum^{l_{s}}_{t=1}x_{s,t}e^{k}_{s,t})e^{k}_{i,j})\cdot\varphi\nonumber\\
&&+o(1)
\end{eqnarray}
hold uniformly for $h\in\overline{B_{\Lambda_{k}}(0,\delta_{0})}$
and  $\varphi\in\Pi_{k}$ satisfying $||\varphi||\leq 1$.
Furthermore, if $s\neq i,$ then by (\ref{tapainqinshan}) and the
fact that $|b^{i}_{k}-b^{s}_{k}|\rightarrow\infty$ as
$k\rightarrow\infty,$ we get that the following two limits
\begin{equation}\label{bvvmmmkooi}
\lim_{k\rightarrow\infty}\int_{Q_{k}}f'(x,u^{s}_{k}+w_{k}^{s}(\sum^{l_{s}}_{t=1}x_{s,t}e^{k}_{s,t})
+\sum^{l_{s}}_{t=1}x_{s,t}e^{k}_{s,t})\cdot
((w_{k}^{i})'(\sum^{l_{i}}_{t=1}x_{i,t}e^{k}_{i,t})e^{k}_{i,j})\cdot\varphi=
0\end{equation}   and
\begin{equation}\label{dakhkhltuiu}
\lim_{k\rightarrow\infty}\int_{Q_{k}}f'(x,u^{s}_{k}+w_{k}^{s}(\sum^{l_{s}}_{t=1}x_{s,t}e^{k}_{s,t})
+\sum^{l_{s}}_{t=1}x_{s,t}e^{k}_{s,t})\cdot e^{k}_{i,j}\cdot\varphi=
0\end{equation}
 holds uniformly for
$\sum^{l_{s}}_{t=1}x^{2}_{s,t}\leq\delta^{2}_{0}$. By
$(\ref{vvvacdseeewee})-(\ref{dakhkhltuiu})$,  we get that as
$k\rightarrow\infty,$ the following two equalities
\begin{eqnarray}\label{qianshanniano}
&&\int_{Q_{k}}f'(x,\sum^{n}_{s=1}u^{s}_{k}+\sum^{n}_{s=1}w_{k}^{s}(h_{s})
+\sum^{n}_{s=1}h_{s})\cdot
((w_{k}^{i})'(h_{i})e^{k}_{i,j})\cdot\varphi\nonumber\\
&=&\int_{Q_{k}}f'(x,u^{i}_{k}+w_{k}^{i}(h_{i}) +h_{i})\cdot
((w_{k}^{i})'(h_{i})e^{k}_{i,j})\cdot\varphi+o(1)\end{eqnarray} and
\begin{eqnarray}\label{hongjunbu}
&&\int_{Q_{k}}f'(x,\sum^{n}_{s=1}u^{s}_{k}+\sum^{n}_{s=1}w_{k}^{s}(h_{s})
+\sum^{n}_{s=1}h_{s})\cdot
e^{k}_{i,j}\cdot\varphi\nonumber\\
&=&\int_{Q_{k}}f'(x,u^{i}_{k}+w_{k}^{i}(h_{i}) +h_{i})\cdot
e^{k}_{i,j}\cdot\varphi+o(1)
\end{eqnarray}
hold uniformly for
$\sum^{n}_{s=1}h_{s}\in\overline{B_{\Lambda_{k}}(0,\delta_{0})}$ as
$k\rightarrow\infty.$
 By (\ref{dadifeiying}),
(\ref{qianshanniano}), (\ref{hongjunbu}) and the fact that
$\displaystyle(w_{k}^{s})'(\sum^{l_{s}}_{t=1}x_{s,t}e^{k}_{s,t})
e^{k}_{i,j}=0$ if $s\neq i$, we get that
\begin{eqnarray}\label{chongminrenhenben}
&&\int_{Q_{k}}\nabla(\sum^{n}_{s=1}(w_{k}^{s})'(\sum^{l_{s}}_{t=1}x_{s,t}e^{k}_{s,t})
e^{k}_{i,j})\cdot\nabla\varphi
+\int_{Q_{k}}V(x)\cdot(\sum^{n}_{s=1}(w_{k}^{s})'(\sum^{l_{s}}_{t=1}x_{s,t}e^{k}_{s,t})
e^{k}_{i,j})\cdot\varphi\nonumber\\
&=&\int_{Q_{k}}\nabla((w_{k}^{i})'(\sum^{l_{i}}_{t=1}x_{i,t}e^{k}_{i,t})
e^{k}_{i,j})\cdot\nabla\varphi
+\int_{Q_{k}}V(x)\cdot((w_{k}^{i})'(\sum^{l_{i}}_{t=1}x_{i,t}e^{k}_{i,t})
e^{k}_{i,j})\cdot\varphi\nonumber\\
&=&\int_{Q_{k}}f'(x,u^{i}_{k}+w_{k}^{i}(\sum^{l_{i}}_{t=1}x_{i,t}e^{k}_{i,t})
+\sum^{l_{i}}_{t=1}x_{i,t}e^{k}_{i,t})\cdot
e^{k}_{i,j}\cdot\varphi\nonumber\\
&&+\int_{Q_{k}}f'(x,u^{i}_{k}+w_{k}^{i}(\sum^{l_{i}}_{t=1}x_{i,t}e^{k}_{i,t})
+\sum^{l_{i}}_{t=1}x_{i,t}e^{k}_{i,t})\cdot
((w_{k}^{s})'(\sum^{l_{s}}_{t=1}x_{s,t}e^{k}_{s,t})e^{k}_{i,j})\cdot\varphi.\nonumber\\
&=&\int_{Q_{k}}f'(x,\sum^{n}_{s=1}u^{s}_{k}+\sum^{n}_{s=1}w_{k}^{s}(\sum^{l_{s}}_{t=1}x_{s,t}e^{k}_{s,t})
+\sum^{l_{s}}_{t=1}x_{s,t}e^{k}_{i,t})\cdot
e^{k}_{i,j}\cdot\varphi\nonumber\\
&&+\int_{Q_{k}}f'(x,\sum^{n}_{s=1}u^{s}_{k}+\sum^{n}_{s=1}w_{k}^{s}(\sum^{l_{s}}_{t=1}x_{s,t}e^{k}_{s,t})
+\sum^{l_{s}}_{t=1}x_{s,t}e^{k}_{s,t})\cdot
(\sum^{n}_{s=1}(w_{k}^{s})'(\sum^{l_{s}}_{t=1}x_{s,t}e^{k}_{s,t})e^{k}_{i,j})\cdot\varphi\nonumber\\
&&+o(1).
\end{eqnarray}
Hence by (\ref{chongminrenhenben}), $\displaystyle\lim_{k\rightarrow
\infty}||u_{k}-\sum^{n}_{s=1}u^{s}_{k}||_{k}= 0$ and the result (1)
of this Lemma, we have
\begin{eqnarray}\label{kejinhhhhhbg}
&&\int_{Q_{k}}\nabla(\sum^{n}_{s=1}(w_{k}^{s})'(\sum^{l_{s}}_{t=1}x_{s,t}e^{k}_{s,t})
e^{k}_{i,j})\cdot\nabla\varphi
+\int_{Q_{k}}V(x)\cdot(\sum^{n}_{s=1}(w_{k}^{s})'(\sum^{l_{s}}_{t=1}x_{s,t}e^{k}_{s,t})
e^{k}_{i,j})\cdot\varphi\nonumber\\
&&-\int_{Q_{k}}f'(x,u_{k}+w_{k}(\sum^{n}_{s=1}h_{s})+\sum^{n}_{s=1}h_{s})\cdot
(\sum^{n}_{s=1}(w_{k}^{s})'(h_{s})\cdot\varphi\nonumber\\
&=&\int_{Q_{k}}f'(x,u_{k}+w_{k}(\sum^{n}_{s=1}h_{s})+\sum^{n}_{s=1}h_{s})\cdot
e^{k}_{i,j}\cdot\varphi +o(1).
\end{eqnarray}
holds uniformly for
$\displaystyle\sum^{n}_{s=1}h_{s}\in\overline{B_{\Lambda_{k}}(0,\delta_{0})}$
as $k\rightarrow\infty.$

Minus (\ref{chunmiambujiuexiao}) by (\ref{kejinhhhhhbg}), we get
that for any $\varphi\in\Pi_{k}$,
\begin{eqnarray}\label{qinermemduan}
&&\left((\widetilde{T}_{k}\nabla^{2}J_{k}(u_{k}+w_{k}(\sum^{n}_{s=1}h_{s})+\sum^{n}_{s=1}h_{s}))(w'_{k}
(\sum^{n}_{l=1}h_{l})e^{k}_{i,j}-(\sum^{n}_{s=1}(w^{s}_{k})'
(h_{s})e^{k}_{i,j})),\
\varphi\right)_{k}\nonumber\\
&&=o(1),\ \mbox{as}\ k\rightarrow\infty.
\end{eqnarray}

By (\ref{tapainqinshan}), we deduce that as $k\rightarrow\infty,$
for any $1\leq \nu\leq l_{t}$ and $ 1\leq t\leq n$, we have
\begin{equation}\label{dudiaohanjianxue}
(\sum^{n}_{s=1}(w^{s}_{k})' (h_{s})e^{k}_{i,j},\
e^{k}_{t,\nu})_{k}\rightarrow 0
\end{equation}
holds uniformly for $h\in \overline{B_{\Lambda_{k}}(0,\delta_{0})}$.
Thus  as $k\rightarrow\infty,$ the limit
\begin{equation}\label{hongchengyouuuujhdfer}
||\sum^{n}_{s=1}(w^{s}_{k})'
(h_{s})e^{k}_{i,j}-\widetilde{T}_{k}(\sum^{n}_{s=1}(w^{s}_{k})'
(h_{s})e^{k}_{i,j})||_{k}\rightarrow 0
\end{equation}
holds uniformly for $h\in \overline{B_{\Lambda_{k}}(0,\delta_{0})}$.
By (\ref{hongchengyouuuujhdfer}) and (\ref{qinermemduan}), we get
that
\begin{eqnarray}\label{qinermemduantt5r}
&&\left((\widetilde{T}_{k}\nabla^{2}J_{k}(u_{k}+w_{k}(\sum^{n}_{s=1}h_{s})+\sum^{n}_{s=1}h_{s}))(w'_{k}
(\sum^{n}_{s=1}h_{s})e^{k}_{i,j}-(\sum^{n}_{s=1}\widetilde{T}_{k}(w^{s}_{k})'
(h_{s})e^{k}_{i,j})),\
\varphi\right)_{k}\nonumber\\
&&=o(1),\ \mbox{as}\ k\rightarrow\infty.
\end{eqnarray}
Choose $\varphi=(w'_{k}
(\sum^{n}_{s=1}h_{s})e^{k}_{i,j}-(\sum^{n}_{s=1}\widetilde{T}_{k}(w^{s}_{k})'
(h_{s})e^{k}_{i,j}))$ in (\ref{qinermemduantt5r}) and by  Lemma
\ref{xiaojiazhi}, we deduce that as $k\rightarrow\infty,$
\begin{equation}\label{xianqianchongkkku}
||(w'_{k}
(\sum^{n}_{l=1}h_{l})e^{k}_{i,j})-\widetilde{T}_{k}(\sum^{n}_{s=1}(w^{s}_{k})'
(h_{s})e^{k}_{i,j})||_{k}\rightarrow 0
\end{equation}
holds uniformly for $h\in \overline{B_{\Lambda_{k}}(0,\delta_{0})}$.

 Thus the result {\bf (2)} of this Lemma follows from
(\ref{xianqianchongkkku}) and (\ref{hongchengyouuuujhdfer})
directly. \hfill$\Box$

\end{document}